\documentclass[11pt]{article}
\usepackage{amsmath,amssymb,amsthm,graphics,environ}
\newsavebox{\measuretikzpicture}
\NewEnviron{scaletikzpicturetowidth}[1]{%
  \def\tikzscale{0.05}
  \savebox{\measuretikzpicture}{\BODY}%
  \pgfmathparse{#1/\wd\measuretikzpicture*\tikzscale}%
  \edef\tikzscale{\pgfmathresult}%
  \BODY
}
\usepackage{rotating,pdflscape,afterpage}
\usepackage[sc,osf]{mathpazo}
\AtBeginDocument{
	\DeclareSymbolFont{AMSb}{U}{msb}{m}{n}
	\DeclareSymbolFontAlphabet{\mathbb}{AMSb}
}
\usepackage{graphicx}
\usepackage{extpfeil}
\usepackage{mathabx} 
\usepackage{float}
\usepackage{caption}
\usepackage{subcaption}
\newcommand{\mockalph}[1]{\!}
\usepackage{enumitem}
\usepackage[nouppercase]{scrlayer-scrpage}
\usepackage{avant} 
\usepackage{mathrsfs}
\usepackage{mathtools}
\usepackage{faktor,xfrac}
\usepackage[all]{xypic}\xyoption{rotate}
\newdir{ >}{{}*!/-7pt/\dir{>}}
\newdir{i}{{}*!/-7pt/\dir{(}}
\usepackage{thmtools}
\usepackage[utf8]{inputenc}
\usepackage[T1]{fontenc}
\usepackage{tocloft}
\makeatletter
\renewcommand{\l@figure}{\@dottedtocline{1}{1em}{3.5em}}
\renewcommand{\l@table}{\@dottedtocline{2}{1em}{3.5em}}
\makeatother
\newcommand*{\noaddvspace}{\renewcommand*{\addvspace}[1]{}}
\addtocontents{lof}{\protect\noaddvspace}

\usepackage[hyphens]{url}		%
\usepackage{hyperref}	
\usepackage{cleveref} 	

\makeatletter
\let\c@figure\c@table
\let\c@equation\c@table
\makeatother

\numberwithin{table}{subsection}
\numberwithin{equation}{subsection}
\numberwithin{figure}{subsection}

\newtheorem{theorem}[table]{Theorem}

\newtheorem{proposition}[table]{Proposition}

\newtheorem{corollary}[table]{Corollary}

\newtheorem{lemma}[table]{Lemma}

\newtheorem{claim}[table]{Claim}

\theoremstyle{definition}

\newtheorem{definition}[table]{Definition}

\newtheorem{construction}[table]{Construction}
\newtheorem{notation}[table]{Notation}

\newtheorem{observation}[table]{Observation}

\newtheorem{conjecture}[table]{Conjecture}

\newtheorem{discussion}[table]{}

\theoremstyle{remark}
\newtheorem{fact}[table]{Fact}
\newtheorem{examples}[table]{Examples}
\newtheorem{example}[table]{Example}
\newtheorem{exercise}[table]{Exercise}

\newtheorem{problem}[table]{Problem}
\newtheorem{histrmks}[table]{Historical remarks}
\newtheorem{remark}[table]{Remark}
\newtheorem{remarks}[table]{Remarks}

\theoremstyle{plain}
\newtheorem*{thm*}{Theorem}
\newtheorem*{theorem*}{Theorem}
\newtheorem*{prop*}{Proposition}
\newtheorem*{proposition*}{Proposition}
\newtheorem*{lemma*}{Lemma}
\newtheorem*{corollary*}{Corollary}
\newtheorem*{cor*}{Corollary}

\theoremstyle{definition}
\newtheorem*{definition*}{Definition}
\newtheorem*{defn*}{Definition}
\newtheorem*{QQ*}{Question}
\newtheorem*{obs*}{Observation}
\newtheorem*{notation*}{Notation}

\theoremstyle{remark}
\newtheorem*{rmk*}{Remark}
\newtheorem*{remark*}{Remark}

\newtheorem*{examples*}{Examples}
\newtheorem*{example*}{Example}
\newtheorem*{EG*}{Example}
\newtheorem*{EGs*}{Examples}
\newtheorem*{fact*}{Fact}
\newtheorem*{prob*}{Problem}

\usepackage{etoolbox}
\usepackage[usenames,dvipsnames,svgnames,table]{xcolor}
\newcommand		{\defd}[1]	{\textcolor{RoyalBlue}{\textbf{\textit{#1}}}}
\newcommand		{\defm}[1]	{\textcolor{RoyalBlue}{#1}}

\tracingpatches
\makeatletter
\patchcmd{\@setref}{\bfseries ??}{\bfseries\color{red} FIX ME!}{}{}
\patchcmd{\@setcite}{\bfseries ?}{\bfseries\color{red} FIX ME!}{}{}
\patchcmd{\@setcref}         {??}{\color{red} FIX ME!}{}{}
\patchcmd{\@setcref}         {??}{\color{red} FIX ME!}{}{}
\patchcmd{\@setcrefrange}    {??}{\color{red} FIX ME!}{}{}
\patchcmd{\@setcrefrange}    {??}{\color{red} FIX ME!}{}{}
\patchcmd{\@setcrefrange}    {??}{\color{red} FIX ME!}{}{}
\patchcmd{\@setcrefrange}    {??}{\color{red} FIX ME!}{}{}
\patchcmd{\@setcrefrange}    {??}{\color{red} FIX ME!}{}{}
\patchcmd{\@setcrefrange}    {??}{\color{red} FIX ME!}{}{}
\patchcmd{\@setnamecref}     {??}{\color{red} FIX ME!}{}{}
\patchcmd{\@setnamecref}     {??}{\color{red} FIX ME!}{}{}
\patchcmd{\@setcpageref}     {??}{\color{red} FIX ME!}{}{}
\patchcmd{\@setcpageref}     {??}{\color{red} FIX ME!}{}{}
\patchcmd{\@setcpagerefrange}{??}{\color{red} FIX ME!}{}{}
\patchcmd{\@setcpagerefrange}{??}{\color{red} FIX ME!}{}{}
\patchcmd{\@setcpagerefrange}{??}{\color{red} FIX ME!}{}{}
\patchcmd{\@setcpagerefrange}{??}{\color{red} FIX ME!}{}{}
\patchcmd{\@setcpagerefrange}{??}{\color{red} FIX ME!}{}{}
\patchcmd{\@cref}            {??}{\color{red} FIX ME!}{}{}
\def\blx@citation@entry#1#2{%
	\blx@bibreq{#1}%
	\ifinlist{#1}{\blx@cites}
	{}
	{\listgadd{\blx@cites}{#1}%
		\blx@auxwrite\@mainaux{}{\string\abx@aux@cite{#1}}}%
	\ifinlistcs{#1}{blx@segm@\the\c@refsection @\the\c@refsegment}
	{}
	{\listcsgadd{blx@segm@\the\c@refsection @\the\c@refsegment}{#1}}%
	\blx@ifdata{#1}%
	{}%
	{\ifcsdef{blx@miss@\the\c@refsection}%
		{\ifinlistcs{#1}{blx@miss@\the\c@refsection}%
			{{\bfseries\color{red} cite:} }%
			{\blx@logreq@active{#2{#1}}}}%
		{\blx@logreq@active{#2{#1}}}}}
\def\blx@citeadd#1{%
	\ifcsdef{blx@keyalias@\the\c@refsection @#1}
	{\edef\blx@realkey{\csuse{blx@keyalias@\the\c@refsection @#1}}}
	{\def\blx@realkey{#1}}%
	\expandafter\blx@citation\expandafter{\blx@realkey}\blx@msg@cundefon
	\expandafter\blx@ifdata\expandafter{\blx@realkey}
	{\advance\blx@tempcnta\@ne
		\listeadd\blx@tempa{\blx@realkey}}
	{\ifnum\blx@tempcntb>\z@MUlticitedelim\fi
		\expandafter\abx@missing\expandafter{\blx@realkey}%
		\advance\blx@tempcntb\@ne}}

\DeclarePairedDelimiterX{\pmodx}[1]{(}{)}{{\operator@font mod}\mkern6mu#1}
\renewcommand{\pmod}{%
	\allowbreak
	\if@display\mkern18mu\else\mkern8mu\fi
	\pmodx
}
\makeatother

\makeatletter
\newcommand{\oset}[3][0ex]{%
	\raisebox{.175ex}{$%
		\mathrel{\mathop{#3}\limits^{
				\vbox to#1{\kern-2\ex@
					\hbox{$\scriptstyle#2$}\vss}}}
		$}%
}
\makeatother

\newcommand{\myred}{BrickRed}
\hypersetup{
	colorlinks			= true, 	
	urlcolor			= \myred, 
	linkcolor			= \myred, 	
	citecolor			= PineGreen 	
}

\usepackage{xspace}
\usepackage[left=2.54cm,right=2.54cm,top=2.54cm,bottom=2.54cm]{geometry}
\usepackage{tikz}
\usetikzlibrary{matrix,fit,backgrounds,calc,arrows,shapes} 
\tikzstyle{image}=[rectangle,fill=Red!20,inner sep=-2pt]
\tikzstyle{nonzero}=[rectangle,fill=Navy!20,inner sep=0pt]
\tikzstyle{nonzerosm}=[rectangle,fill=Navy!20,inner sep=-2pt]

\newcommand{\dgeq}{\rotatebox[origin=c]{-45}{$\geq$}}
\newcommand{\ugeq}{\rotatebox[origin=c]{+45}{$\geq$}}

\makeatletter
\newbox\xrat@below
\newbox\xrat@above
\newcommand{\xrightarrowtail}[2][]{%
	\setbox\xrat@below=\hbox{\ensuremath{\scriptstyle #1}}%
	\setbox\xrat@above=\hbox{\ensuremath{\scriptstyle #2}}%
	\pgfmathsetlengthmacro{\xrat@len}{max(\wd\xrat@below,\wd\xrat@above)+.6em}%
	\mathrel{\tikz [>->,baseline=-.55ex]
		\draw (0,0) -- node[below=-2pt] {\box\xrat@below}
		node[above=-2pt] {\box\xrat@above}
		(\xrat@len,0) ;}}
\newbox\xrat@below
\newbox\xrat@above
\renewcommand{\xtwoheadrightarrow}[2][]{%
	\setbox\xrat@below=\hbox{\ensuremath{\scriptstyle #1}}%
	\setbox\xrat@above=\hbox{\ensuremath{\scriptstyle #2}}%
	\pgfmathsetlengthmacro{\xrat@len}{max(\wd\xrat@below,\wd\xrat@above)+.6em}%
	\mathrel{\tikz [->>,baseline=-.55ex]
		\draw (0,0) -- node[below=-2pt] {\box\xrat@below}
		node[above=-2pt] {\box\xrat@above}
		(\xrat@len,0) ;}}
\makeatother
\newcommand{\xmono}{\xrightarrowtail}
\newcommand{\mono}{\xmono{\phantom{\ \, }}}
\newcommand{\xepi}{\xtwoheadrightarrow}
\newcommand{\epi}{\xepi{\phantom{\ \, }}}

\makeatletter
\newcommand{\presectionskip}{-1.5\baselineskip}
\newcommand{\postsectionskip}{0.3\baselineskip}
\usepackage{titlesec}
\renewcommand{\section}{\@startsection
	{chapter}{0}{0mm}
	{\presectionskip}
	{\postsectionskip}
	{\sffamily\huge}}
\renewcommand{\section}{\@startsection
	{section}{1}{0mm}
	{\presectionskip}
	{\postsectionskip}
	{\sffamily\LARGE}}
\renewcommand{\subsection}{\@startsection
	{subsection}{2}{0mm}
	{\presectionskip}
	{\postsectionskip}
	{\sffamily\Large}}
\renewcommand{\subsubsection}{\@startsection
	{subsubsection}{3}{0mm}
	{\presectionskip}
	{\postsectionskip}
	{\sffamily\normalsize}}
\renewcommand{\@seccntformat}[1]{\csname the#1\endcsname.\quad}
\newcommand\HUGE{\@setfontsize\Huge{30}{47}} 
\titleformat{\chapter}[display]
{\sffamily\Large}
{Chapter {\HUGE\normalfont\thechapter}}    
{1em}
{\huge}
\titleformat{\section}
{\sffamily\Large}
{\normalfont{\@setfontsize\Large{18}{47}\thesection.}}    
{1em}
{}
\titleformat{\subsection}
{\sffamily\large}
{\normalfont{\Large\thesubsection}.}    
{1em}
{}
\titleformat{\subsubsection}
{\sffamily\normalsize}
{\normalfont{\large\thesubsubsection}.}    
{1em}
{}
\makeatother

\makeatletter
\def\smallunderbrace#1{\mathop{\vtop{\m@th\ialign{##\crcr
				$\hfil\displaystyle{#1}\hfil$\crcr
				\noalign{\kern3\p@\nointerlineskip}%
				\tiny\upbracefill\crcr\noalign{\kern3\p@}}}}\limits}
\makeatother

\renewcommand{\SS}{\textsection}
\newcommand{\bthm}{\begin{theorem}}
	\newcommand{\ethm}{\end{theorem}}
\newcommand{\bprop}{\begin{proposition}}
	\newcommand{\eprop}{\end{proposition}}
\newcommand{\bcor}{\begin{corollary}}
	\newcommand{\ecor}{\end{corollary}}
\newcommand{\bconj}{\begin{conjecture}}
	\newcommand{\econj}{\end{conjecture}}
\newcommand{\blem}{\begin{lemma}}
	\newcommand{\elem}{\end{lemma}}
\newcommand{\bclm}{\begin{claim}}
	\newcommand{\eclm}{\end{claim}}
\newcommand{\bpf}{\begin{proof}}
	\newcommand{\epf}{\end{proof}}
\newcommand{\bdetails}{\begin{details}}
	\newcommand{\edetails}{\end{details}}
\newcommand{\bdefi}{\begin{definition}}
	\newcommand{\edefi}{\end{definition}}
\newcommand{\bdefn}{\begin{definition}}
	\newcommand{\edefn}{\end{definition}}
\newcommand{\bex}{\begin{example}}
	\newcommand{\eex}{\end{example}}
\newcommand{\bprob}{\begin{problem}}
	\newcommand{\eprob}{\end{problem}}
\newcommand{\bob}{\begin{observation}}
	\newcommand{\eob}{\end{observation}}
\newcommand{\bexer}{\begin{exercise}}
	\newcommand{\eexer}{\end{exercise}}
\newcommand{\bexers}{\begin{exercises}}
	\newcommand{\eexers}{\end{exercises}}
\newcommand{\brmk}{\begin{remark}}
	\newcommand{\ermk}{\end{remark}}
\newcommand{\bhist}{\begin{histrmks}}
	\newcommand{\ehist}{\end{histrmks}}
\newcommand{\brmks}{\begin{remarks}}
	\newcommand{\ermks}{\end{remarks}}
\newcommand{\bverb}{\begin{verbatim}}
	\newcommand{\everb}{\end{verbatim}}
\newcommand{\bntn}{\begin{notation}}
	\newcommand{\entn}{\end{notation}}
\newcommand{\bfct}{\begin{fact}}
	\newcommand{\efct}{\end{fact}}
\newcommand{\bfcts}{\begin{facts}}
	\newcommand{\befcts}{\end{facts}}
\newcommand{\benum}{\begin{enumerate}}
	\newcommand{\eenum}{\end{enumerate}}
\newcommand{\bitem}{\begin{itemize}}
	\newcommand{\eitem}{\end{itemize}}

\renewcommand	{\o}		{\circ}

\renewcommand	{\:}		{\colon}

\newcommand		{\quotientmed}[2]	{{\raisebox{.2em}{$#1$}}\ \!\!\big/\!\!\ 
	{\raisebox{-.2em}{$#2$}}}

\newcommand		{\SSS}	{Serre spectral sequence\xspace}

\newcommand		{\exterior}	{\Lambda}
\newcommand		{\ext}		{\exterior}

\renewcommand	{\th}		{^{\mathrm{th}}}
\newcommand		{\st}		{^{\mathrm{st}}}

\makeatletter
\newcommand{\subalign}[1]{%
	\vcenter{%
		\Let@ \restore@math@cr \default@tag
		\baselineskip\fontdimen10 \scriptfont\tw@
		\advance\baselineskip\fontdimen12 \scriptfont\tw@
		\lineskip\thr@@\fontdimen8 \scriptfont\thr@@
		\lineskiplimit\lineskip
		\ialign{\hfil$\m@th\scriptstyle##$&$\m@th\scriptstyle{}##$\crcr
			#1\crcr
		}%
	}
}
\makeatother

\newcommand		{\eqn}[1]			{\begin{align*} #1 \end{align*}}

\newcommand		{\quation}[1]		{\begin{equation} #1 \end{equation}}

\newcommand		{\case}[1]			{\begin{cases} #1 \end{cases}}

\newcommand		{\bs}				{\bigskip}

\newcommand		{\mn}				{\mspace{-2mu}}
\newcommand		{\mnn}				{\mspace{-1mu}}
\newcommand		{\dsp}				{\displaystyle}

\newcommand		{\nd}			{\noindent}
\newcommand		{\ol}			{\overline}
\newcommand		{\os}			{\overset}
\newcommand		{\us}			{\underset}

\newcommand		{\wt}			{\widetilde}

\newcommand		{\mr}		{\mathrm}
\newcommand		{\bb}		{\mathbb}

\newcommand		{\mc}		{\mathcal}
\newcommand		{\f}		{\mathfrak}

\renewcommand	{\a}		{\alpha}
\renewcommand	{\b}		{\beta}
\renewcommand	{\d}		{\delta}
\newcommand		{\g}		{\gamma}
\renewcommand	{\epsilon}	{\varepsilon}
\newcommand		{\vk}		{\varkappa}
\newcommand		{\e}		{\epsilon}
\newcommand		{\z}		{\zeta}

\newcommand		{\vt}		{\vartheta}
\renewcommand	{\l}		{\lambda}

\newcommand		{\s}		{\sigma}

\newcommand		{\vp}		{\varphi}
\newcommand		{\w}		{\omega}

\newcommand		{\W}		{\Omega}

\newcommand		{\D}		{\Delta}

\DeclareSymbolFont{cmletters}{OT1}{cmr}{m}{n}
\DeclareMathSymbol{\Ups}{\mathalpha}{cmletters}{"7}
\renewcommand	{\Upsilon}{\Ups}

\newcommand		{\ft}		{\f t}
\newcommand		{\fg}		{\f g}

\newcommand		{\F}		{\bb F}

\newcommand		{\Z}		{\bb Z}
\newcommand		{\Q}		{\bb Q}
\newcommand		{\R}		{\bb R}
\newcommand		{\C}		{\bb C}

\renewcommand 	{\H}	{H^*}

\let\union\cup%
\renewcommand	{\cup}		{\mspace{-1mu}\smile\mspace{-1mu}}
\let\inter\cap%
\newcommand		{\Union}	{\bigcup}

\newcommand		{\less}		{\setminus}
\newcommand		{\sub}		{\subseteq}
\newcommand		{\subn}		{\subsetneq}

\DeclareRobustCommand{\lq}	{\text{\reflectbox{$/$}}}		

\DeclareMathOperator{\id}	{id}

\renewcommand	{\-}		{^{-1}}
\renewcommand	{\o}		{\circ}
\renewcommand	{\.}		{\cdot}
\newcommand		{\x}		{\times}
\newcommand		{\xu}[3]	{\smash{{#2}\us{#1}\times{#3}}}

\DeclareMathOperator*{\otimesvariable}{%
	\mathchoice {\raisebox{.85pt}{$\displaystyle\otimes$}}
	{\raisebox{.85pt}{$\otimes$}}
	{\raisebox{0.7pt}{$\scriptstyle\otimes$}}
	{\raisebox{0.2pt}{$\scriptscriptstyle\otimes$}}
}
\newcommand		{\tensor}		{\otimesvariable}
\newcommand		{\ox}			{\tensor}

\newcommand		{\direct}		{\oplus}

\newcommand		{\+}			{\direct}
\newcommand		{\Tensor}		{\bigotimes}
\newcommand		{\Direct}		{\bigoplus}

\newcommand		 {\diag}		{\operatorname{diag}}

\DeclareMathOperator{\tr}		{tr }
\DeclareMathOperator{\rk}		{rk }
\DeclareMathOperator{\im}		{im }

\DeclareMathOperator{\Tor}		{Tor}
\DeclareMathOperator{\Stab}		{Stab }

\DeclareMathOperator{\Ad}		{Ad }

\DeclareMathOperator{\Aut}		{Aut }

\renewcommand	{\O}			{\mr{O}}

\newcommand		{\SO}			{\mr{SO}}

\newcommand		{\U}			{\mr{U}}
\newcommand		{\SU}			{\mr{SU}}
\newcommand		{\Sp}			{\mr{Sp}}

\newcommand		{\longto} 		{\longrightarrow}
\newcommand		{\lt}			{\longto}
\newcommand		{\xtoo}			{\xrightarrow} 
\newcommand		{\lmt}			{\longmapsto}

\newcommand		{\from}			{\leftarrow}
\newcommand		{\longfrom}		{\longleftarrow}
\newcommand		{\lf}			{\longfrom}

\newcommand		{\inc}		{\hookrightarrow}
\newcommand		{\xinc}		{\xhookrightarrow}
\newcommand		{\longinc}		{\xinc[]{\ \ \ \ }}

\newcommand		{\longepi}	{\xepi[]{\ \ \ \ }}

\newcommand		{\eqto}		{\xrightarrow{=}}
\newcommand		{\simto}		{\xrightarrow{\sim}}

\newcommand		{\longsimto}	{\os\sim\longto}
\newcommand		{\isoto}		{\longsimto}

\newcommand		{\homeoto}		{\os\homeo\to}

\newcommand		{\ceq}			{\coloneqq}
\newcommand		{\eqc}			{\eqqcolon}

\newcommand		{\hmt}			{\simeq}
\newcommand		{\iso}			{\cong}
\newcommand		{\homeo}		{\approx}

\usepackage{todonotes}
\usepackage{csquotes}
\renewcommand{\mkbegdispquote}[2]{\itshape}

\newcommand{\Det}{\mathrm{Det}}

\newcommand{\Orb}{\mathcal{O}}

\newcommand{\GZ}{Gelfand--Zeitlin\xspace}
\newcommand{\TGZS}{topological Gelfand--Zeitlin system\xspace}
\newcommand{\Zf}{\mathbb Z[\sfrac 1 {\,2}]}
\newcommand{\lup}[1]{\l^{\smash{(#1)}}}

\newcommand{\xip}[1]{\xi^{\smash{(#1)}}}
\newcommand{\xipko}{\xip{k+1}}
\newcommand{\genpt}{\vt}
\newcommand{\genup}[1]{\genpt^{\smash{(#1)}}}

\newcommand{\theinterval}{\g}

 \newcommand{\oH}{\ol{H}{}^*}

\newcommand{\Sq}{\mathrm{Sq}}

\newcommand\parallelogram[1][2]{%
\raisebox{-.5pt}{\tikz{\draw[scale = .08,line width=1pt] (0,0)--(3,0)--(4,3)--(1,3)--(0,0)--(3,0);}}}

\newcommand\longpara[1][2]{%
\raisebox{-.5pt}{\tikz{\draw[scale = .08,line width=1pt] (0,0)--(6,0)--(7,3)--(1,3)--(0,0)--(6,0);}}}

\newcommand\rpara{\reflectbox{\parallelogram}}
\newcommand\longrpara{\reflectbox{\longpara}}
 
\newcommand\mtrap[1][2]{%
\raisebox{-0.5pt}{\tikz{\draw[scale = .08,line width=1pt] (0,0)--(4,0)--(3,3)--(1,3)--(0,0)--(4,0);}}}
 
\newcommand\wtrap[1][2]{%
\raisebox{-0.5pt}{\tikz{\draw[scale = .08,line width=1pt] (0,3)--(4,3)--(3,0)--(1,0)--(0,3)--(4,3);}}}

\newcommand\MMtrap[1][2]{%
\raisebox{-.25pt}{\tikz{\draw[scale = .08,line width=1pt] (0,0)--(4,0)--(3,3)--(1,3)--(0,0);
	\draw[scale = .08,line width=1pt] (.5,1.5)--(3.5,1.5);}}%
}

\newcommand\WWtrap[1][2]{%
\raisebox{-.25pt}{\tikz{\draw[scale = .08,line width=1pt] (1,0)--(3,0)--(4,3)--(0,3)--(1,0);
	\draw[scale = .08,line width=1pt] (.5,1.5)--(3.5,1.5);}}%
}

\newcommand\WMtrap[1][2]{%
\raisebox{-.25pt}{\tikz{\draw[scale = .08,line width=1pt] (0.5,0)--(3.5,0)--(4,1.5)--(3.5,3)--(0.5,3)--(0,1.5)--(0.5,0);
\draw[scale = .08,line width=1pt] (0,1.5)--(4,1.5);}}%
}

\newcommand\MWtrap[1][2]{%
\raisebox{-.25pt}{\tikz{\draw[scale = .08,line width=1pt] (0,0)--(4,0)--(3.5,1.5)--(0.5,1.5)--(0,0);
\draw[scale = .08,line width=1pt] (0.5,1.5)--(0,3)--(4,3)--(3.5,1.5);}}%
}

\newcommand\hex[1][2]{%
	\raisebox{-.5pt}{\tikz{\draw[scale = .08,line width=1pt] (0.75,-.25)--(3.25,-.25)--(4,1.5)--(3.25,3.25)--(0.75,3.25)--(0,1.5)--(0.75,-0.25);}}%
}

\begin{document}

\title{\vspace{-1em}\huge The topology of Gelfand--Zeitlin fibers}
\author{Jeffrey D.~Carlson and Jeremy Lane}
		
	\maketitle
	
\begin{abstract}
 We prove several new results about the topology of
 fibers of Gelfand--Zeitlin systems on unitary and orthogonal coadjoint orbits,
 at the same time finding a unifying framework recovering 
 and shedding light on essentially all known results.
 We find completely explicit descriptions of the 
 diffeomorphism type of the fiber in many instances
 a direct factor decomposition of the fiber, 
and a torus factor
corresponding to the action given by the Thimm trick.
The new description also gives us a weak local normal form
for a coadjoint orbit,
which we use to define a topological toric degeneration,
new in the orthogonal case,

 We also compute the first three homotopy groups 
(new in the orthogonal case)
 and cohomology rings of a fiber (new in both cases).
 All these descriptions
 can be read in a
 straightforward manner from the combinatorics of the associated
 Gelfand--Zeitlin pattern. 
\end{abstract}

\section{Introduction}\label{sec:intro}

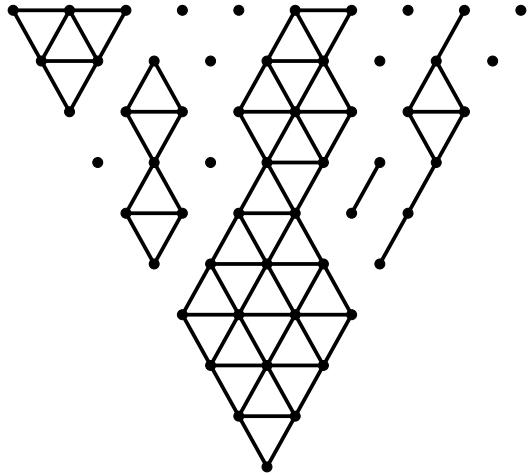
\begin{figure}
	\centering
	\begin{tikzpicture}[scale =.5,line cap=round,line join=round,>=triangle 45,x=1cm,y=1cm,scale=1.5, every node/.style={scale=1.5}]
	\begin{scriptsize}
	
	\draw [fill=black] (1+1/2,9*.9) circle (2.5pt);
	\draw [fill=black] (2+1/2,9*.9) circle (2.5pt);
	\draw [fill=black] (3+1/2,9*.9) circle (2.5pt);
	\draw [fill=black] (4+1/2,9*.9) circle (2.5pt);
	\draw [fill=black] (5+1/2,9*.9) circle (2.5pt);
	\draw [fill=black] (6+1/2,9*.9) circle (2.5pt);
	\draw [fill=black] (7+1/2,9*.9) circle (2.5pt);
	\draw [fill=black] (8+1/2,9*.9) circle (2.5pt);
	\draw [fill=black] (9+1/2,9*.9) circle (2.5pt);
	\draw [fill=black] (10+1/2,9*.9) circle (2.5pt);
	
	\draw [fill=black] (1+2/2,8*.9) circle (2.5pt);
	\draw [fill=black] (2+2/2,8*.9) circle (2.5pt);
	\draw [fill=black] (3+2/2,8*.9) circle (2.5pt);
	\draw [fill=black] (4+2/2,8*.9) circle (2.5pt);
	\draw [fill=black] (5+2/2,8*.9) circle (2.5pt);
	\draw [fill=black] (6+2/2,8*.9) circle (2.5pt);
	\draw [fill=black] (7+2/2,8*.9) circle (2.5pt);
	\draw [fill=black] (8+2/2,8*.9) circle (2.5pt);
	\draw [fill=black] (9+2/2,8*.9) circle (2.5pt);
	
	\draw [fill=black] (1+3/2,7*.9) circle (2.5pt);
	\draw [fill=black] (2+3/2,7*.9) circle (2.5pt);
	\draw [fill=black] (3+3/2,7*.9) circle (2.5pt);
	\draw [fill=black] (4+3/2,7*.9) circle (2.5pt);
	\draw [fill=black] (5+3/2,7*.9) circle (2.5pt);
	\draw [fill=black] (6+3/2,7*.9) circle (2.5pt);
	\draw [fill=black] (7+3/2,7*.9) circle (2.5pt);
	\draw [fill=black] (8+3/2,7*.9) circle (2.5pt);
	
	\draw [fill=black] (1+4/2,6*.9) circle (2.5pt);
	\draw [fill=black] (2+4/2,6*.9) circle (2.5pt);
	\draw [fill=black] (3+4/2,6*.9) circle (2.5pt);
	\draw [fill=black] (4+4/2,6*.9) circle (2.5pt);
	\draw [fill=black] (5+4/2,6*.9) circle (2.5pt);
	\draw [fill=black] (6+4/2,6*.9) circle (2.5pt);
	\draw [fill=black] (7+4/2,6*.9) circle (2.5pt);
	
	\draw [fill=black] (1+5/2,5*.9) circle (2.5pt);
	\draw [fill=black] (2+5/2,5*.9) circle (2.5pt);
	\draw [fill=black] (3+5/2,5*.9) circle (2.5pt);
	\draw [fill=black] (4+5/2,5*.9) circle (2.5pt);
	\draw [fill=black] (5+5/2,5*.9) circle (2.5pt);
	\draw [fill=black] (6+5/2,5*.9) circle (2.5pt);
	
	\draw [fill=black] (1+6/2,4*.9) circle (2.5pt);
	\draw [fill=black] (2+6/2,4*.9) circle (2.5pt);
	\draw [fill=black] (3+6/2,4*.9) circle (2.5pt);
	\draw [fill=black] (4+6/2,4*.9) circle (2.5pt);
	\draw [fill=black] (5+6/2,4*.9) circle (2.5pt);
	
	\draw [fill=black] (1+7/2,3*.9) circle (2.5pt);
	\draw [fill=black] (2+7/2,3*.9) circle (2.5pt);
	\draw [fill=black] (3+7/2,3*.9) circle (2.5pt);
	\draw [fill=black] (4+7/2,3*.9) circle (2.5pt);
	
	\draw [fill=black] (1+8/2,2*.9) circle (2.5pt);
	\draw [fill=black] (2+8/2,2*.9) circle (2.5pt);
	\draw [fill=black] (3+8/2,2*.9) circle (2.5pt);
	
	\draw [fill=black] (1+9/2,1*.9) circle (2.5pt);
	\draw [fill=black] (2+9/2,1*.9) circle (2.5pt);
	
	\draw [fill=black] (1+10/2,0*.9) circle (2.5pt);
	
	\draw [line width=0.5mm] (1+1/2,9*.9) -- (2+1/2,9*.9) -- (3+1/2,9*.9);
	\draw [line width=0.5mm] (1+1/2,9*.9) -- (1+2/2,8*.9) -- (2+1/2,9*.9) -- (2+2/2,8*.9)-- (3+1/2,9*.9);
	\draw [line width=0.5mm] (1+2/2,8*.9) -- (2+2/2,8*.9) -- (1+3/2,7*.9) -- (1+2/2,8*.9);
	
	\draw [line width=0.5mm] (3+2/2,8*.9) -- (2+3/2,7*.9) -- (3+3/2,7*.9) -- (3+2/2,8*.9);
	\draw [line width=0.5mm] (2+3/2,7*.9) -- (2+4/2,6*.9) -- (3+3/2,7*.9);
	\draw [line width=0.5mm] (2+4/2,6*.9) -- (1+5/2,5*.9) -- (2+5/2,5*.9) -- (2+4/2,6*.9);
	\draw [line width=0.5mm] (1+5/2,5*.9) -- (1+6/2,4*.9) -- (2+5/2,5*.9);
	
	\draw [line width=0.5mm] (6+1/2,9*.9) -- (7+1/2,9*.9) -- (6+2/2,8*.9) -- (6+3/2,7*.9)--(5+4/2,6*.9)--(4+5/2,5*.9)--(4+6/2,4*.9)--(4+7/2,3*.9)--(3+8/2,2*.9)--(2+9/2,1*.9)--(1+10/2,0*.9);
	\draw [line width=0.5mm] (1+10/2,0*.9)--(1+9/2,1*.9)--(1+8/2,2*.9)--(1+7/2,3*.9)--(2+6/2,4*.9)--(3+5/2,5*.9)--(4+4/2,6*.9)--(4+3/2,7*.9)--(5+2/2,8*.9)--(6+1/2,9*.9);
	\draw [line width=0.5mm] (5+2/2,8*.9) --(6+2/2,8*.9);
	\draw [line width=0.5mm] (4+3/2,7*.9)--(5+3/2,7*.9)--(6+3/2,7*.9);
	\draw [line width=0.5mm] (4+4/2,6*.9)--(5+4/2,6*.9);
	\draw [line width=0.5mm] (3+5/2,5*.9)--(4+5/2,5*.9);
	\draw [line width=0.5mm] (2+6/2,4*.9)--(3+6/2,4*.9)--(4+6/2,4*.9);
	\draw [line width=0.5mm] (1+7/2,3*.9)--(2+7/2,3*.9)--(3+7/2,3*.9)--(4+7/2,3*.9);
	\draw [line width=0.5mm] (1+8/2,2*.9)-- (2+8/2,2*.9)-- (3+8/2,2*.9);
	\draw [line width=0.5mm] (1+9/2,1*.9)--(2+9/2,1*.9);
	\draw [line width=0.5mm] (6+1/2,9*.9)--(6+2/2,8*.9);
	\draw [line width=0.5mm] (5+2/2,8*.9)--(5+3/2,7*.9)--(5+4/2,6*.9);
	\draw [line width=0.5mm] (4+4/2,6*.9)--(4+5/2,5*.9);
	\draw [line width=0.5mm] (3+5/2,5*.9)--(3+6/2,4*.9)--(3+7/2,3*.9)--(3+8/2,2*.9);
	\draw [line width=0.5mm] (2+6/2,4*.9)--(2+7/2,3*.9)--(2+8/2,2*.9)--(2+9/2,1*.9);
	\draw [line width=0.5mm] (6+2/2,8*.9)--(5+3/2,7*.9)--(4+4/2,6*.9);
	\draw [line width=0.5mm] (4+5/2,5*.9)--(3+6/2,4*.9)--(2+7/2,3*.9)--(1+8/2,2*.9);
	\draw [line width=0.5mm] (4+6/2,4*.9)--(3+7/2,3*.9)--(2+8/2,2*.9)--(1+9/2,1*.9);
	
	\draw [line width=0.5mm] (9+1/2,9*.9) -- (8+2/2,8*.9)--(8+3/2,7*.9)--(7+4/2,6*.9)--(6+5/2,5*.9)--(5+6/2,4*.9);
	\draw [line width=0.5mm] (8+2/2,8*.9)--(7+3/2,7*.9)--(7+4/2,6*.9);
	\draw [line width=0.5mm] (7+3/2,7*.9)--(8+3/2,7*.9);
	
	\draw [line width=0.5mm] (6+4/2,6*.9)--(5+5/2,5*.9);
	
	\end{scriptsize}
	\end{tikzpicture}
	
	\caption{Example of a GZ pattern associated to a fiber of a GZ system on a non-regular coadjoint orbit of $\U(10)$.}
	\label{fig:unitary-GZ-pattern}
	\end{figure}

\emph{Gelfand--Zeitlin systems} are a family of completely integrable systems
named for their connection to Gelfand--Zeitlin canonical
bases~\cite{GuilleminSternberg1983a},\footnote{\ The name ``Zeitlin'' is
variously romanized in the literature, also appearing as ``Cetlin'' and
``Tsetlin.'' } 
orthogonal Lie groups. 
The fibers of their moment maps, or
\emph{Gelfand--Zeitlin} \emph{fibers}, are interesting from several
perspectives, such as geometric quantization
\cite{GuilleminSternberg1983a,HamiltonKonno2014}, Floer theory
\cite{NishinouNoharaUeda2010,NoharaUeda2016,ChoKimOh2020}, and the topology of
integrable systems on symplectic manifolds \cite[Problem 2.9]{Alexey2018}. 
The
moment map images of \GZ systems on unitary and orthogonal coadjoint orbits
are polytopes known as \emph{Gelfand--Zeitlin polytopes} whose faces are
enumerated by combinatorial diagrams called \emph{Gelfand--Zeitlin
patterns}\footnote{\ Some authors prefer to use ladder diagrams, which are
equivalent.} (see Figure~\ref{fig:unitary-GZ-pattern}). 
The \GZ (henceforth
\emph{GZ}) fiber over a given point is naturally associated with the GZ pattern
of the face of the GZ polytope containing the point in its 
interior.

In recent work, Cho, Kim and Oh observed that both unitary and orthogonal GZ
fibers are total spaces of certain towers of fiber
bundles~\cite{ChoKimOh2020,ChoKim2020}. 
They used this description to show that
every unitary GZ fiber decomposes as a direct product of a torus and a space
whose first and second homotopy groups are trivial. 
They were also able to
explicitly describe the fibers in each stage of their towers as certain
products of spheres determined by the combinatorics of the associated pattern.
In particular, their results give a combinatorial formula for the dimension of
a fiber in terms of the associated pattern. 
Unitary GZ fibers were also
studied by Bouloc, Miranda, and Zung from a very different
perspective~\cite{BoulocMirandaZung2018}. 
In addition to recovering an
equivalent dimension counting formula, they show that every unitary GZ fiber
admits a product decomposition, finer than Cho--Kim--Oh's, in which direct
factors are enumerated by connected components of the associated GZ pattern. 
In
the case where a connected component is a diamond, they are able to simplify
their description of the corresponding factor to show that it is diffeomorphic
to a unitary group.

Although these recent results have greatly improved our understanding of GZ
fibers, they do not provide a clear description of the topology of GZ fibers in
all cases. 
For example, it is not clear how to extract from these results a
simple, useful description of the diffeomorphism type of the GZ fiber
corresponding to more complicated patterns such as that of 
\Cref{fig:unitary-GZ-pattern}.
This brings us to our first result. 
We modify Cho--Kim--Oh's tower construction 
to show that every unitary and orthogonal GZ fiber can also 
be expressed as a balanced product,
that is, a quotient of the form 
\begin{equation}\label{eq:intro-balanced-product}
	\begin{split}
		&
\quotientmed{H_\alpha \times H_{\alpha+1} \times \cdots \times H_{n} \,}
	    {\, L_\alpha \times L_{\alpha+1} \times \cdots \times L_{n} }\mathrlap, \\
&
(h_\alpha,\, h_{\alpha+1},\, \dots, \, h_n) 
	\sim 
\big(
h_\alpha\ell_\alpha^{-1}, \,
\psi_\alpha(\ell_\alpha)h_{\alpha+1}\ell_{\alpha+1}^{-1},\,
\dots,\,
\psi_{n-1}(\ell_{n-1})h_n\ell_n^{-1}
\big)
\mathrlap,
	\end{split}
\end{equation}
where the groups $H_k$ and $L_k$ are certain products of unitary and orthogonal
groups determined by the GZ pattern and $\psi_k\colon L_k \lt H_{k+1}$ are
certain injective homomorphisms (see Theorem \ref{thm:explicit-description}).
It follows immediately from this description that every GZ fiber is a direct
product whose factors correspond to connected components of the GZ pattern, a
result which is new in the orthogonal case (see Corollary
\ref{thm:component-product}). 
More importantly, we are able to systematically
simplify the expressions~\eqref{eq:intro-balanced-product} in terms of the
combinatorics of the associated GZ pattern. 
As a result of this 
simplification, we are able to provide relatively succinct
descriptions of the diffeomorphism types of all GZ fibers,
even those with relatively complex GZ
patterns such as in Figure~\ref{fig:unitary-GZ-pattern} (see
Example~\ref{intro-example} below). 
We illustrate these results with a number
of examples (Examples~\ref{eg:tori}--\ref{eg:intro example diffeo type}). 
In particular, our results recover the descriptions of elliptic non-degenerate
singular fibers and multi-diamond singular fibers given by Bouloc \emph{et al.}
as special cases. 

Our decomposition also allows us to extract torus 
factors from a GZ fiber with almost no effort,
not only in the unitary case, but also now in the orthogonal case
(\Cref{sec:circles}).
These torus factors are associated with a free torus action 
on each torus, which we are able to identify with the
action induced by the Thimm trick (\Cref{thm:Thimm}).
Moreover, our decomposition leads to weak local expressions
for a coadjoint orbit (\Cref{sec:face,sec:ray})
which allow us to parlay the extraction of our tori
into a topological model for a toric degeneration (\Cref{thm:det}).
That such degenerations exist was known for unitary GZ systems~\cite{NishinouNoharaUeda2010},
but an analogue for orthogonal GZ systems seems to remain open.


The decomposition allows us to extract the first three homotopy
groups of a unitary or orthogonal Gelfand--Zeitlin fiber from 
the GZ pattern with little effort (\Cref{sec:homotopy}).
We find that orthogonal fibers split as a product of a torus 
and a space whose first
and second homotopy groups have the forms $(\Z/2)^s$ and $\Z^f$ respectively.
We also show that the third homotopy groups of unitary and orthogonal GZ fibers
are free abelian
(\Cref{thm:pi-3-U,thm:pi-3-O}).
With a more substantial effort,
we can also compute the cohomology rings of arbitrary fibers of
unitary and orthogonal GZ systems 
(\Cref{sec:cohomology-U,sec:cohomology-O}).
We may briefly summarize these results as follows.

\begin{theorem}
	The integral cohomology of a unitary \GZ fiber is an exterior algebra 
	on odd-degree generators.
	The cohomology of an orthogonal \GZ fiber over $\Zf$ is also an
	exterior algebra; the integral and mod-$2$ cohomology groups are
	isomorphic to those of a product of real Stiefel manifolds.
\end{theorem}

\noindent The degrees of the relevant generators and the precise relevant
Stiefel manifolds are again determined in a straightforward manner from the
associated GZ pattern; see Theorems~\ref{thm:main-U} and~\ref{thm:main-O} for
more detail. 
The cohomology rings of these Stiefel manifolds is also 
known~\cite{cadekmimuravanzura}, so our description
gives much but not all of the integral cohomology ring structure for the GZ fiber.
These cohomology computations follow inductively from an analysis
of the Serre spectral sequences associated to Cho--Kim--Oh's towers.
Surprisingly, the spectral sequence associated to each bundle collapses in the
unitary case, and moreover, as the cohomology of base and fiber for each bundle
are exterior algebras, there is no extension problem, additively or
multiplicatively. 
In the orthogonal case, the spectral sequences associated to the 
sphere bundles do not typically collapse,
but those associated to allied Stiefel manifold bundles do. 

\begin{example}\label{intro-example}
 Consider the GZ pattern in \Cref{fig:unitary-GZ-pattern}. 
 Using our results, one can immediately read from this 
 pattern that an associated GZ fiber is diffeomorphic to
\[
	(S^1)^7 \,\x\, (S^3)^3 \,\x\,
	\U(2)\backslash \big(\U(4) \times \U(3) \big)/\U(2),
\]
has integral cohomology ring isomorphic to 
\[
\ext[z_{1,1},
	z_{1,2},
	z_{1,3},
	z_{1,4},
	z_{1,5},
	z_{1,6},
	z_{1,7},
	z_{3,1},
	z_{3,2},
	z_{3,3},
	z_{5,1},
	z_{5,2},
	z_{7,1}],
\qquad
|z_{m,j}| = m,
\]
and has $\pi_3 \cong \Z^3$.
See Examples~\ref{eg:intro example diffeo type}
and~\ref{eg:intro example cohom}
for more details.
\end{example}

We end this introduction with a discussion of the broader context and
motivation for our work. Gelfand--Zeitlin systems are useful and interesting because they share
many features with toric integrable systems (convexity and global action-angle
coordinates) but have non-toric singularities.
They form part of a larger family of examples with
similar behaviour that also includes bending flow systems on polygon space
\cite{KM}, Goldman systems on moduli space \cite{Gold}, and integrable systems
constructed by toric degeneration~\cite{harada-kaveh} (of which GZ systems are
an example \cite{NishinouNoharaUeda2010}).
This family of examples continues to grow.
For instance, the construction of integrable systems by toric
degeneration was recently extended~\cite{HL},
and as a consequence,
it was shown that every coadjoint orbit of every 
compact, connected Lie group admits
integrable systems that share the features of GZ systems mentioned above. 

There are many questions regarding the non-toric singularities of these
systems which are currently unanswered.
For example, it is known in
the cases of bending flow systems \cite{Bouloc2018} and GZ systems
\cite{BoulocMirandaZung2018,ChoKimOh2020} that the non-toric fibers are
isotropic (this result is expected to have applications in Lagrangian Floer
theory \cite{ChoKimOh2020}), but it is not known whether non-toric fibers of
integrable systems constructed by toric degeneration are isotropic in general
(this is conjectured to be true by Bouloc \emph{et
al.}~\cite{BoulocMirandaZung2018}).
Similarly, very little is known in general
about the topology of non-toric fibers of integrable systems constructed by
toric degeneration. There is also the problem of finding local normal forms
for the non-toric singularities of these systems, as motivated by Bolsinov
\emph{et al.} for the case of GZ systems~\cite[Problem 2.9]{Alexey2018}.
At the moment, very little is known about this problem,
even in the case of the more concrete GZ systems,
with the exception of several examples in low
dimension studied by Alamiddine~\cite{Alamiddine2009}.
Progress on any of these questions would have further applications in symplectic topology,
geometric quantization, and the topology of integrable systems.

To give a particularly concrete example of one such application, we end by
describing one of the original motivations for studying GZ systems (and a
potential future application of this work). Guillemin and Sternberg observed
that if integral points in the boundary of an integral GZ polytope are formally
included in the Bohr--Sommerfeld set, then the dimension of the resulting
Bohr--Sommerfeld quantization equals the dimension of the K\"ahler
quantization of the coadjoint orbit \cite{GuilleminSternberg1983a}. In ongoing
work, we plan to use our description of the topology of the singular GZ fibers,
and a description of a local model for the GZ systems in neighbourhoods of said
fibers, to give a more principled justification of Guillemin and Sternberg's
observation, i.e.,~to prove that a singular fiber is Bohr--Sommerfeld if and
only if it is integral. Moreover, we hope to extend the results of
Hamilton--Kono~\cite{HamiltonKonno2014} by showing that holomorphic sections in
the K\"ahler quantization that correspond to boundary points of the GZ polytope
converge under a deformation of complex structure to distributional sections
supported on the singular fibers over the same boundary points.

\medskip

\noindent\emph{Organization of the paper.}
In Section~\ref{sec:Gelfand--Zeitlin-systems}, 
we recall the description of GZ fibers as towers. 
In particular, we recall results of Cho--Kim--Oh
which describe the bundles of homogeneous spaces that occur at each stage 
in terms of the associated GZ pattern~\cite{ChoKimOh2020,ChoKim2020}. 
We note a minor correction 
to Cho and Kim's description~\cite{ChoKim2020} of the bundles in the orthogonal case 
(although this correction does not impact their final results,
it does play an important role in our account).
Section~\ref{sec:product} contains our description of the GZ fibers as biquotients (Theorem~\ref{thm:biquotient}). 
This presentation shows that GZ fibers decompose as direct products
of factors indexed by components of the associated GZ pattern.
We provide examples demonstrating how easily our biquotient description can be read from the associated GZ pattern.
We also demonstrate how this description enables
us to recover a direct torus factor, 
and relatedly a torus action
and a (topological) toric degeneration;
this last depends on the weak local normal forms 
over the interior of a face and the interiors of two adjoining faces 
developed in \Cref{sec:local}.
In \Cref{sec:homotopy},
we use this description to determine the first three homotopy groups
of a GZ fiber.
In Section~\ref{sec:cohomology-U}, 
we state and prove the main theorem (Theorem \ref{thm:main-U})
regarding cohomology of GZ fibers in the unitary case. 
Finally, Section \ref{sec:cohomology-O} contains our cohomological results in the orthogonal case. 

\medskip

\noindent\emph{Acknowledgements.}
Both authors were supported by a Fields Postdoctoral Fellowship
      associated with the 2020 Thematic Program in Toric Topology 
during the original collaborative work that lead to this manuscript.

\section{Gelfand--Zeitlin systems and Gelfand--Zeitlin fibers}\label{sec:Gelfand--Zeitlin-systems}

This section recalls the Gelfand--Zeitlin systems on unitary and orthogonal
coadjoint orbits as well as the description of their fibers as towers. In
order to avoid repeating ourselves, we begin in Section \ref{sec:general-case}
by giving a general description of GZ systems and fibers that applies equally
in both cases. We recall the details that are specific to the unitary and
orthogonal cases in Sections~\ref{sec:GZ-fibers-unitary-case} and
\ref{sec:GZ-fibers-orthogonal-case} respectively. 

For further details regarding GZ systems, we direct the reader to the original
papers of Guillemin and Sternberg
\cite{GuilleminSternberg1983a,GuilleminSternberg1983b} and the recent papers of
Cho--Kim--Oh describing the GZ fibers as towers of bundles of homogeneous
spaces \cite{ChoKim2020, ChoKimOh2020}. As mentioned in the introduction, an
alternative description of the GZ fibers is given by Bouloc \emph{et
al.}~\cite{BoulocMirandaZung2018}. We also remark that the Ph.D.~thesis of
Milena Pabiniak~\cite{Pabiniak2012} contains many useful details.

\subsection{Some known facts and set-up in the general case}\label{sec:general-case}

This section recalls why fibers 
of both unitary and orthogonal \GZ systems
are towers of bundles of homogeneous spaces,
in terms of a much more general set-up.

\begin{definition}\label{def:GZ-abstract}
Let positive integers $\defm \alpha < \defm n$ 
be given along with $\defm{G_k}$ ($\alpha \leq k \leq n+1$) topological groups
connected by injective homomorphisms $\defm\vp = \defm{\vp_k}\: G_k \lt G_{k+1}$
($\alpha \leq k \leq n$)
and $G_k$-spaces $\defm{V_k}$ ($\alpha \leq k \leq n+1$)
connected by maps
$\defm\Phi = \defm{\Phi_{k+1}}\: V_{k+1} \lt V_k$
($\a \leq k \leq n$).
We write $\defm{\Orb_{\genup k}} \ceq G_k \. \genup k$
for the orbit of $\genup k \in V_k$
and suppose that
\bitem
\item
  Each $\Phi_{k+1}$ is $\vp_k$-equivariant in the sense that 
  for all $g_k \in G_k$ and $\genup{k+1} \in V_{k+1}$,
  \[
  \Phi_{k+1}\big(\vp_k(g_k)\genup{k+1}\big) = g_k \Phi_{k+1}(\genup{k+1});
  \]
\item
  The group $\vp_k({G_k}) \leq G_{k+1}$ acts transitively
  on any nonempty intersection 
  $
  \defm{\Orb_{\genup{k+1},\genup k}} \ceq
  \Orb_{\genup{k+1}} \inter \Phi\-(\Orb_{\genup k})
  $.
\eitem
Fix $\defm{\genup{n+1}} \in V_{n+1}$ and write 
$\defm\Orb = \defm{\Orb_{\genup{n+1}}}$.
Taking orbits of iterated $\Phi$-images yields a natural map 
\eqn{
\defm\Psi\: \mc O &\lt \prod_{k=\alpha}^{n} V_{k}/G_{k} \eqc \defm\D,\\
\genup{n+1} &\lmt \big(\Orb_{ \Phi^{n+1-k}(\genup {n+1})}\big)_{k=\a}^n .
}
We call the map $\Psi$ a \defd{\TGZS}
and will attempt to characterize 
the \defd{Gelfand--Zeitlin fiber}
$\Psi\-(p)$
over a point $\defm p = (\defm{\bar\lambda^{( k )}})_{k=\a}^{n} \in \D \inter \im \Psi$.
\end{definition}

\begin{definition}\label{def:GZ-system}
If 
$\vp_k\: G_k \lt G_{k+1}$ are injective Lie group homomorphisms,
$V_k = \fg_k^*$ the coadjoint representations,
and $\Phi_{k+1} = \vp_k^*\: \fg_{k+1}^* \lt \fg_k^*$
the cotangent maps,
$\Psi$ from \Cref{def:GZ-abstract}
is called a \defd{Gelfand--Zeitlin system}.
We will call a map defined in the general context
\emph{smooth}
if it becomes smooth in the case of a GZ system.
The main cases are the 
\defd{unitary} and 
\defd{orthogonal GZ systems} 
with $G_k = \U(k)$ or $\SO(k)$ 
and $\a = 1$ or $2$ (respectively)
and each $\vp_k$ a block inclusion 
$g 
\lmt 
\big[\begin{smallmatrix}
	g&\\ &1
\end{smallmatrix}\big].
$

Fix a maximal torus $\defm{T_k}$ in each $G_k$, 
and fix positive Weyl chambers $\defm{\ft_{k,+}^*}$ 
in the linear duals~$\ft^*_k$ of their Lie algebras.
Since each $\Ad^*_{G_k}$-orbit on $V_k = \fg_k^*$
meets $\ft_{k,+}^*$ in precisely one point,
we identify $\ft_{k,+}^*$ with the orbit space
$V_k/G_k = \fg_k^*/G_k$,
making $\Psi$ a map $\Orb \lt \prod_{k=\a}^n \ft_{k,+}^*$.
We will write $\smash{\defm{\lambda^{(k)}}}$ for points of $\ft_{k,+}^*$,
so that we are interested in the spaces $\Psi\-(\lup n,\dots,\lup 1)$.

\end{definition}


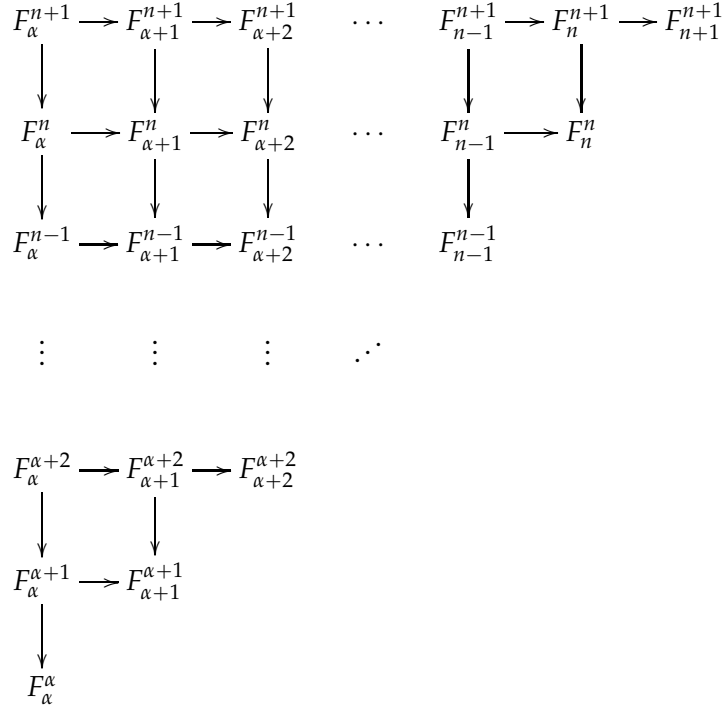
\begin{figure}
\centerline{
	\xymatrix@C=1.3em@R=2em{
		{F^{n+1}_\alpha} \ar[d]\ar[r]		&
		F^{n+1}_{\alpha+1}\ar[d]\ar[r]	&
		F^{n+1}_{\alpha+2}\ar[d]		& 
		\cdots 				&
		F^{n+1}_{n-1} \ar[r]\ar[d]	& 
		F^{n+1}_n\ar[r]\ar[d]		& 
		F^{n+1}_{n+1}			\\	
		F^{n}_\alpha \phantom{{}^{1}}\ar[r]\ar[d]
							&
		F^{n}_{\alpha+1} \ar[r]\ar[d]		&
		F^{n}_{\alpha+2}\ar[d]		& 
		\cdots				&
		F^n_{n-1} \ar[d]\ar[r]		& 
		F_n^n 				\\
		F^{n-1}_\alpha 		\ar[r]		&
		F^{n-1}_{\alpha+1} 	\ar[r]		&
		F^{n-1}_{\alpha+2} 			& 
		\cdots				&
		F_{n-1}^{n-1} 			\\
		 \vdots				&
		 \vdots 				&
		 \vdots 				& 
		 \reflectbox{$\ddots$}		\\
		F^{\alpha+2}_\alpha \ar[d] \ar[r]		&
		F^{\alpha+2}_{\alpha+1} \ar[r]\ar[d]	& 
		F_{\alpha+2}^{\alpha+2}			\\
		F^{\alpha+1}_\alpha\ar[d]\ar[r]		& 
		F_{\alpha+1}^{\alpha+1}			\\
		F_\alpha^\alpha				\\
	}
}
\caption{Gelfand--Zeitlin fibers can be viewed as the total space $F^{n+1}_\alpha$ 
of a tower of fiber bundles (on the left side of the diagram) 
which is constructed through a sequence of pullbacks.}
\label{fig:big-diagram}
\end{figure}

\begin{remark}\label{rmk:GZ-system-known-facts}
Guillemin and Sternberg showed
that a unitary or orthogonal Gelfand--Zeitlin system $\Psi$ 
generates a Hamiltonian action of $T_n \times \dots \times T_\alpha$ 
on an open, dense subset of $\mathcal{O}_{\lambda^{(n+1)}}$, 
with respect to the canonical Kostant--Kirillov--Souriau symplectic structure 
on $\mathcal{O}_{\lambda^{(n+1)}}$%
~\cite{GuilleminSternberg1983a}.
They also showed
in these cases that
the map $\Psi$ defined in \eqref{def:GZ-abstract}
is a GZ system on any coadjoint orbit $\mathcal{O}_{\lambda^{(n+1)}}$~\cite{GuilleminSternberg1983b}.
In both cases, the resulting densely defined torus action is completely integrable,
or toric, for every $\lambda^{(n+1)}$ in $\ft_{n+1,+}^*$. 
It was shown in a more general setting by the second author that the dense subset 
where $\Psi$ generates a torus action is connected and the fibers of $\Psi$ 
are all connected~\cite{Lane2018}. 
\end{remark}

\begin{discussion}[The staircase and pullbacks]\label{rmk:lambda-staircase}

Following \Cref{def:GZ-system},
for each $p = (\bar\lambda^{( k )}) \in \D$
we fix orbit representatives $\defm{\lup k} \in V_k$ of $\bar\lambda^{( k )}$ for all $k$.
(In the cases of unitary and orthogonal systems,
there will be very natural such choices, amounting to a continuous section.)
To say a point $\genup{n+1} \in V_{n+1}$ lies in $\Psi\-(p)$ means
that each $\defm{\genup k} \ceq \Phi^{n+1-k}(\genup{n+1})$
lies in ${\Orb_{\lup k}}$,
and that each $\genup {k+1}$ lies in ${\Orb_{\lup{k+1},\lup k}}$.
Conversely, we can ask when 
a system $(\z^{(k)}) \in \prod {\Orb_{\lup{k+1},\lup k}}$
comes from some $\genup{n+1}$.
Let $\defm{\Psi_i}\: V_k \lt V_i/G_i$
be the composite of $\Phi^{k-i}\: V_k \lt V_i$ and the projection,
for any $i \leq k$.
For $\a \leq j \leq k \leq n+1$,
we denote by $\defm{\smash{F^k_j}} \sub V_k$ the space of points $\z^{(k)}$ with
$\Psi_i(\z^{(k)}) = \lup i$ for $j \leq i \leq k$.
Then $\Psi\-(p) = F^{n+1}_\alpha$.
Following Cho--Kim--Oh~\cite{ChoKimOh2020},
these all fit together into the system
shown in \Cref{fig:big-diagram},
in which all horizontal maps are inclusion and all vertical maps are $\Phi$.
Each rectangle in \Cref{fig:big-diagram} in fact represents a fiber product
in a natural way.\footnote{\
  Indeed,
  for $i < j < k < m$,
  if $\z^{(k)} \in \smash{F^k_i} \sub \smash{F^k_j}$ is the image under $\Phi^{m-k}$
  of $\z^{(m)} \in \smash{F^m_j}$,
  then it means $\Psi_s(\z^{(m)}) = \lup s$ for $i \leq s \leq k$
  as well as for $j \leq s \leq m$.
}
Since each rectangle is a fiber product,
this means $\Psi\-(p)$
can be calculated inductively as the an iterated fiber product if we 
can just understand the maps in \Cref{fig:big-diagram}.
\end{discussion}

%

\begin{discussion}[Stabilizers and the vertical maps]\label{rmk:stabilizer-bundle}
We work in the context of a \TGZS
with orbit representatives $\lup k$ chosen as in \Cref{rmk:lambda-staircase},
summarizing the work of Cho--Kim--Oh.\footnote{\ 
	With (N.B.) somewhat different notation and emphases.
	For example, they do not seem to explicitly 
	state \eqref{eq:HkLk-bundle} is a bundle of homogeneous spaces,
	though it follows from their analysis, 
	but this will be crucial for us.
}
For each $k \leq n$, 
fix $\defm{\xip{k+1}} \in \Orb_{\lup{k+1}} \inter \Phi\-(\lup k)$.\footnote{\ 
	This is possible as,
	since we assume $p$ lies in $\im \Psi$,
	we have $\lup k \in G_k \. \Phi(\xip{k+1})$
	for some $\xip{k+1} \in G_{k+1}\.\lup{k+1}$,
	and hence
	there is some $g_k$ with $\lup k = g_k \.\Phi(\xip{k+1}) = 
		\Phi\big(\vp(g_k)\xip{k+1}\big)$.
}
We write 
\[
\defm{L_k} \ceq \vp\-(\Stab_{G_{k+1}} \xip{k+1}) 
	\leq 
\Stab_{G_k} \lup k \eqc \defm{H_k}
\mathrlap.
\]
Then $\Orb_{\lup k} = G_k \. \lup k$ 
is diffeomorphic to $G_k/H_k$ by the orbit--stabilizer theorem.
Since we have assumed the action of 
$\vp(G_k)$ on $\Orb_{\lup{k+1},\lup k}$ is transitive,
we see $\Orb_{\lup{k+1},\lup k}$ is diffeomorphic to $G_k/L_k$,
and the fiber of the restriction $\Phi\: \Orb_{\lup{k+1},\lup k} \lt \Orb_{\lup k}$
over $\lup k$ itself can be identified with $H_k/L_k$.
In the case of a GZ system, since all groups are Lie groups,
it is well known that 
\begin{equation}\label{eq:HkLk-bundle}
	H_k/L_k \to G_k /L_k \to G_k/H_k
\end{equation}
is a fiber bundle.
Thus from \Cref{fig:big-diagram}, in which every rectangle now represents
a pullback of fiber bundles,
we see that $\Psi\-(p)$ is an iterated fiber bundle
with fibers $H_\a/L_\a$, \ldots, $H_n/L_n$.
\end{discussion}

\subsection{The unitary case}\label{sec:GZ-fibers-unitary-case}
	
This section establishes notation and describes the groups $L_k\leq H_k \leq G_k$ 
of \Cref{rmk:stabilizer-bundle} for 
unitary \GZ systems 
following Cho--Kim--Oh~\cite{ChoKimOh2020}.


\begin{discussion}\label{rmk:unitary-setup}
Recall from \Cref{def:GZ-system}
that in this case $G_k  = \U(k)$ for $1 \leq k \leq n+1$.
The coadjoint representation 
$\f u(k)^*$ can be identified with
the conjugation action of $\U(k)$ on the space of $k\times k$ Hermitian matrices.\footnote{\ 
	It would be more natural, but less convenient,
	to consider the space skew-Hermitian matrices,
	which is equivariantly isomorphic under multiplication by $i$.
}
Then $\Phi_{k+1} = \vp_k^*$ becomes the truncation map extracting 
from a $(k+1) \x (k+1)$ matrix
its upper-left $k \x k$ submatrix.
We select the diagonal matrices of 
$\U(k)$ as the maximal torus $T_k$;
then the positive Weyl chamber $\ft_{k,+}^*$ 
can be identified with the set of points 
$\smash{\lambda^{(k)} = 
(\defm{\smash{\lambda^{(k)}_1}}, \dots, 
\defm{\smash{\lambda^{(k)}_k})}} \in \R^k$
with $\smash{\lambda^{(k)}_1} \geq \dots \geq \smash{\lambda^{(k)}_k}$,
so that 
the coadjoint orbit $\mathcal{O}_{\lambda^{(k)}}$ is
the set of $k\times k$ Hermitian matrices with eigenvalues,
listed in weakly decreasing order,
$\smash{\lambda^{(k)}_1} , \dots , \smash{\lambda^{(k)}_k}$.
\end{discussion}

\begin{discussion}[The GZ polytope and GZ patterns]
It is a result of linear algebra that $\mathcal{O}_{\lambda^{(k+1)},\,
\lambda^{(k)}}$ is non-empty if and only if $\lambda^{(k)}$ and
$\lambda^{(k+1)}$ satisfy \defd{interlacing inequalities}
\begin{equation}\label{eq:interlacing-U}
 \begin{array}{@{}c@{\;}c@{\;}c@{\;}c@{\;}c@{\;}c@{\;}c@{\;}c@{\;}c@{\;}c@{\;}c@{\;}c@{\;}c@{\;}c@{\;}c@{\;}c@{\;}c@{}}
 \lambda_1^{(k+1)} & & & & \lambda_{2}^{(k+1)} & & & &
 \lambda_{3}^{(k+1)} & & & & \lambda_{k}^{(k+1)} & && &
 \lambda_{k+1}^{(k+1)}. \\ &\dgeq & & \ugeq & &\dgeq & & \ugeq &
 &\dgeq&&\ugeq&&\dgeq &&\ugeq & \\ && \lambda_{1}^{(k)} &&&&
 \lambda_{2}^{(k)} && & &\cdots &&&& \lambda_{k}^{(k)} &&\\ \end{array}
 \end{equation}
 The image of the GZ system on $\mathcal{O}_{\lambda^{(n+1)}}$
 is the set $\D = \D_{\lup{n+1}}$
 of those points in $\prod_{k=1}^n \R^k \iso \R^{n(n+1)/2}$ 
 whose coordinates
 $\smash{\lup{k}_j}$ ($1 \leq j \leq k \leq n$)
 satisfy the interlacing
 inequalities \eqref{eq:interlacing-U} for all
 $k$~\cite{GuilleminSternberg1983a}. 
 Since each inequality cuts out a half-space,
 it follows $\D \subn \R^{n(n+1)/2}$ is a polytope,
 called the \defd{unitary Gelfand--Zeitlin 
 polytope} associated to $\lambda^{(n+1)}$.
 We 
 represent these inequalities in a triangular array, 
 such as the following pattern for the case $n=2$.
\begin{equation}\begin{array}{@{}c@{\;}c@{\;}c@{\;}c@{\;}c@{\;}c@{\;}c@{\;}c@{\;}c@{\;}c@{\;}c@{\;}c@{\;}c@{\;}c@{\;}c@{\;}c@{\;}c@{}}
 \lambda_{1}^{(3)} & & & & \lambda_{2}^{(3)} & & & & \lambda_{3}^{(3)} \\
 &\dgeq & & \ugeq & &\dgeq & & \ugeq & \\
 && \lambda_{1}^{(2)} &&&& \lambda_{2}^{(2)} && \\
 &&&\dgeq && \ugeq && &\\
 &&&& \lambda_{1}^{(1)} &&&& 
\end{array}	
\end{equation}
To each point in $p = (\lup k) \in \Delta_{\lambda^{(n+1)}}$, 
we associate its
\defd{unitary Gelfand--Zeitlin pattern},
a graded planar graph constructed 
from the triangular array of inequalities,
by replacing each number $\smash{\lup{k}_j}$ with a vertex
and drawing an edge between two nearest-neighbour vertices 
(in the same row or adjacent rows)
if and only if the associated $\lup k _j$ are equal.%
\footnote{\ 
	These are equivalent to the \emph{ladder diagrams} of Cho--Kim--Oh.
	To convert from one of our GZ patterns to a ladder diagram, draw a 
	grid in gray such that each $\smash{\lup k _j}$ lies in a (diamond) box,
	and for each strict inequality, which will cross one grey box edge,
	color that edge black. 
	Then remove the top row of entries $\smash{\lup{n+1}_j}$
	and rotate the diagram by $-\pi/4$ so that $\smash{\lup 1 _1}$ is at
	the lower-left.
}
Such a graph will look something like \Cref{fig:unitary-GZ-pattern}.
Note that the pattern of $p$ depends only on the face of
$\Delta$ containing $p$ in its interior 
and not on the particular values $\lup k_j$.
\end{discussion}

\begin{discussion}[Rows and the groups $H_k$]\label{rmk:unitary-row}
\defd{Row $k$} of a unitary GZ pattern is the full subgraph on the
row of vertices corresponding to the numbers 
$\lup k _1$, \ldots, $\lup k _k$;
this linear graph
will comprise a number (call it $\defm{\ell} = \ell_k$)
of component segments each corresponding to a distinct entry of $\lup k$.
If $\defm{d_i}$ ($1 \leq i \leq \ell$) 
is the number of vertices in the $i\th$ component from the left,
the centralizer 
of the canonical diagonal orbit representative $\lup k$ from \Cref{rmk:unitary-setup}
is the block-diagonal subgroup 
	\begin{equation}\label{eq:J-U}
		H_k = \U(d_1) \oplus \dots \oplus \U(d_{\ell}) \leq \U(k).
	\end{equation}
\end{discussion}

\begin{discussion}[Representatives $\xip{k+1}$ and groups $L_k$; 
	  trapezoids and parallelograms]\label{rmk:xi-explicit}
	  Continuing with the notation of \Cref{rmk:unitary-row},
	  let $\defm{\mu_i}$ be the distinct entries of $\l^{(k)}$,
	in decreasing order.
	An element of $ \Phi\-(\l^{(k)}) \inter \Orb_{\lup{k+1}}$,
	as in \Cref{rmk:stabilizer-bundle},
	is then constrained to be of the block form 
	\begin{equation}\label{eq:Y}
	  \xip{k+1} 
	  	=
	  \left[
	  \begin{array}{c|c|c|c}
	    \diag(\mu_1,\cdots,\mu_1)&  		&		&z_1\\
	    \hline
	    	&\ddots&	& \vdots	\\
	    \hline
	    &	& \diag(\mu_{\ell},\ldots,\mu_{\ell})		&z_\ell \\
	    \hline
	    \bar z_1^\top		&\cdots		& \bar z_{\ell}^\top&	c
	  \end{array} 
	  \right]\mathrlap.
	\end{equation} 
	Here, because trace is a conjugacy invariant,
	we must have $\defm c = \tr \l^{(k+1)} - \tr \l^{(k)} \in \R$,
	and we will show momentarily that the
	squares of the Euclidean norms $\defm{r_i} = \|z_i\|$
	of $\defm{z_i} \in \C^{d_i}$ 
	are also rational functions of (and hence continuous in)
	the entries $\lup{k}$ and $\lup{k+1}$.

For an element $h_k = u_1 \+ \cdots \+ u_\ell \in H_k$,
one checks the conjugation action of $\vp(h_k)$ on $\xip{k+1}$
fixes the diagonal entries and replaces each column vector $z_i$
with $u_i\.z_i$.
Since $\U(d_i)$ acts transitively on the sphere of radius $r_i$ in $\C^{d_i}$,
we see $H_k/L_k$ is diffeomorphic to the product of spheres
$\prod_{r_i \neq 0} S^{2d_i-1}$.
We usually choose $\xip{k+1}$
so that $z_i = [r_i\ 0\ \cdots\ 0]^\top$ (whether or not $r_i = 0$).
For this choice, 
the group $L_k \leq H_k$ such that $L_k \+ [1]$ centralizes $\xip{k+1}$
is the block sum 
\begin{equation}\label{eq:K-U}
  L_k = K_1 \oplus \cdots \oplus K_{\ell} \leq \U(d_1) \oplus \dots \oplus \U(d_{\ell}) =  H_k
  \mathrlap,
\end{equation}
where $\defm{K_i} \leq \U(d_i)$ 
is $\U(d_i)$ itself when $r_i = 0$
and $[1] \+ \U(d_i-1)$ when $r_i > 0$.

This dichotomy can be interpreted usefully in terms of the GZ pattern.
Consider the full subgraph of the GZ pattern on the vertices of rows $k$ and $k+1$.
A connected component of this subgraph 
will be called a
\defd{\mtrap-shape} if it has more vertices in row $k$, a
\defd{\wtrap-shape} if it has more in row $k+1$, and
a \defd{\parallelogram-shape} if it contains the same number in each row.%
\footnote{\
  The terminology of
  \parallelogram-shape/\wtrap-shape/\mtrap-shape is equivalent to the
  N-block/W-block/M-block trichotomy used for the ladder diagrams 
  of Cho--Kim--Oh~\cite{ChoKimOh2020}.
}
Since 
a nonzero $r_i$ forces down the multiplicity of $\mu_i$
as an eigenvalue 
of $\xip{k+1}$,
while $r_i = 0 $ guarantees an multiplicity
of at least $d_i$,
it follows that a \mtrap-shape corresponds to $K_i = [1] \+ \U(d_i)$
while a \parallelogram- or \wtrap-shape corresponds to $K_i = \U(d_i)$.
Hence sphere factors of $H_k/L_k$ arise from \mtrap-shapes in the pattern.

We will need the formula for $r_i^2$.	
Write $\defm{\l_s}$ for the distinct entries of $\lup{k+1}$
and equate the characteristic polynomials of $\lup{k+1}$ and $\xip{k+1}$
(expanding by minors along the last column for~$\xi$).
Dividing out common factors $x - \smash{\lup{k+1}_j}$
and evaluating at $x = \mu_i$ shows that if $\mu_i$ belongs to a \wtrap- or \parallelogram-shape
in the $k\th$ and $(k+1)\st$ rows of the GZ pattern,
then $r_i = 0$,
and otherwise
\quation{\label{eq:r}
r_i^2 = - \frac{ \prod_w (\mu_i - \l_w) }{ \prod_{m \neq i} (\mu_i - \mu_m) }
\mathrlap,
}
where the indices $w$ run over \wtrap-shapes and the $m$ run over \mtrap-shapes.
\end{discussion}

\brmk\label{rmk:xi-prime}
In exceptional situations later,
	we will use a representative $'\xip{k+1}$ as in \Cref{rmk:xi-explicit}
	which is like our usual choice of $\xip{k+1}$ 
	except that $z'_i = [0\ \cdots\ 0 \ r_i]^\top$
	and the corresponding block factor of 
	its stabilizer $L'_k$ is $\U(d_i-1)\+[1]$.
	Note that $'\xip{k+1} = \vp_k(h) \xip{k+1}$
	for a block-unitary matrix $h = u_1 \+ \cdots \+ u_\ell  \in H_k \in \prod_{i=1}^\ell \SU(d_i)$.\footnote{\ 
	For example,
	this can be taken such that $u_{s}$ is the identity matrix
	for $s \neq i$ 
	and $u_i \in \U(d_i)$ 
	is $\pm$-signed permutation matrix
	with a $1$ in the first column.
	}
\ermk


\begin{example}\label{eg:HkLk}
 Consider the unitary GZ pattern in \Cref{fig:unitary-GZ-pattern}. 
 Looking at rows $4$ and $5$, we see that $H_4 = \U(4)$ and $L_4 =
 [1]\oplus \U(3)$. 
 Looking at rows $6$ and $7$, we see that
 \[ 
 H_6 = \U(2)\oplus
 \U(2) \oplus \U(1) \oplus \U(1),
 \quad L_6 = \big([1]\oplus \U(1)\big) \oplus \U(2)
 \oplus\U(1) \oplus \U(1).
 \]
 Looking at the penultimate row $9$, we see $H_9
 \iso \U(2) \+ \U(1)^{\oplus 2} \+ \U(2) \+ \U(1)^3$ whereas $L_9 \iso \U(2)
 \+ \U(0)^{2} \+ \U(2) \+ \U(0) \+ \U(1) \+ \U(0)$. 
 We only look at row $10$
 in conjunction with row $9$ to determine what the \mtrap-, \wtrap-, and
 \parallelogram-shapes are; while $H_{10}$ is defined, $L_{10}$ is not. 
 Correspondingly, it does not make sense to refer to \mtrap-shapes in row $10$,
 since there is no row $11$. 
\end{example}

\subsection{The orthogonal case}\label{sec:GZ-fibers-orthogonal-case}

This section establishes notation and describes the groups $L_k\leq H_k \leq G_k$ 
of \Cref{rmk:stabilizer-bundle} for 
orthogonal \GZ systems,
mostly following Cho--Kim~\cite{ChoKim2020}.%
\footnote{\ 
	The description of the subgroups $H_k$ and $L_k$ given
	here diverges slightly from theirs~\cite[Lem.~6.5]{ChoKim2020},
	which gives blocks $\SO(2d_i)$
	where we find $\U(d_i)$.
	As $\SO(2d_i)/\SO(2d-1)$ 
	and $\U(d_i)/\U(d_i-1)$ are both diffeomorphic to $S^{2d_i-1}$,
	the quotients $H_k/L_k$ they give are fortunately correct anyway,
	so this minor error does not affect the rest of their results,
	but the distinction will be important for us.
}
Throughout this section, we write $\defm{[k]}$ for $\lfloor k/2\rfloor$, 
where $\defm{\lfloor \cdot \rfloor}$ is the floor function.

\begin{discussion}
Recall from \Cref{def:GZ-system} that in this case $G_k = \SO(k)$ for $2 \leq k \leq n+1$.
The coadjoint representations 
$\f {so}(k)^*$ can be identified with
the conjugation action of $\SO(k)$ on the space of $k\times k$ skew-symmetric matrices.
Again $\Phi_{k+1} = \vp_k^*$ becomes the truncation map extracting 
from a $(k+1) \x (k+1)$ matrix
its upper-left $k \x k$ submatrix.
We select the block-diagonal matrices $\SO(2)^{\+[k]}$ as the maximal torus $T_k$
for $k$ even, and $\SO(2)^{\+[k]} \+ [1]$ for $k$ odd,
and write $\defm{\b(\l)} \ceq
\smash{\big[\begin{smallmatrix}
	0&-\l\\ \l&\phantom{-}0
      \end{smallmatrix}\big]}$ for $\l \in \R$.
Under this identification, 
$\f t_k^*$ is
for $k$ even
the space of matrices of the block-diagonal form
$\smash{\Direct_{j = 1}^{[k]} \b(\l_i)}$ for various $\l_i \in \R$,
and for $k$ odd of matrices $\smash{\Direct_{j = 1}^{[k]} \b(\l_i) \+ [1]}$.
 The standard positive Weyl chamber $\ft_{k,+}^*$ 
is the set of such matrices with 
	\begin{equation}\label{eq:sequence-O-even}
	  \begin{gathered}
	  \begin{aligned}
	    \phantom{\qquad k \mbox{ even},}
	\lambda_1^{(k)} \geq
	\lambda_2^{(k)} \geq
	\cdots \geq
	\lambda_{[k]-1}^{(k)} &\geq 
	\smash{	\big|\lambda_{[k]}^{(k)}\big| }
	,
&	\qquad k \mbox{ even},
	\\
	\phantom{\qquad k \mbox{ odd},}
	\lambda_1^{(k)} \geq
	\lambda_2^{(k)} \geq
	\cdots \geq
	\lambda_{[k]-1}^{(k)} &\geq 
	\lambda_{[k]}^{(k)}
	\geq 0,
&	\qquad k \mbox{ odd}.
	\end{aligned}
      \end{gathered}
\end{equation}
For brevity, we abusively write the standard conjugacy class representatives
in $\ft_{k,+}^*$ again as $\defm{\lup k}$.
Each of these $\lup k_j$ corresponds to the pair of complex conjugate eigenvalues $\pm i \lup k_j$;
for $k$ odd there is also one additional instance of $0$ in the list of eigenvalues.
\end{discussion}


\begin{discussion}[The GZ polytope and GZ patterns]
  The interlacing inequalities characterizing for which $\lambda^{(k)}$ and $\lambda^{(k+1)}$
  the intersection $\mathcal{O}_{\lambda^{(k+1)},\,\lambda^{(k)}}$ is non-empty again\footnote{\ 
  For nonemptiness, see Pabiniak~\cite[App.~B]{Pabiniak2012}.
}
  depend on the parity of $k$.
  For $k$ odd, they are
	\begin{equation}\label{eq:interlacing-O-even}
	\begin{array}{@{}c@{\;}c@{\;}c@{\;}c@{\;}c@{\;}c@{\;}c@{\;}c@{\;}c@{\;}c@{\;}c@{\;}c@{\;}c@{\;}c@{\;}c@{\;}c@{\;}c@{}c@{\;}c@{}c@{\;}c@{}c@{\;}c@{}}
 	\lambda_1^{(k+1)} & & & & \lambda_{2}^{(k+1)} & & & & \lambda_{3}^{(k+1)} & & & & \lambda_{[k+1]-1}^{(k+1)} & && & \big|\lambda_{[k+1]}^{(k+1)}\big| &&\\
	&\dgeq & & \ugeq & &\dgeq & & \ugeq & &\dgeq&&\ugeq&&\dgeq &&\ugeq & &\dgeq&\\
	&& \lambda_{1}^{(k)} &&&& \lambda_{2}^{(k)} && & &\cdots &&&& \lambda_{[k]}^{(k)} 
	&& && 0 \\
	\end{array}	
	\end{equation}
whereas for $k$ even, they are
	\begin{equation}\label{eq:interlacing-O-odd}
	\begin{array}{@{}c@{\;}c@{\;}c@{\;}c@{\;}c@{\;}c@{\;}c@{\;}c@{\;}c@{\;}c@{\;}c@{\;}c@{\;}c@{\;}c@{\;}c@{\;}c@{\;}c@{}c@{\;}c@{}c@{\;}c@{}}
 	\lambda_1^{(k+1)} &&&& \lambda_{2}^{(k+1)} &&&& \lambda_{3}^{(k+1)} &&&& \lambda_{[k+1]}^{(k+1)} &&&& 0\mathrlap.\\
	&\dgeq && \ugeq && \dgeq && \ugeq && \dgeq && \ugeq && \dgeq &&{\ugeq} \\
	&& \lambda_{1}^{(k)} &&&& \lambda_{2}^{(k)} && & &\cdots &&&& \big|\lambda_{[k]}^{(k)}\big| &&\\
	\end{array}	
	\end{equation}
In both cases, the pattern of inequalities can be extended 
symmetrically by including $0$ and the negatives of the $\lup k_j$,
representing $(1/i)$ times all eigenvalues of all truncations.\footnote{\ 
  One may omit the right side and the zero in the bottom row without loss of information, 
  but we opt to keep these features because they help simplify some descriptions later on.
}
This is thus again represented by a triangular array,
such as the following for $n = 4$.
\begin{equation}\label{eq:gen-4}
\begin{array}{@{}c@{\;}c@{\;}c@{\;}c@{\;}c@{\;}c@{\;}c@{\;}c@{\;}c@{\;}c@{\;}c@{\;}c@{\;}c@{\;}c@{\;}c@{\;}c@{\;}c@{\;}c@{\;}c@{\;}c@{\;}c@{}}
	\lambda_{1}^{(5)} && & & \lambda_{2}^{(5)} & & & & 0 & & & & -\lambda_{2}^{(5)} && & & -\lambda_{1}^{(5)} \\
	&\dgeq & & \ugeq & &\dgeq & & \ugeq & & \dgeq & & \ugeq&& \dgeq && \ugeq\\
	&& \lambda_{1}^{(4)} & & & & |\lambda_{2}^{(4)}| & & & & -|\lambda_{2}^{(4)}| &&&& - \lambda_{1}^{(4)}\\
	&&&\dgeq & & \ugeq & & \dgeq & & \ugeq & & \dgeq && \ugeq \\
	&&&& \lambda_{1}^{(3)} &&&& 0&& && -\lambda_{1}^{(3)} \\
	&&&&&\dgeq && \ugeq &&\dgeq && \ugeq \\
	&&&&&& |\lambda_{1}^{(2)}| &&&& -|\lambda_{1}^{(2)}| &\\
	&&&&&&&\dgeq && \ugeq && &&&&\\
	&&&&&&&& 0 &&&
\end{array}	
\end{equation}

The \defd{orthogonal Gelfand--Zeitlin (GZ) polytope} $\D$ 
is the set of points in $\smash{\R^{\sum_{k=1}^n [k]}}$
for which the coordinates $\lup k_j$ ($1 \leq j \leq [k] < k \leq n$)
satisfy all the interlacing inequalities. 
One passes to a codimension-one face every time a strict inequality
becomes an equality.


\defd{Orthogonal Gelfand--Zeitlin (GZ) patterns} 
are obtained from inequality triangles much as in the unitary case.
The resulting triangles are symmetric about the vertical midline.
For visual salience we colour vertices corresponding to $0$ white,
although we know vertices along the midline correspond to $0$ anyway.
The full subgraph on the black vertices to the left of the
vertical midline is the \defd{positive part} of the orthogonal GZ pattern
and that on the right the \defd{negative part}.
The same pattern obtains regardless of whether a given 
$\lup {2k}_k$ is positive or negative,
so that although the pattern of $p \in \D$ depends only on the face 
of $\D$ whose interior contains $p$,
the pattern no longer corresponds to just one face.
\end{discussion}

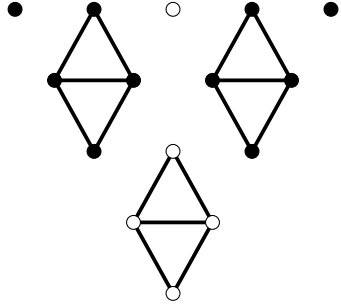
\begin{figure}
\centering
\begin{tikzpicture}[scale =.7,line cap=round,line join=round,>=triangle 45,x=1cm,y=1cm,scale=1.5, every node/.style={scale=1.5}]
	\begin{scriptsize}
		
	\draw [line width=0.5mm] (2-2/2,2*.9)-- (1-2/2,2*.9) -- (2-3/2,3*.9) -- (2-2/2,2*.9) -- (1-1/2,1*.9)--(1-2/2,2*.9);
	\draw [line width=0.5mm] (1,0)--(1+1/2,1*.9) --(2,0) -- (1,0);
	\draw [line width=0.5mm] (3-2/2,2*.9) -- (4-2/2,2*.9) -- (4-3/2,3*.9) -- (3-2/2,2*.9) -- (1+3/2,1*.9) -- (4-2/2,2*.9);
	\draw [line width=0.5mm] (1,0) -- (1.5,-1*.9) -- (2,0);
	
	\draw [fill=black] (1-3/2,3*.9) circle (2.5pt);
	\draw [fill=black] (2-3/2,3*.9) circle (2.5pt);
	\draw [fill=white] (3-3/2,3*.9) circle (2.5pt);
	\draw [fill=black] (4-3/2,3*.9) circle (2.5pt);
	\draw [fill=black] (5-3/2,3*.9) circle (2.5pt);

	\draw [fill=black] (1-2/2,2*.9) circle (2.5pt);
	\draw [fill=black] (2-2/2,2*.9) circle (2.5pt);
	\draw [fill=black] (3-2/2,2*.9) circle (2.5pt);
	\draw [fill=black] (4-2/2,2*.9) circle (2.5pt);
	
	\draw [fill=black] (1-1/2,1*.9) circle (2.5pt);
	\draw [fill=white] (1+1/2,1*.9) circle (2.5pt);
	\draw [fill=black] (1+3/2,1*.9) circle (2.5pt);
	
	\draw [fill=white] (1,0) circle (2.5pt);
	\draw [fill=white] (2,0) circle (2.5pt);

	\draw [fill=white] (1.5,-1*.9) circle (2.5pt);
	\end{scriptsize}
\end{tikzpicture}
\caption{An orthogonal GZ pattern associated to a fiber of a GZ system on a coadjoint orbit of $\SO(5)$.}
\label{fig:orthogonal-GZ-pattern}
\end{figure}

\begin{example}
  If in \eqref{eq:gen-4}
  the values are $\lambda_1^{(5)} = 2$,\,
$\lambda_2^{(5)} = \lambda_1^{(4)} = -\lambda_2^{(4)} = \lambda_1^{(3)} = 1$,\,
and $\lambda_1^{(2)} =0$, 
then the GZ pattern is given by \Cref{fig:orthogonal-GZ-pattern}. 
\end{example}

\begin{discussion}[Rows and the groups $H_k$]\label{rmk:orthogonal-row}
\defd{Row $k$} of an orthogonal GZ pattern is again the full subgraph on this 
row of vertices,
comprising some number of component segments.
Let $\defm{\ell}$ denote the number of connected components of the positive part,
each corresponding to a distinct absolute value of an entry of $\lup k$,
let $\defm{d_i}$ ($1 \leq i \leq \ell$) 
be the number of vertices in the $i\th$ component from the left of the positive part,
and let $\defm{d_0}$ be the
number of vertices in the white component.

Then the centralizer $H_k$ of $\lup k$ in $\SO(k)$
is again determined by the GZ pattern and seen to lie in the block sum
$\Direct_{i=1}^\ell \SO(2d_i) \+ \SO(d_0)$.
For $0 < i < \ell$ a block must commute with $\b(1)^{\+d_i}$,
forcing it to lie in  $\Sp(2d_i,\R) \inter \SO(2d_i)$,
which is the image of the embedding $\U(d_i) \lt \SO(2d_i)$
induced by $a + ib \lmt 
\big[\begin{smallmatrix}
	a&-b\\b&\phantom{-}a
\end{smallmatrix}\big]\:
\C \lt \R^{2 \x 2}$.
The same holds for the $\ell\th$ block if $k$ is odd
or $k = 2m$ is even with $\lup {2m} _m \neq -\lup{2m}_{m-1}$.
If however $k = 2m$ with $\lup{2m}_m = -\lup{2m}_{m-1}$,
then the element of $\ft_{2m,+}^*$ 
we are denoting $\lup{2m}$
ends in the block
$\b(\lup{2m}_{m-1})^{\+ d_\ell -1}  \+ \b(-\lup{2m}_{m-1})$,
which is the conjugate of
$\b(\lup{2m}_{m-1})^{\+ d_\ell}$
by the permutation matrix
 $\defm o = [1]^{\oplus d_{\ell} -2} \oplus 
 \big[\begin{smallmatrix}
	0&1\\1&0
 \end{smallmatrix}\big]$;
it follows the corresponding block of $H_k$
is the conjugate $\U(d_\ell)^o$.%
\footnote{\ 
 	Hence this block is not conjugate to $\U(d_\ell)$ in $\SO(k)$ 
	but becomes so in $\SO(k+1)$.
	This dichotomy will not be reflected in the quotients $H_k/L_k$,
	but will affect our description of $L_k$ and 
	of the horizontal maps $G_k/L_k \lt G_{k+1}/H_{k+1}$
	in \Cref{fig:big-diagram}.
}
If we write
$\defm{\U_{\ell}}$ for 
the $\ell\th$ block of $H_k$
in all cases, then we have
\begin{equation}\label{eq:J-O}
	H_k = 
	\U(d_1) \oplus \cdots \oplus \U(d_{\ell-1}) \oplus \U_{\ell} \oplus \SO(d_0)
\end{equation}
where $\SO(1) \iso \SO(0) \ceq \{1\}$
and $d_0$ is always congruent to $k$ modulo $2$.
\end{discussion}

\begin{discussion}[Representatives $\xip{k+1}$ and groups $L_k$]\label{rmk:xi-SO-explicit}
  The range of options for $\xip{k+1}$ in \Cref{rmk:stabilizer-bundle}
	for the orthogonal case is a variant of \eqref{eq:Y}, 
	with $\b(\mu_i)^{\+ d_i}$ replacing $\diag(\mu_i,\ldots,\mu_i)$
	for $0 < i < \ell$, 
	the $o$-conjugate replacing it for $i = \ell$,
	and an extra block of zeroes corresponding to 
	the $\SO(d_0)$ factor of $H_k$.
	The $z_i$ are replaced by arbitrary vectors in $\R^{2d_i}$
	or $\R^{d_0}$ of fixed norm $r_i$,
	which continue to be continuous in $p \in \D$
	and $c$ by $0$.
	We continue to take $z_i = [r_i \ 0\ \cdots\ 0]^\top$
	when possible
	and $z'_\ell = [0\ \cdots\ 0\ r_\ell]^\top$ when necessary.
	The standard choice of $\xip{k+1}$
	has stabilizer
 	\begin{equation}\label{eq:K-O} L_k = K_1 \oplus \cdots \oplus
			K_{\ell_k} \oplus K_{0} \leq
   \U(d_1) \oplus \dots \oplus \U_{\ell} \+ \SO(d_0) =  H_k\mathrlap,
	\end{equation}
	where the blocks are given as follows,
	with the conventions $\SO(0)  \ceq \{1\} \eqc \U(0)$.
\bitem
	\item
	For $1\leq i < \ell$, and $i = \ell$ with $\U_\ell = \U(d_\ell)$,
	the $K_i$ are as in the unitary case.
	\item
		If
		$\U_\ell = \U(d_\ell)^o$,
		then  $K_\ell$ is $\U_\ell$ 
		for a \parallelogram-shape and
		$\big([1] \oplus \U(d_\ell - 1)\big)^o$ 
		for a \mtrap-shape.
%
	\item
	  	The group
		${K_{0}}$ is $\SO(d_0)$
		for a 
		\wtrap- or \parallelogram-shape
		and
		$ [1] \+ \SO(d_0-1)$
		for a \mtrap-shape.
  \end{itemize}
  The alternative choice $'\xipko$ is again related to 
  $\xipko$ by $'\xipko = \vp_k(h)\xipko$
  for some $h = \Direct_{i=1}^\ell u_i \+ u_0
  \in \Direct_{i=1}^{\ell-1} \SU(d_i) \+ \mr{S}\U_\ell \+ \SO( d_0)$,
  where $\mr{S}\U_\ell$ is either $\SU(d_\ell)$ or $\SU(d_\ell)^o$.
  \end{discussion}


\begin{example}
Consider the orthogonal GZ pattern in Figure~\ref{fig:orthogonal-GZ-pattern}.
From rows $4$ and $5$ we see that $H_4 = \U(2) \times \U(2)$ 
and $L_4 = \big([1]\oplus\U(1)\big)\times \big([1]\oplus\U(1)\big)$,
from rows $3$ and $4$ we see that $H_3 = L_3 = \U(1) \times \U(1)$, and
from rows $2$ and $3$ we see that $H_2 = \SO(2)$ and $L_2 = \{1\}$.
\end{example}

\section{The structure of Gelfand--Zeitlin fibers}\label{sec:product}

We begin in \Cref{sec:coordinates} by 
identifying a system of coordinates on $\Psi\-(p)$ in $\prod G_k/L_k$
and connecting this to the staircase \Cref{fig:big-diagram}.
Unpacking this in \Cref{sec:modification} allows us to identify
the horizontal maps in the diagram abstractly and to change them without altering
the diffeomorphsim type of the final pullback.
In \Cref{sec:horizontal-explicit}, 
we explicitly identify these maps in the cases of unitary and orthogonal \GZ systems.
In \Cref{sec:balanced-product}, 
we rewrite the bundle towers of \Cref{rmk:stabilizer-bundle}
as balanced products
(Theorem~\ref{thm:explicit-description});
these are equivalent to biquotients (\Cref{thm:biquotient}).
A direct product decomposition of GZ fibers indexed by pattern components,
follows immediately from this description (\Cref{thm:component-product}). 
Moreover, these factors generically simplify substantially 
(\Cref{rmk:tensor-discussion}, \Cref{rmk:pinch}, \Cref{rmk:opportunism}).
In \Cref{sec:circles},
we derive a toral direct factor 
for both unitary and orthogonal GZ fibers 
in an elementary manner from the balanced product description
and show the resulting torus action is the same
given by the Thimm trick (\Cref{thm:Thimm}).
Interspersed throughout,
we give a number of examples demonstrating these results.
In particular, we recover Bouloc, Miranda, and Zung's
description of regular, elliptic, and (multi-)diamond fibers of unitary GZ
systems as special cases.

\subsection{Coordinate systems and transition functions}\label{sec:coordinates}

In this subsection,
we are again in the general set-up of a \TGZS $\Psi$,
as in \Cref{def:GZ-abstract} and \Cref{rmk:stabilizer-bundle},
trying again to characterize when a list of elements
of $\prod_{k=\a}^n {\Psi\mn_k}\-(\lup k)$ 
is the list of iterated $\Phi$-images of some $\genup{n+1} \in \Psi\-(p)$.

\begin{discussion}[The staircase, $G_k/L_k$ coordinates, and conjugating elements $a_k$]\label{rmk:a}


  The staircase of \Cref{fig:big-diagram},
  interpreted through the lens of \Cref{rmk:stabilizer-bundle},
  represents $\Psi\-(p) \iso F^{n+1}_\a$ 
  as the iterated fiber product
  \quation{\label{eq:Omega}
      \vphantom{X_{X_{X_{X_{X_X}}}}}
      \xu{{G_{\alpha+1}/H_{\alpha+1}}} {G_\alpha/L_\alpha }
      { G_{\alpha+1}/L_{\alpha+1} }
      \xu{G_{\alpha+2}/H_{\alpha+2}}{}
      {\cdots}
      \xu{G_n/H_n}{}{ G_n/L_n}\mathrlap,
  }
a   subset of the direct product $ \prod_{k=\a}^n G_k/L_k$
  comprising certain lists $(g_k L_k) $.

  What are these $g_k$?
Consider a point $\genup{n+1} \in \Psi\-(p)$
and recall that we fixed a choice of representatives $\xip{k+1} \in \Orb_{\lup{k+1}} \inter \Phi\-(\lup k) = \Orb_{\lup{k+1},\lup k}$.
The transitivity condition
along with the assumption that $\genup{n+1}$ lies in $F^{n+1}_n = \Orb_{\lup{n+1},\lup n}$
implies that $\vp(g_n)\xip{n+1} = \genup{n+1}$ for some $g_n \in G_n$.
Similarly, that $\Phi(\genup{ n+1})$ lies in $F^{n}_{n-1} = \Orb_{\lup n, \lup{n-1}}$
implies there is some $g_{n-1} \in G_{n-1}$ with $\vp(g_{n-1})\xip n = \Phi(\genup{ n+1})$ 
and so on, so that $(g_k L_k)$ 
corresponds to the list $\big(\vp(g_{k})\xip {k+1}\big) = \big(\Phi^{n+1-k}(\genup{n+1})\big)$
 of iterated truncations of $\genup{n+1}$.

 Which lists occur is determined by
a set of coherence conditions 
  involving the projections $G_k/L_k \longepi G_k/H_k$
  and the maps \[G_k/L_k
\simto \vp(G_k)\. \xi_{k+1}		
\longinc
	 G_{k+1} \. \lup{k+1}
 \simto G_{k+1}/H_{k+1}\] 
 we now characterize.
\end{discussion}

\begin{lemma}\label{thm:iotak}
	Make identifications 
	$F^k_k = \mathcal{O}_{\lambda^{(k+1)},\lambda^{(k)}}\cong G_k/L_k$ 
	and $F^{k+1}_k = \mathcal{O}_{\lambda^{(k+1)}}\cong G_{k+1}/H_{k+1}$ 
	by fixing $\xip{k+1}\in \mathcal{O}_{\lambda^{(k+1)},\lambda^{(k)}}$ 
	and $\lambda^{(k+1)}\in \mathcal{O}_{\lambda^{(k+1)}}$ 
	as in \Cref{rmk:stabilizer-bundle},
	and let $a_{k+1} \in G_{k+1}$ be such that $\l^{(k+1)} = 
	a_{k+1}\xip{k+1}$. 
	Then $F^{k+1}_k \lt F^{k+1}_{k+1}$ is identified with 
	the map $G_k/L_k \lt G_{k+1}/H_{k+1}$
	given by  
\quation{\label{eq:iotaa}
\defm{\iota_{a_{k+1},1}} \: gL_k \longmapsto 
\varphi_k(g) a_{k+1}\- H_{k+1}\mathrlap. 
}
\end{lemma}
\begin{discussion}\label{rmk:compatibility-conditions}
The compatibility conditions for a point $(g_k L_k)$ to lie in the fiber product
are then
\quation{\label{eq:compat}
\phantom{ , \qquad (\a \leq k \leq n-1)} 
	\vp(g_k)a_{k+1}\-H_{k+1} = g_{k+1} H_{k+1}
	\qquad (\a \leq k \leq n-1)\mathrlap.
}
From a point in the fiber product \eqref{eq:Omega},
the element  $\genup{n+1} \in \Orb$ can be recovered
again as $\vp(g_n)\xip{n+1} = \vp(g_n)a_{n+1}\-\lup{n+1}$.\footnote{\ 
Since the horizontal maps
in \Cref{fig:big-diagram} are all injections,
the coordinates $g_k L_k$ for $k < n$
in some sense merely ensure that $g_n$ satisfies the compatibility conditions.
}
\end{discussion}


\begin{discussion}[The effects of changing $\xi$]\label{rmk:change-of-xi}
We fixed the representatives $\lup k$ 
 and groups  $H_k = \Stab_{G_k} \l^{(k)}$
in \Cref{rmk:lambda-staircase} permanently,
and chose $\xip{k+1}$  
in \Cref{rmk:stabilizer-bundle}
provisionally,
with the warning in \Cref{rmk:xi-prime} 
that we would later need to consider different choices.
If change 
$\xip{k+1}$ to $'\xip{k+1} = \vp_k(h)\xip{k+1}$ for some $h \in H_k/L_k$,
fixing the other representatives,
then in 
\eqref{eq:Omega} and \eqref{eq:compat}
we get 
\[
	L'_k = hL_k h\-,\qquad\qquad
	g'_k = g_k h\-,\qquad\qquad
	a'_{k+1} = a_{k+1} \vp_k(h)\-\mathrlap,
\]
the other $L_j$, $g_j$, and $a_j$ being unchanged.
\end{discussion}

\subsection{A new pullback}\label{sec:modification}



%

\begin{construction}\label{rmk:modification}
  We might consider changing $\iota_{a,1}$ in \eqref{eq:iotaa} to
\[
\defm{\iota_{a,b}}\: gL_k \longmapsto  {b}\varphi_k(g)a\- H_k
\]
for some $b \in G_{k+1}$.
Since
$\iota_{a,{b}}$ is equal to the composition of $\iota_{a,1}$
with left translation by ${b}$ on $G_{k+1}/H_{k+1}$,
which is an automorphism 
of the bundle $G_{k+1}/L_{k+1} \lt G_{k+1}/H_{k+1}$,
the pullbacks of this bundle by $\iota_{a,{b}}$ and $\iota_{a,1}$ are diffeomorphic.
It follows that the resulting final pullbacks~$F_\alpha^{n+1} \iso \Psi\-(p)$
in \Cref{fig:big-diagram} will be diffeomorphic if we make such a change
for any number of $k$.
In particular, we would like to instead take ${b}_{k+1} = a_{k+1} $
for all $k$,
allowing us to reanalyze
the horizontal maps in \Cref{fig:big-diagram}
as composites
\quation{\label{eq:conj-display}
    \frac{G_k}{L_k}				\isoto
    \frac{ a \vp(G_k) a\-} {a \vp(L_k) a\-} 		\lt
    \frac{G_{k+1}}{H_{k+1}}			\mathrlap,
}
the first factor being induced by the map $c_a \o \vp\: G_k \lt G_{k+1}$
taking $g$ to $a\vp(g)a\-$
and the second by the inclusion.
\end{construction}

\begin{notation}\label{def:Omega}
  We will write 
  $\defm{\W^r_p}$ for $F_\a^{n+1}$ in \eqref{eq:Omega}
  under the maps $g_k L_k \lmt \vp_k(g_k)a_{k+1}\-H_{k+1}$
  and~$\defm{\W^c_p}$ for the (diffeomorphic) pullback 
 under the maps $g_k L_k \lmt a_{k+1}\vp_k(g_k)a_{k+1}\-H_{k+1}$.
\end{notation}

\begin{discussion}[Coordinate conversion from $\W^r_p$ to $\W^c_p$]\label{rmk:explicit-modification}
When the maps $G_k/L_k \lt G_{k+1}/H_{k+1}$
are altered
per \Cref{rmk:modification},
the conditions on coset lists $(\hat g_k L_k) \in \prod G_k/L_k$
from \Cref{rmk:compatibility-conditions}
become instead
  \[
  	\qquad 		\phantom{(\alpha \leq k \leq n-1)}
	a_{k+1}	\vp_k(\hat g_k)a_{k+1}\-H_{k+1} = \hat g_{k+1}H_{k+1}
	\qquad 		(\alpha \leq k \leq n-1)
	\mathrlap.
  \]
For any choice of $a_\a \in G_\a$ (of course, $a_\a = 1$ is most natural),
and omitting the $\vp_k$ from the notation to make the pattern clearer,
an explicit diffeomorphism is $(g_k L_k) \lmt (\hat g_k L_k) $ with
  \[
      \hat g_\a 		 = a_\a g_\a ,\qquad
      \hat g_{\a+1} 	 = a_{\a+1} a_\a \.g_{\a+1},\qquad
      \cdots,\qquad
      \hat g_{n}	= a_{n} 	\cdots a_\a \.g_{n}.
  \]
\end{discussion}

\begin{discussion}\label{rmk:latitude}
We actually have more latitude 
to adjust the horizontal maps in \Cref{fig:big-diagram}
without changing the diffeomorphism type of $F_\a^{n+1}$.
Indeed,
given any automorphism $\s$ of $G_{k+1}$ taking $H_{k+1}$ to itself,
we may exchange the map $\vp_k$ for $\sigma \o \vp_k$,
changing the map 
of \eqref{eq:conj-display}
to $gL_k \lmt \s\big({a}\varphi_k(g){a}\-)H_{k+1}$
and leaving the diffeomorphism type of the fiber unchanged.
\end{discussion}



\subsection{Continuity of transition functions}\label{sec:horizontal-explicit}

Now that we have a description of the horizontal maps 
in \Cref{fig:big-diagram},
we are able to say what the elements $a_{k+1} \in G_k$
are for unitary and orthogonal \GZ fibers
and as a result to describe the maps $c_{a_{k+1}} \o \vp_k$
figuring in \Cref{thm:iotak}.
It will be important to us later that,
like our choice of elements $\genup{k+1}$,
the $a_{k+1}$ can locally be chosen continuous in $p$,
at least for unitary and orthogonal \GZ systems.

\begin{discussion}\label{rmk:a-explicit}
For a unitary GZ system,
we claim a natural choice of matrices $a = a_{k+1}$ conjugating $\xi = \xip{k+1}$ to $\l^{(k+1)}$
	is continuous.
	Note that the assumed equation  $\l^{(k+1)} = \Ad^*_a \xi = a \xi a\-$,
	rearranged to $\xi a\- = a\- \l^{(k+1)}$,
	means precisely that the $j\th$ column of $a\- = \bar a^\top$ lies in the
	$\lup{k+1}_j$-eigenspace of $\xi$.
	Hence, to define~$a$,
	it will be enough to find an orthonormal frame of such eigenvectors.
	Every instance of an eigenvalue of $\l^{(k)}$
	also occuring as an eigenvalue of $\lup{k+1}$,
	say $\smash{\lup{k+1}_j}$,
	corresponds to a 
	 $j\th$ row of the ``boring'' form
	 $\smash{\lup{k+1}_j}\.(\d^j)^\top
	 = [0\, \cdots \, 0 \, \smash{\lup{k+1}_j} \, 0 \, \cdots \, 0]^\top$
	in $\xi$,
	so we may take the column vector $\d^j$ 
	the $j\th$ column of $a\-$.
	If the corresponding component in the $k\th$ and $(k+1)\st$ rows of the GZ pattern
	is a \parallelogram- or \mtrap-shape, 
	then we have found a basis for the corresponding eigenspace.
	Otherwise, in the case of a \wtrap-shape, we have selected all but one basis eigenvector
	and the orthogonal complement to the others within the $\lup{k}_j$-eigenspace of $\xi$
	is $1$-dimensional, and so contains precisely a circle of candidate eigenvectors.
	From our convention for choosing the $z_i$ in \eqref{eq:Y}
	to have only first coordinate nonzero,
	it will be 
	the first column vector in $\a\-$ (counting from the left)
	in the $\lup{k}_j$-eigenspace of $\xi$
	that has this non-$\d^j$ form.

	Obviously the $\d^j$ are continuous in $\xi$ and orthonormal,
	but it could \emph{a priori} be that the new
	$w\th$ columns $\defm{v_w} = [v^1_w\, \cdots\, v^{k+1}_w]^\top$,
	which we only specified up to a constant in $\U(1)$,
	cannot be continously defined.
	To see they can be, note that the entries $v^j_w$ for which 
	the $j\th$ row of $\xi$ is $\smash{\lup{k+1}_j}\.\d^j$ 
	must be $0$ for $v_w$ to be orthogonal to the columns $v_j = \d^j$ we have already chosen.
	For the remaining $j$,
	the $i\th$ row of $\xi$ has the ``interesting'' form
	$[0\,\cdots\,0\,\mu_i\,0\,\cdots\,0\,r_i]$.
	Multiplying the $j\th$ row of $\xi$ 
	with the $w\th$ column $v_w$ of $a\-$ we are determining,
	which we assume is a $\lup{k+1}_w$-eigenvector of $\xi$,
	we find the constraint $\smash{v^j_w } = r_i(\lup{k+1}_w-\mu_i)\- \. v^{k+1}_w$,
	meaning the other nonzero entries are uniquely determined as real multiples of the 
	bottom entry $v^{k+1}_w$. 
	Temporarily set this to $1$ and call the resulting vector $v'_w$.
	Since $R = 1^2 + \sum_m (\sfrac{r_i}{(\l^{(k+1)}_w-\mu_m)})^2$ is a continuous function of~$\xi$ 
	(here the sum is over interesting rows of $\xi$, corresponding to \mtrap-shapes),
	the vector $v_w = \frac 1{\sqrt R} v'_w$ is a continuous function of $\xi$ as well.
	The matrix $a'$ is actually orthogonal, having only real entries.
	In case $\det a' = -1$, we can get a special orthogonal matrix $b$ instead
	 by substituting the last column $v_{k+1}$ with $-v_{k+1}$.
	Finally, set $a = b\-$.

	We will later need to show that the elements $\xip{k+1}$ 
	and hence the elements $a_{k+1}$
	can also be chosen so as to be continuous over the union $E \union F$
	of the interior $E \sub \D$ of a face
	and that  of an adjoining codimension-one face.
To see this, 
we work matrix coefficient by matrix coefficient in \eqref{eq:Y}.
Most entries of $\xip{k+1}$ are simply $0$, which is surely continuous.
Along the first $k$ entries of the diagonal,
we have convergence by assumption,
since these entries are $\lup{k}(q)$.
The last diagonal entry $c(q) = \tr \lup{k+1}(q)-\tr\lup{k}(q)$
is hence also continuous.

The remaining entries, $r_i$ in the discussion in \Cref{rmk:xi-explicit},
yield to a case analysis.
For $r_i$ corresponding to a \mtrap-shape not labeled by $\mu^\pm$,
the expression \eqref{eq:r} has $ r_i(q) \to r_i(p)$ as $q \to p$
since we have $\mu^\pm \to \mu$ and no factor approaching $0$ appears in denominator.
For $r_i$ corresponding to a \mtrap-shape labeled by $\mu^\pm$,
we subdivide by case.
For \rpara\rpara\ and \rpara\wtrap,
and symmetrically \parallelogram\parallelogram\ and \wtrap\parallelogram, 
there are no such $r_i$.
For \mtrap\rpara, $\mu^+$ does not appear in \eqref{eq:r}
and $r_i(q) \to r_i(p)$ since $\mu^- \to \mu$.
For \parallelogram\mtrap, symmetrically, we use $'\xip{k+1}$ in place of $\xip{k+1}$
and again $r_i \to 0$,
this time since $\mu^+ \to \mu$.
For \mtrap\wtrap\, the square $r_i(q)^2$ 
factors as a rational function of $\mu^-$ alone
and the factor $\mu^- - \mu^+ \to 0$,
so indeed $r_i(q) 
\to 
0 = r_i(p)$.
For \wtrap\mtrap, we take $'\xip{k+1}(q)$ as noted above,
and the discussion of $r_i$ is the same as for \mtrap\wtrap.

The description and continuity of the $\xipko $
are analogous in the  orthogonal case, 
but more notationally involved.
\end{discussion}

\begin{discussion}\label{rmk:horizontal-unitary}
In the unitary case,
the map $c_a \o \vp\: L_k \lt H_{k+1}$ of \eqref{eq:conj-display}
	sends the factor $K_j$ of $L_k$ 
	associated to the $j\th$ connected 
	component of the $k\th$ row to the 
	unitary factor $\U(e)$ of $H_{k+1}$ corresponding to the same connected component of
	the GZ pattern.
	In the case of a \mtrap- or \parallelogram-shape 
	($\smash{e = \rk K_j} = d_j-1$ or $d_j$, respectively)
	the composite map $K_j \inc \smash{L_k \os{c_a}\to H_{k+1}} \epi \U(e)$
	is the natural isomorphism, respectively 
	$[1] \+ \U(e) \simto \U(e)$
	or 
	$\U(e) \eqto \U(e)$.
	In the case of a \wtrap-pattern ($e = d_j+1$),
	this composition is the block-inclusion $[1] \+ \U(d_j) \inc \U(1+d_j)$.
	Note that although the element $a$ and hence the automorphism $c_a \in \Aut G_{k+1}$
	vary with $p$, the map $c_a \o \vp\: L_k \lt H_{k+1}$ as just described 
	is the same for any two points in the interior of the same face of $\D$.
\end{discussion}

\begin{discussion}\label{rmk:horizontal-orthogonal}
Again in the orthogonal case,
the map $c_a \o \vp\: L_k \lt H_{k+1}$ of \eqref{eq:conj-display}
	sends the factor $K_j$ of $L_k$ 
	associated to the $j\th$ connected 
	component of the $k\th$ row to the 
	factor of $H_{k+1}$ corresponding to the same connected component of
	the GZ pattern.
	For the white component, these factor inclusions 
	$\SO(d_0) \inc L_k \os{c_a}\to H_{k+1} \epi \SO(e)$
	are as expected by analogy with the unitary case,
	either the identity map in the case of a \mtrap-
	or \parallelogram-shape,
	and $\SO(d_0) \simto [1] \+ \SO(d_0) \inc \SO(1+d_0)$
	in the case of a \wtrap-shape.

	The numerology of the inclusions
	$\U(d) \inc L_k \mono H_{k+1} \epi \U(e)$ of unitary
	factors indexed by the same $\l$ is the same as in the unitary case described in
	\Cref{rmk:horizontal-unitary},
	with the same dependence on \mtrap-, \parallelogram-, or \wtrap-shapes
	occuring in the GZ pattern, 
	but the embeddings themselves may be different.
	The embeddings between \emph{standardly embedded} unitary factors
	are the standard embeddings, as in that section.

	By contrast, the embedding of 
	a nonstandard unitary factor $K_\ell$ in a standard 
	unitary factor $\U(e)$ of $H_{k+1}$
	is given by
	$\smash{o\-\U(e)o \os{c_o}\lt \U(e)}$
	in the cases of both a \parallelogram-shape and a \mtrap-shape
	(and a \wtrap-shape cannot occur).
	The embedding of a standard $K_\ell$ in a nonstandard 
	factor of $H_{k+1}$ (necessarily $K_\ell = \U(d_\ell)$
	in this case since a \mtrap-shape cannot occur)
	is given by $\U(d_\ell) \os{c_o}\lt \U(d_\ell)$
	for a \parallelogram-shape,
	but for a \wtrap-shape is given by
	$\U(d_\ell) \simto \U(d_\ell) \+ [1] \inc \U(d_\ell + 1)^o$.
	Note that unlike in all other cases,
	this is an \emph{upper-left} rather than lower-right
	block inclusion, and that the image is invariant
	under conjugation by $o$, which is equal to the identity
	matrix except in the lower-right $2\x 2$ block.
	To deal with this case, we employ \Cref{rmk:latitude}
	to modify the map $G_k/L_k \lt G_{k+1}/H_{k+1}$
	by an automorphism $\s$ of $G_k$ preserving $H_k$.
	In this case, we let $\s$ be conjugation by $o\-Po$,
	where $P$ is the matrix 
	inducing the permutation exchanging the upper-left $2d_\ell \x 2d_\ell$
	with the lower-right $2 \x 2$ submatrix.
	The factor in the replaced $H_{k+1}$ is again $\U(d_\ell+1)^o$,
	and the injection from $K_\ell = \U(d_\ell)$
	is now 
	$\U(d_\ell) \simto [1] \+ \U(d_\ell) \inc \U(1+d_\ell)$
	followed by $c_o$.
	
	Note again that the restriction of the maps $c_a \o \vp$ or $\s \o c_a \o \vp$
	to $L_k$
	is the same for every $p$ in the interior of a face of $\D$
	even though the conjugating elements $a$ vary from point to point.
\end{discussion}

\subsection{Iterated balanced products and biquotients}\label{sec:balanced-product}

It is in this section we undergo the perspective shift that 
we have claimed will simplify our lives.
We begin with a straightforward lemma. 

\bdefn
Suppose a group $G$ acts on the right on a space $X$
and on the left on a space $Y$.
In this case, we define the \defd{balanced product} $\defm{X \ox_G Y}$ 
to be the quotient $(X \x Y)/G$ under the diagonal action. 
The elements of ${X \ox_G Y}$, which we denote $\defm{x \ox y}$,
satisfy $xg \ox y = x \ox gy$.
\edefn

\begin{lemma}\label{thm:pullback-tensor}
Let $\vp\: G \lt \wt G$ be a map of topological groups
and let $L \leq G$ and $\wt L \leq \wt H \leq \wt G$ be closed subgroups such that $\vp(L) \leq \wt H$.
Suppose $\wt L$ acts continuously on a topological space $X$.
Then the map 
\eqn{
\vk\:\xu{\wt G/\wt H} {G/L}{\wt G} \ox_{\wt L} X 
	&\lt 
G \ox_L \wt H \ox_{\wt L} X,\\
(gL,\wt g \ox x) 
	&\lmt
g \ox \vp(g)\- \wt g \ox x\mathrlap,
}
is a well-defined homeomorphism.
If the homomorphisms are of compact Lie groups 
and $X$ is a smooth manifold with
smooth $\wt L$-action, then $\vk$ is a diffeomorphism.
\end{lemma}
\begin{proof}
Since the fiber product ${G \x_{\smash{\wt G/\wt H}} \wt G}$
is the set of pairs $(g,\wt g) \in G \x \wt G$ 
such that $\vp(g)\-\wt g$ lies in ${\wt H}$,
the assignments
$(g,\wt g) \lmt \big(g,\vp(g)\-\smash{\wt g}\big)$
and $(g,\wt h) \lmt \big(g,\vp(g)\wt h\big)$
define a continuous map
$G \x_{\smash{\wt G/\wt H}} \smash{\wt G} \lt G \x \smash{\wt H}$
and its continuous inverse.
One checks using coset representatives that these 
descend to well-defined maps
$\kappa\: G/L \x_{\smash{\wt G/\wt H}} \smash{\wt G} \lt G \ox_L \smash{\wt H}$ 
and $\kappa\-$,
necessarily homeomorphisms,
and then that 
$\kappa \x \id_X$ and 
$\kappa\- \x \id_X$ 
descend to a well-defined~$\vk$ and~$\vk\-$. 
Under the smoothness hypotheses,
the original homeomorphism and its inverse are smooth,
and as all the actions in question are free and proper,
so the orbit spaces are again smooth manifolds
and the quotient maps are smooth submersions,
meaning~$\vk$ and~$\vk\-$ are smooth as well.
\end{proof}

Here, then, is the perspective shift.

\begin{theorem}\label{thm:explicit-description}
Let $\Psi$ be a \TGZS, $p \in \D$ a point,
$L_k \leq H_k \leq G_k$ 
the groups discussed in \Cref{rmk:stabilizer-bundle},
and $a_{k+1}$ for $\a +1 \leq k \leq n+1$ elements as discussed in \Cref{rmk:a}.
Then there is a diffeomorphism between $\Psi^{-1}(p)$ and
\[
\defm{T_p}
	\ceq
G_\alpha \,\ox_{L_\alpha}\, H_{\alpha+1} \,\ox_{L_{\alpha+1}}\, \cdots \,\ox_{L_{n-1}}\, H_{n}/L_{n}
\]
where the balanced products with respect to the groups $L_k$
are taken with respect to the natural inclusion of $L_k$ into $H_k$
and the embedding of $L_k$ into $H_{k+1}$ via the homomorphism 
$c_{a_{k+1}} \o \vp_k$.
\end{theorem}

In both the orthogonal and unitary cases, we have $G_\a = H_\a$,
so that in fact $G_k$ does not figure.

\begin{proof}
We have a diffeomorphism $\Psi\-(p) \isoto \W^r_p$ 
by \Cref{rmk:lambda-staircase}
and a diffeomorphism $\W^r_p \lt \W^c_p$ 
by \Cref{rmk:modification}.
Consider~\ref{fig:big-diagram},
with the horizontal maps induced by $c_{a_{k+1}}$ 
as discussed in \Cref{rmk:modification},
so that $\W_p^c = F_\a^{n+1}$.
We will prove diffeomorphisms
\[
\phantom{(\a \leq j < k)\quad } F_j^{k} \cong G_j \,\ox_{L_j}\, H_{j+1} \,\ox_{L_{j+1}}\,
\cdots \,\ox_{L_{k-2}}\, H_{k-1}\,\ox_{L_{k-1}}\,\ast , \quad (\a \leq j < k \leq n+1)
\]
so the result will follow on taking $j = \a$ and $k = n+1$.
For any fixed vertical index $k$, 
the claim is proved 
by decreasing (i.e., leftward) induction on the horizontal index $j$.
For the base case,
when $j = k -1$, 
the expression simplifies to $F_{k-1}^{k} \cong G_{k-1}/L_{k-1}$. 
For the induction step, assume the diffeomorphism for $j$
and consider the tall rectangle $\smash{F_{[j-1,j]}^{[j,k]}}$
comprising $k-j$ vertically stacked boxes in \eqref{fig:big-diagram}.
We have $\smash{F^j_{j-1} = G_{j-1}/L_{j-1}}$ 
and $\smash{F^j_j = G_j/H_j}$,
and by induction
an expression $G_j \ox_{L_j} X$ for~$\smash{F_j^k}$.
For this $X$ and $(G,L) = (G_{j-1},L_{j-1})$ and 
$\smash{(\wt G,\wt H, \wt L)} = (G_j,H_j,L_j)$,
an application of \Cref{thm:pullback-tensor}
completes the induction.
\end{proof}

%

Since the groups $H_k$ and $L_k$ depend only on the face of $\D$ containing $p$ in its interior,
we recover in arguably more explicit form Cho--Kim--Oh's  observation 
that the diffeomorphism type of a \GZ fiber $\Psi\-(p)$
is determined entirely by the open face of $\D$ containing $p$.
\begin{discussion}[Coordinate conversion from $\W^c_p$ to $T_p$]\label{rmk:explicit-conversion}
An inductive calculation 
shows the diffeomorphism $\W^c_p \lt T_p$ 
takes $ (\hat g_k L_k)$ to
\quation{\label{eq:hat-tensor-conversion}
\hat g_\a \otimes (c_{a_{\a+1}} \o \vp_\a)(\hat g_\a)\- \.\hat g_{\a+1} 	
	\otimes (c_{a_{\a+2}} \o \vp_{\a+1})(\hat g_{\a+1})\- \.\hat g_{\a+2} 	
	\otimes \cdots
  	\otimes (c_{a_n} \o \vp_{n-1})(\hat g_{n-1})\- \.\hat g_{n} L_n
	.\ 
	}
	It follows that starting with a point $\vp_n(g_n)\xip{n} = \vp_n(g_n)a_{n+1}\- \l^{(n+1)}$
of $\mc O_{\l^{(n+1)}}$,
converting it to the point $(g_k L_k) \in \W^r_p$,
converting this in turn to the point $(\hat g_k L_k) \in \W^c_p$ 
via \Cref{rmk:explicit-modification},
and then finally converting to a point of $T_p$,
one recovers the point
\quation{\label{eq:tensor-conversion}
g_\a 	\otimes a_{\a+1}\vp_\a(g_\a)\- g_{\a+1}
	\otimes a_{\a+2}\vp_{\a+1}(g_{\a+1})\- g_{\a+2}
	\otimes \cdots
	\otimes a_n \vp_{n-1}(g_{n-1})\- g_n L_n 
	\in T_p
	\mathrlap.
}

Although the elements $a_{k+1}(p) \in G_{k+1}$ continuously vary, 
in \Cref{thm:only-standard}, 
we fixed $p$ and then explicitly identified in \Cref{rmk:horizontal-unitary} 
and \Cref{rmk:horizontal-orthogonal}
the effect of the conjugation maps $L_k \lt H_{k+1}$
on block components,
and found that for our standard choice of $\xip{k+1}$ they were block inclusions
of a standard form.
Since we had fixed $p$ before this discussion
found that these injections were of this form,
they are actually independent of $p$.
It follows that the description of $T_p$ in \Cref{thm:explicit-description}
does not depend on the choice of representatives $\xip{k+1}$
and within the interior of any given face of $\D$,
does not depend on $p$ either.
The same follows for the map $\W^c_p \lt T_p$ of \eqref{eq:hat-tensor-conversion}.
We will discuss in \Cref{sec:ray} 
how much this changes in transitioning betwen faces.
\end{discussion}

\begin{discussion}[The effects on $T_p$ of changing $\xi$]\label{rmk:tensor-change-of-xi}
On the other hand, the diffeomorphism $\W^r_p \lt T_p$
\emph{does} depend on $p$ because the values $a_k$ appear once each
in \eqref{eq:tensor-conversion}.
Indeed, in \Cref{rmk:change-of-xi},
we discussed the effects of a change of $\xi_{k+1}$
to $\vp(h)\xi_{k+1}$ in the description of $\W^r_p$.
Applying \Cref{rmk:explicit-conversion},
owing to cancellation of the new $h$'s in $g'_k$ and $a'_{k+1}$,
in the balanced product representation,
the point in \eqref{eq:tensor-conversion} 
then becomes
\quation{\label{eq:tensor-h-conversion}
g_\a \otimes \cdots
\otimes a_k \vp_{k-1}(g_{k-1})\- g_k \textcolor{red}{ h\- }
	\otimes \cdots
	\otimes a_n \vp_{n-1}(g_{n-1})\- g_n L_n \in T_p\mathrlap,
}
with only the displayed $k\th$ tensor-coordinate changed.
In particular, the balanced product presentation of \Cref{thm:explicit-description}
does depend on the choice of orbit representatives $\xip{k+1}$.
\end{discussion}

\begin{remark}
 Interestingly, Bouloc--Miranda--Zung~\cite[Thm.~4.16]{BoulocMirandaZung2018}
 obtain a different description of a GZ fiber as a biquotient.
We know only special cases where their expressions and ours 
simplify to the same thing (\Cref{eg:diamonds,eg:tori}).
 \end{remark}

\begin{remark}
Another instance of the balanced iterated product construction 
is the \emph{generalized Bott towers}
of Kaji, Kuroki, Lee, and Suh~\cite[Defs.~3.1, 3.5]{KKLS}.
Inspecting cohomology, one sees a generalized Bott tower 
is almost never a GZ fiber.
\end{remark}

\brmk\label{rmk:biquotient}
This sort of iterated balanced product 
can be expressed perhaps less usefully as a two-sided quotient.
\ermk

\bdefn
When $A$ and $B$ are closed subgroups of a Lie group $G$
such that the action $(a,b)\.g \ceq agb\-$ of $A \x B$ on $G$ is free,
the orbit space is called a \defd{biquotient} and denoted $\defm{A\backslash G/B}$.
\edefn

\begin{proposition}[Biquotient description]\label{thm:biquotient}
    The expression $T_p$ of \Cref{thm:explicit-description}
    is diffeomorphic to 
 a biquotient of Lie groups.
    Explicitly, let
	\begin{equation}
	\begin{split}
	  G & = \textstyle\prod_{k=1}^n H_k,
		\\
		A & = \big\{ (1,\ell_2,\ell_2,\ell_4,\ell_4, \ldots ) 
			\in G \colon \ell_{2i} \in L_{2i} 
		 \big\},
		\\
		B & = \big\{ (\ell_1,\ell_1, \ell_3,\ell_3, \ldots ) 
			\in G \colon \ell_{2i-1} \in L_{2{i-1}} 
		\big\},
	\end{split}
	\end{equation}
	where we have suppressed the maps $\vp_k$ for legibility.
Then  a diffeomorphism $A\lq G/B \lt T_p$ is given by    
\[
A(g_1,g_2,g_3,g_4,\cdots)B
	\lmt
	g_1 \otimes g_2\- \otimes g_3 \otimes g_4\- \otimes \cdots \otimes \smash{g_n^{(-1)^n}L_n}
\mathrlap.\qedhere
\]
      \end{proposition}

\subsection{Diagram components, direct factors, and simplifications }\label{sec:simplify}

We now explain 
how the expression $T_p$ of \Cref{thm:explicit-description}
can be simplified:
it factors and those factors tend to telescope.
Afterwards we will give several concrete examples.

\begin{corollary}\label{thm:component-product}
Every unitary or orthogonal \GZ fiber 
is diffeomorphic to a direct product of factors indexed 
by certain connected components of the pattern. 
For a unitary GZ pattern, each component contributes one factor.
For an orthogonal GZ pattern, each component 
of the positive part and each white component contributes one direct factor.
\end{corollary}

\begin{proof}
  By \Cref{thm:explicit-description}, $\Psi\-(p)$
  is diffeomorphic to the iterated balanced product $T_p$.
The groups $L_k$ and $H_k$ decompose into direct products of factors which
correspond to connected components of the $k$th row of the GZ pattern.
\eqref{eq:J-U} and
\eqref{eq:J-O}
show the 
inclusions $L_k \to H_k$ 
send each block of $L_k$ into the block of $H_k$
corresponding to the same component of the GZ pattern,
and 
\eqref{eq:K-U} and 
\eqref{eq:K-O}
show the same for the injections $c_{a_{k+1}} \o \vp_k\: L_k \lt H_{k+1}$.
Now observe that given an indexed list 
$M_j \from K_j \to N_j$ of diagrams of groups,
one has $\dsp{\textstyle{\prod M_j}}
\ox_{\prod K_j} {\textstyle\prod N_j} \iso {\textstyle\prod_{}} \, (M_j \ox_{K_j} N_j)$.
\end{proof}
\begin{notation}\label{def:FSO}
	In the product decomposition of an orthogonal GZ fiber
	given by \Cref{thm:component-product},
	we write 
	$\defm{F_\U}$ for the product of the factors 
	corresponding to white components
	in the corresponding \GZ pattern
	and $\defm{F_\SO}$ for
	the product of the factors corresponding to white components.
	Hence $F \iso F_\U \x F_\SO$.
\end{notation}
\begin{discussion}[The two types of direct factor]
It is often possible to obtain much 
more concise expressions for the factors described in \Cref{thm:component-product}.
Suppose we are given 
a connected component of a GZ pattern (unitary or orthogonal) 
beginning in row $\defm{k_b}$ and ending in row $\defm{k_t}$.
For any integer $k \in [k_b, k_t]$, 
write $\defm{w}(k)$ for the \defd{width}
of row $k$ of the component, 
i.e., the number of vertices in the component in this row.
The width $w(k_b)$
of the bottom row in the connected component is always $1$,
and $w(k_t)$ is also $1$ unless 
$k_t = n+1$, meaning the component
ends in the top row of the GZ pattern.
There are two different cases to analyze,
depending whether or not $k_t < n+1$.

If $k_t \leq n$,
then 
the direct factor corresponding to the given connected component is of the form
\begin{equation}\label{eq:description-v1}
	J_{k_b} \ox_{K_{k_b}} J_{k_b+1} \ox_{K_{k_b+1}} \dots \ox_{K_{k_t-2}} J_{k_t-1} \ox_{K_{k_t-1}} J_{k_t}
\end{equation}
where the factors $K_k \leq J_k$ for $k_b \leq k \leq k_t$ 
are the factors of $L_k$ and $H_k$ (respectively) described 
by \eqref{eq:J-U} and \eqref{eq:K-U} in the unitary case
and \eqref{eq:J-O} and \eqref{eq:K-O} in the orthogonal case.
Note that $K_{k_t}$ is trivial since $k_t \leq n$. 
For a black component (one of a pair, in the orthogonal case),
we have $J_k \iso \U\big(w(k)\big)$,
and for a white component,
we have $J_k \iso \SO\big(w(k)\big)$.
On the other hand, if
$k_t = n+1$, there is no group $J_{n+1}$, and 
the factor corresponding to 
the given connected component has the form
\begin{equation}\label{eq:top-row-description-v1}
	J_{k_b} \ox_{K_{k_b}} J_{k_b+1} \ox_{K_{k_b+1}} \dots \ox_{K_{n-1}} J_{n}/K_{n}.
\end{equation}
\end{discussion}

\begin{proposition}\label{thm:only-standard}
	In the case of a unitary \GZ pattern,
	the groups $J_k$ and $K_k$
	in the above descriptions may be taken to be unitary
	groups and each inclusion $J_k \lt K_k$
	and $J_k \lt K_{k+1}$
	can be taken to be the identity or 
	a lower-right block inclusion $\U(d) \simto [1] \oplus \U(d) \inc \U(1+d)$. 
	In the case of an orthogonal \GZ pattern,
	the same holds for each factor of $F_\U$,
	and for a factor of $F_\SO$, 
	the inclusions can each be taken to be the identity or of the form
	$\SO(d) \simto [1] \oplus \SO(d) \inc \SO(1+d)$. 
\end{proposition}
\begin{proof}
For the case of a unitary GZ fiber, this follows immediately
from our description of the groups $H_k$ and $L_k$ and the
maps $c_a \o \vp\: L_k \lt H_{k+1}$.
In the orthogonal case, 
The description is similarly immediate for $F_\SO$
and components of $F_\U$ in which all $J_k$ and $K_k$ are standard unitary.

The argument for components with any nonstandard unitary $J_k$
is more involved.
We observe an iterated balanced product as in 
\eqref{eq:description-v1} or
\eqref{eq:top-row-description-v1}
is determined by a zig-zag diagram
\[ J_{k_b} \lf
K_{k_b}	\lt
J_{k_{b+1}} \lf
K_{k_{b+1}} \lt
\cdots \lf
K_{k_t -1} \lt
J_{k_t } \lf
K_{k_t}\mathrlap,
\]
where $K_{k_t} = 1$ unless $t = n+1$,
and that an isomorphism of such diagrams,
meaning a sequence of Lie group isomorphisms $J_k \isoto J'_k$
and $K_k \isoto K'_K$
connecting such a zig-zag with another 
in such a way that all squares commute,
induces a diffeomorphism of the one iterated balanced
product with the other.

Now if $k$ is even and $J_k = \U(d)^o$ is nonstandard, 
our description of $H_k$, $L_k$, 
and $c_A\: L_{k-1} \lt H_k$
tells us the following.
First, $K_k$ is again $\U(d)^o$ or is $\big([1] \oplus \U(d-1)\big)^o$;
let us agree to write $\U(e)^o$ for either.
The unitary group $J_{k+1} \iso K_k$ is standard,
with the map
$K_k \lt J_{k+1} = \U(e)$ an isomorphism given by conjugation $c_o$ 
by the transposition matrix $o$.
Following our modification of the map $K_{k-1} \lt J_k$,
there are two possibilities for $K_{k-1}$,
either $\U(d)$ or $\U(d-1)$,
and the map is given by either the identity or the lower-right block inclusion into $\U(d)$,
followed in both cases by $c_o$.
Let us agree to write $K_{k-1} = \U(f)$ in either case
and write $\defm i$ for the inclusion $\U(f) \lt \U(d)$, which may be the identity.
Then we have the following isomorphism of zig-zags:
\[
\xymatrix{
\cdots & \ar[l] \U(f) \ar[r]^{c_o \o\, i} \ar[d]^{=} & \U(d)^o \ar[d]^{c_o}& \ar@{_{(}->}[l] \ \U(e)^o \ar[d]^{c_o} \ar[r]^{c_o} & \U(e)\ar[d]^{=} & \ar[l]\cdots\\
\cdots & \ar[l] \U(f) \ar[r]_{i} & \U(d) & \ar@{_{(}->}[l] \ \U(e) \ar[r]_{=} & \U(e) & \ar[l]\cdots\mathrlap,
}
\]
where the maps not displayed are the identity, and the displayed part of the bottom zig-zag is standard.
Applying one of these isomorphisms for each nonstandard $J_k$,
we arrive at a zig-zag where all $J_k$ and $K_k$ are standard unitary,
and all maps $K_k \lt J_{k+1}$ are the identity or the inclusion of $\U(d)$ in $\U(d+1)$
as the lower-right block. 
\end{proof}

%


\begin{discussion}[Telescopy]\label{rmk:tensor-discussion}
We may further simplify the expression \eqref{eq:description-v1} as follows. 
For each integer $k \in [k_b, k_t - 1]$,
the expression contains a factor $J_k \ox_{K_k} J_{k+1}$
determined by the pair of consecutive rows $k$, $k+1$. 
When these rows form an \mtrap-shape or a \parallelogram-shape, 
then $K_k = J_{k+1}$ so 
we have the simplification $J_k \ox_{K_k} J_{k+1} \cong J_k$,
whereas if they form a \wtrap-shape, then $K_k = J_{k}$, 
so we have $J_k \ox_{K_k} J_{k+1} \cong J_{k+1}$. 
Iterating the procedure, 
one can always simplify the expression 
\eqref{eq:description-v1} into one of the following forms.

First assume $k_t \leq n$,
meaning the corresponding
GZ pattern component ends before the top row.
If each row is
of width $1$, 
then the expression collapses to one factor,
$\U(1)$ for a black component or 
$\SO(1) = {*}$ for a white component.
Otherwise, 
take a subsequence 
\[
\defm{\wt \ell_1} \leq \defm{\ell_1} \leq 
\defm{\wt \ell_2} \leq \defm{\ell_2} \leq 
\cdots \leq
\defm{\wt\ell_{r-1}} \leq \defm{\ell_{r-1}} \leq 
\defm{\wt \ell_r}\mathrlap,\]
of $[k_b+1,k_t-1]$ such that $\defm{M_i} = w(\wt \ell_i)$ 
and $\defm{m_i} = w(\ell_i)$ are, respectively, the sequences of local maximum 
and minimum values of the function $w$ on the interval $[k_b,k_t]$. 
Then the expression~\eqref{eq:description-v1} telescopes to
$J_{\wt\ell_1} \,\ox_{K_{\ell_1}}\,
J_{\wt\ell_2} \,\ox_{K_{\ell_2}}\, \cdots \,\ox_{K_{\smash{\ell_{r-1}}}}\, 
J_{\wt\ell_r}$,
which is
\begin{equation}\label{eq:description-v2}
\U(M_1) \ox_{\U(m_1)}
\cdots \ox_{\U(m_{r-1})} 
\U(M_r)
\qquad\qquad\quad\mbox{or}\quad\qquad\qquad
\SO(M_1) \ox_{\SO(m_1)}
\cdots \ox_{\SO(m_{r-1})} 
\SO(M_r)
\end{equation}
depending on whether the component is black or white.
In the case of a white component, this expression corresponds
to a decomposition of the component into maximal hexagons and pentagons.

The analysis in the 
$k_t = n+1$ expression is similar.
In the degenerate case where each row is of width $1$,
the expression collapses to $\U(1)/\U(1)$ or $\SO(1)/\SO(1)$,
which is a point either way.
Otherwise, 
much as before,
the expression~\eqref{eq:top-row-description-v1} telescopes
to $J_{\wt\ell_1} \,\ox_{K_{\ell_1}}\,
J_{\wt\ell_2} \,\ox_{K_{\ell_2}}\, \cdots \,\ox_{K_{\smash{\ell_{r-1}}}}\, 
J_{\wt\ell_r}/K_n$. 
More explicitly, depending on whether it is a black or white component respectively, this equals
\begin{equation}\label{eq:top-row-description-v2}
\U(M_1) \ox_{\U(m_1)}
\cdots \ox_{\U(m_{r-1})} 
\U(M_r)/K_n
\ \quad\quad\quad\mbox{or}\quad\quad\quad \ 
\SO(M_1) \ox_{\SO(m_1)}
\cdots \ox_{\SO(m_{r-1})} 
\SO(M_r)/K_n
\end{equation}
where $M_i = w(\tilde \ell_i)$ and $m_i = w(\ell_i)$ as before.
The group $K_n$ is determined by the shape that occurs
in the connected component between the top two rows (see examples below).
\end{discussion}

\begin{discussion}[Stiefel bundles]\label{rmk:Stiefel}
  It will be useful later that \eqref{eq:description-v1}
  and \eqref{eq:top-row-description-v1} 
  give an iterated
fiber bundle construction for each of the direct factors
described in \Cref{thm:component-product},
with base $J_{\smash{\wt\ell_1}}/K_{\ell_1}$ and fibers
$J_{\smash{\wt\ell_i}}/K_{\ell_i}$ and either 
$J_{\smash{\wt\ell_r}}/K_n$
or
$J_{\smash{\wt\ell_r}}$.
Similarly, \eqref{eq:description-v2} and \eqref{eq:top-row-description-v2}
give another
iterated fiber bundle construction
for the direct factors,
with base $\U(M_1)/\U(m_1)$ or $\SO(M_1)/\SU(m_1)$
and fibers $\U(M_j)/\U(m_j)$ or $\SO(M_j)/\SO(m_j)$
or $\U(M_r)$ or $\SO(M_r)$.
These are all complex Stiefel manifolds $V_{M-m}(\C^M) = \U(M)/\U(m)$
or real Stiefel manifolds $V_{M-m}(\R^M) = \SO(M)/\SO(m)$ 
(with $m = 0$ recovering $\U(M)$ and $m \in \{0,1\}$, $\SO(M)$).
Stiefel manifolds are themselves interpretable as iterated sphere bundles over spheres 
using $\U(k+1)/\U(k) \iso S^{2k+1}$ and $\SO(k+1)/\SO(k) \iso S^k$,
so this iterated bundle description gives a refinement
of the existing iterated sphere bundle descriptions of Cho--Kim--Oh and Cho--Kim%
~\cite{ChoKimOh2020,ChoKim2020}.
\end{discussion}

\brmk\label{rmk:xex}
One might hope therefore that the cohomology ring of a GZ fiber 
should agree with that of the product of the iterated 
Stiefel fibers.
In \Cref{thm:main-U} we will show this is true of unitary GZ fibers.
For orthogonal GZ fibers, 
we will see in \Cref{thm:main-O} that this 
is true over $\Z[\sfrac 1 2]$
and for the cohomology \emph{groups} over $\Z$,
but we find in \Cref{rmk:xex2}
that it is not true for the cohomology ring over $\Z$.
\ermk

\begin{example}[Diamonds in unitary GZ patterns]\label{eg:diamonds} 
	Suppose that a connected component of a unitary GZ pattern 
	is such that the sequence of local minimum and maximum widths, omitting
	first and last rows, contains a single entry
	(necessarily a local maximum, with $w \geq 2$).
	If the component ends before the top row, 
	the corresponding direct factor is $\U(w)$. 
	If it ends in the top row, 
	the corresponding factor is $\U(w)/\U(1) \cong \SU(w)$. 
	In the special case where the sequence of widths
	strictly increases to $w$ and strictly decreases afterward,
	this recovers the existing description of
	``diamond singularities''~\cite[\SS5(e)]{BoulocMirandaZung2018}. 
\end{example}

\begin{example}[Diamonds in orthogonal GZ patterns]\label{ex:diamond} 
	Suppose a white connected component of an orthogonal GZ pattern 
	is a diamond of width $w \geq 2$,
	i.e.,~the sequence of local minimum and maximum widths
	increases to $w$ and then decreases to $1$,
	as in \Cref{fig:V2R3} or \Cref{fig:SO4}.
	Then the corresponding direct product factor is $\SO(w)$. 
\end{example}

\begin{example}[Torus fibers in unitary GZ systems]\label{eg:tori} 
As noted above, a connected component of a unitary GZ pattern 
that has maximum width equal to $1$ and 
ends below the top row
contributes a factor of $\U(1)$, i.e.,~a circle. 
A connected component 
ending in the top row and with $w(k)$ is nondecreasing
contributes $\U\big(w(n)\big)/\U\big(w(n)\big)$, i.e.,~a point. 
In particular, a unitary GZ fiber is a torus
if and only if 
all the connected components of
the associated GZ pattern have one of these forms,
and in this case, the dimension of the torus is 
the number of connected components not connected to the top row. 
This recovers the description of 
``elliptic non-degenerate singular fibers'' given 
by Bouloc \emph{et al.}~\cite[\SS5(c)]{BoulocMirandaZung2018}. 
These are not the only circle factors in general,
as we discuss in \Cref{sec:circles}.
\end{example}

\begin{example}[Torus fibers in orthogonal GZ systems]\label{ex:torus}
We can similarly characterize toral factors of an orthogonal GZ fiber
contributed by factors corresponding to white components of the pattern,
as follows.
A white connected component consisting of an isolated vertex contributes a point. 
A white diamond of width $2$ 
(i.e.,~a white connected component consisting of three rows with widths $1,2,1$) contributes $\SO(2)/\SO(1)$,~i.e.,~a circle. 
A white connected component 
ending in the top row and with $w(k)$ nondecreasing
contributes a point. 
An orthogonal GZ fiber is a torus if and only if all 
the white components have one of these forms and all 
the black components have one of the forms from Example~\ref{eg:tori}.
\end{example}

\begin{example} 
	Consider the orthogonal GZ pattern in Figure~\ref{fig:orthogonal-GZ-pattern}. The black diamond in the positive part
	contributes a factor of $\U(2)/\U(1) \cong \SU(2)$
	and the white diamond of width $2$ contributes a factor of $\SO(2)$. 
	Thus the GZ fiber is diffeomorphic to $\SU(2) \times \SO(2)$.
\end{example}

It will be important to our later analysis
that all connected real and complex Stiefel manifolds (i.e., all but $\O(n)$)
can be realized as GZ fibers.

\begin{figure}
\begin{subfigure}[b]{0.25\textwidth}
	\begin{center}
  	\begin{tikzpicture}[scale =.5,line cap=round,line join=round,>=triangle 45,x=1cm,y=.866cm, every node/.style={scale=1.25}]
	\begin{scriptsize}
	\foreach \j in {0,...,3}
	{
		\foreach \i in {0,...,\j}
		{
			\draw [fill=white] (\i-.5*\j,\j) circle (2.5pt);
		}
	}
	\foreach \j in {0,...,3}
	{
		\draw [line width=1pt] (.5*\j,\j)-- (.5*\j-1.5,3+\j);
		\draw [line width=1pt] (-.5*\j,\j)-- (1.5-.5*\j,3+\j);
	}
	\draw [line width=1pt] (-2,4)-- (-3.5,7);	
	\draw [line width=1pt] (-1.5,5)-- (-2.5,7);	
	\draw [line width=1pt] (-1,6)-- (-1.5,7);	
	\draw [line width=1pt] (2,4)-- (3.5,7);	
	\draw [line width=1pt] (1.5,5)-- (2.5,7);	
	\draw [line width=1pt] (1,6)-- (1.5,7);	
	\draw [line width=1pt] (-2,4)-- (-0.5,7);	
	\draw [line width=1pt] (-2.5,5)-- (-1.5,7);	
	\draw [line width=1pt] (-3,6)-- (-2.5,7);	
	\draw [line width=1pt] (2,4)-- (0.5,7);	
	\draw [line width=1pt] (2.5,5)-- (1.5,7);	
	\draw [line width=1pt] (3,6)-- (2.5,7);	
	
	\draw [line width=1pt] (-.5,1)-- (.5,1);	
	\draw [line width=1pt] (-1.5,3)-- (1.5,3);
	\draw [line width=1pt] (-.5,5)-- (.5,5);
	\draw [line width=1pt] (-1.5,5)-- (-2.5,5);	
	\draw [line width=1pt] (1.5,5)-- (2.5,5);
	\draw [line width=1pt] (-.5,7)-- (-3.5,7);	
	\draw [line width=1pt] (.5,7)-- (3.5,7);	
	
	\draw [line width=1pt] (-1,2)-- (1,2);
	\draw [line width=1pt] (-1,4)-- (1,4);
	\draw [line width=1pt] (-3,6)-- (-1,6);
	\draw [line width=1pt] (3,6)-- (1,6);
	
	\foreach \j in {0,...,3}
	{
		\foreach \i in {0,...,\j}
		{
			\draw [fill=white] (\i-.5*\j,\j) circle (2.5pt);
		}
	}
	
	\foreach \j in {0,...,3}
	{
		\foreach \i in {0,...,\j}
		{
			\draw [fill=black] (2+\i-.5*\j,4+\j) circle (2.5pt);
		}
	}
	
	\foreach \j in {0,...,3}
	{
		\foreach \i in {0,...,\j}
		{
			\draw [fill=black] (-2+\i-.5*\j,4+\j) circle (2.5pt);
		}
	}
	
	\foreach \j in {0,...,2}
	{
		\foreach \i in {0,...,\j}
		{
			\draw [fill=white] (\i-.5*\j,6-\j) circle (2.5pt);
		}
	}
	
	\end{scriptsize}
	\end{tikzpicture}
\end{center}
\caption{$\SO(4)$}
\label{fig:SO4}
\end{subfigure}
\begin{subfigure}[b]{0.2\textwidth}
  \begin{center}
	\begin{tikzpicture}[scale =.5,line cap=round,line join=round,>=triangle 45,x=1cm,y=.866cm, every node/.style={scale=1.25}]
	\begin{scriptsize}
	
	\foreach \j in {0,...,2}
	{
		\draw [line width=1pt] (.5*\j,\j)-- (.5*\j-1,2+\j);
		\draw [line width=1pt] (-.5*\j,\j)-- (1-.5*\j,2+\j);
	}

	\draw [line width=1pt] (-1.5,3)-- (-2,4);
	\draw [line width=1pt] (-1.5,3)-- (-1,4);
	\draw [line width=1pt] (1.5,3)-- (2,4);
	\draw [line width=1pt] (1.5,3)-- (1,4);
	
	\draw [line width=1pt] (-.5,1)-- (.5,1);	
	\draw [line width=1pt] (-.5,3)-- (.5,3);
	
	\draw [line width=1pt] (-1,2)-- (1,2);
	\draw [line width=1pt] (-1,4)-- (-2,4);
	\draw [line width=1pt] (1,4)-- (2,4);
	
	\foreach \j in {0,...,2}
	{
		\foreach \i in {0,...,\j}
		{
			\draw [fill=white] (\i-.5*\j,\j) circle (2.5pt);
		}
	}

	\foreach \j in {0,1}
	{
		\foreach \i in {0,...,\j}
		{
			\draw [fill=black] (-1.5+\i-.5*\j,3+\j) circle (2.5pt);
		}
	}
	
	\foreach \j in {0,1}
	{
		\foreach \i in {0,...,\j}
		{
			\draw [fill=black] (1.5-\i+.5*\j,3+\j) circle (2.5pt);
		}
	}

	\foreach \j in {0,1}
	{
		\foreach \i in {0,...,\j}
		{
			\draw [fill=white] (\i-.5*\j,4-\j) circle (2.5pt);
		}
	}
	
	\end{scriptsize}
	\end{tikzpicture}
	\end{center}
	\caption{$\SO(3) = V_2(\R^3)$}
	\label{fig:V2R3}
\end{subfigure}
\begin{subfigure}[b]{0.25\textwidth}
\begin{center}
	\begin{tikzpicture}[scale =.5,line cap=round,line join=round,>=triangle 45,x=1cm,y=.866cm, every node/.style={scale=1.25}]
	\begin{scriptsize}
	\foreach \j in {0,...,3}
	{
		\foreach \i in {0,...,\j}
		{
			\draw [fill=white] (\i-.5*\j,\j) circle (2.5pt);
		}
	}
	\foreach \j in {0,...,2}
	{
		\draw [line width=1pt] (.5*\j,\j)-- (.5*\j-1.5,3+\j);
		\draw [line width=1pt] (-.5*\j,\j)-- (1.5-.5*\j,3+\j);
	}
	\draw [line width=1pt] (-1.5,3)-- (-.5,5);
	\draw [line width=1pt] (1.5,3)-- (.5,5);
	
	\draw [line width=1pt] (-2,4)-- (-2.5,5);	
	\draw [line width=1pt] (2,4)-- (2.5,5);	
	\draw [line width=1pt] (-2,4)-- (-1.5,5);	
	\draw [line width=1pt] (2,4)-- (1.5,5);
	
	\draw [line width=1pt] (-.5,1)-- (.5,1);	
	\draw [line width=1pt] (-1.5,3)-- (1.5,3);
	\draw [line width=1pt] (-.5,5)-- (.5,5);
	\draw [line width=1pt] (-1.5,5)-- (-2.5,5);	
	\draw [line width=1pt] (1.5,5)-- (2.5,5);
	
	\draw [line width=1pt] (-1,2)-- (1,2);
	\draw [line width=1pt] (-1,4)-- (1,4);
	
	\foreach \j in {0,...,3}
	{
		\foreach \i in {0,...,\j}
		{
			\draw [fill=white] (\i-.5*\j,\j) circle (2.5pt);
		}
	}
	
	\foreach \j in {0,1}
	{
		\foreach \i in {0,...,\j}
		{
			\draw [fill=black] (2+\i-.5*\j,4+\j) circle (2.5pt);
		}
	}
	
	\foreach \j in {0,1}
	{
		\foreach \i in {0,...,\j}
		{
			\draw [fill=black] (-2+\i-.5*\j,4+\j) circle (2.5pt);
		}
	}
	
	\foreach \j in {1,2}
	{
		\foreach \i in {0,...,\j}
		{
			\draw [fill=white] (\i-.5*\j,6-\j) circle (2.5pt);
		}
	}
	
	\end{scriptsize}
	\end{tikzpicture}
      \end{center}
      \caption{$\dsp\smash{\frac{\SO(4)}{\SO(2)}} = V_2(\R^4)$}
	\label{fig:V2R4}
\end{subfigure}
\begin{subfigure}[b]{0.2\textwidth}
  \begin{center}
	\begin{tikzpicture}[scale =.5,line cap=round,line join=round,>=triangle 45,x=1cm,y=.866cm, every node/.style={scale=1.25}]
	\begin{scriptsize}
	\foreach \j in {0,...,3}
	{
		\foreach \i in {0,...,\j}
		{
			\draw [fill=white] (\i-.5*\j,\j) circle (2.5pt);
		}
	}
	\foreach \j in {0,...,2}
	{
		\draw [line width=1pt] (.5*\j,\j)-- (.5*\j-2,4+\j);
		\draw [line width=1pt] (-.5*\j,\j)-- (2-.5*\j,4+\j);
	}
	\draw [line width=1pt] (-2,4)-- (-1,6);
	\draw [line width=1pt] (2,4)-- (1,6);
	
	\draw [line width=1pt] (-1.5,3)-- (0,6);
	\draw [line width=1pt] (1.5,3)-- (0,6);
	
	\draw [line width=1pt] (-2.5,5)-- (-3,6);	
	\draw [line width=1pt] (2.5,5)-- (2,6);	
	\draw [line width=1pt] (-2.5,5)-- (-2,6);	
	\draw [line width=1pt] (2.5,5)-- (3,6);
	
	\draw [line width=1pt] (-.5,1)-- (.5,1);	
	\draw [line width=1pt] (-1.5,3)-- (1.5,3);
	\draw [line width=1pt] (-1.5,5)-- (1.5,5);
	\draw [line width=1pt] (-2,6)-- (-3,6);	
	\draw [line width=1pt] (2,6)-- (3,6);	
	
	\draw [line width=1pt] (-1,2)-- (1,2);
	\draw [line width=1pt] (-2,4)-- (2,4);
	\draw [line width=1pt] (-1,6)-- (1,6);
	\foreach \j in {0,...,4}
	{
		\foreach \i in {0,...,\j}
		{
			\draw [fill=white] (\i-.5*\j,\j) circle (2.5pt);
		}
	}
	
	\foreach \i in {-1,...,2}
	{
		\draw [fill=white] (\i-.5,5) circle (2.5pt);
	}
	\foreach \i in {-1,...,1}
	{
		\draw [fill=white] (\i,6) circle (2.5pt);
	}
	
	\foreach \j in {0,1}
	{
		\foreach \i in {0,...,\j}
		{
			\draw [fill=black] (2.5+\i-.5*\j,5+\j) circle (2.5pt);
		}
	}
	
	\foreach \j in {0,1}
	{
		\foreach \i in {0,...,\j}
		{
			\draw [fill=black] (-2.5-\i+.5*\j,5+\j) circle (2.5pt);
		}
	}
	
	\end{scriptsize}
	\end{tikzpicture}
\end{center}
\caption{$\dsp\smash{\frac{\SO(5)}{\SO(3)}} = V_2(\R^5)$}
\label{fig:V2R5}
\end{subfigure}
\caption{Some orthogonal \GZ patterns}
\label{fig:Stiefel}
\end{figure}

\begin{example}[Stiefel manifolds]\label{ex:Stiefel} 
	Suppose 
	a connected component of a GZ pattern 
	contains vertices in the top row and the sequence of widths 
	weakly increases from $1$ to $w = w(\ell)$ and thereafter
	weakly decreases to $v = k(n+1)$
(for instance the component could be a pentagon),
as in \Cref{fig:V2R4,fig:V2R5}.
Then the corresponding direct factor is 
the complex Stiefel manifold $V_{w-v}(\C^w) \iso \U(w) /\U(v)$ 
if the component is black,
and the real Stiefel manifold $V_{w-v}(\R^w) \iso \O(w)/\O(v) \iso \SO(w)/\SO(v)$ if it is white. 
We saw the unitary groups 
$V_w(\C^w)\iso \U(w)$ 
realized as GZ fibers in \Cref{ex:diamond},
so every complex and real Stiefel manifold except for the disconnected
$V_w(\R^w) \iso \O(w)$ 
occurs as a GZ fiber.
\end{example}

The direct product splitting in \Cref{thm:component-product} 
is still not the end of the story
on simplified descriptions 
of direct factors of GZ fibers.
We make two further observations.

\begin{discussion}[Pinching]\label{rmk:pinch}
First, since $\SO(1) = 1$,
a ``pinch'' row of width $1$ in a white component
corresponds to a direct product splitting
$X \ox_{\SO(1)} Y \iso X \x Y$,
thus leading to a refinement of the direct product decomposition
in \Cref{thm:explicit-description} for an orthogonal \GZ fiber.
\end{discussion}

\begin{discussion}[Opportunistic splittings]\label{rmk:opportunism}
Second, one can often split off Stiefel manifolds as direct factors.
Indeed, if $M_1 \leq M_2$
in the notation of \eqref{eq:description-v2},
then the map
\eqn{
	\U(M_1) \ox_{\U(m_1)} \U(M_2) \ox_{\U(m_2)} Z
	&\lt
	\Big(\quotientmed{\U(M_1)}{\U(m_1)}\Big) \,\mnn\x\, \U(M_2) \ox_{\U(m_2)} Z,\\
	g \ox h \ox z 
	&\lmt
	\big(g\U(m_1), gh \ox z\big)
}
is a diffeomorphism with inverse 
$\big(g\U(m_1), h' \ox z\big) \lmt g \ox g\-h' \ox z$,
and similarly if the groups are special orthogonal.
One 
can repeat this until one encounters a local maximum width $M_p > M_{p+1}$.
For factors corresponding to 
pattern components terminating before the top row,
if $M_r < M_{r-1}$, 
one can also peel Stiefel factors from
the right-hand side.
The simplest examples for which this is not possible are
\[\U(3) \ox_{\U(1)} \U(2)/\U(1)
\qquad\mbox{and}\qquad
\U(3) \ox_{\U(1)} \U(2) \ox_{\U(1)} \U(3)
\mathrlap.
\]
\end{discussion}

The direct factors of \Cref{thm:component-product}
can also be realized as biquotients, 
in a way clear from
the proof of \Cref{thm:biquotient},
\Cref{thm:component-product},
\Cref{thm:only-standard},
and \Cref{rmk:tensor-discussion}:

\begin{proposition}\label{thm:biquotient-factor}
  Every direct product factor of a GZ fiber 
  can be written as a biquotient $A\backslash G/B$ as in \Cref{thm:biquotient}.
\end{proposition}

The following example applies our description to a more 
complicated unitary GZ fiber.
To our knowledge, 
such examples have not been previously been so simply
described in the literature.

\begin{example} \label{eg:intro example diffeo type}
	Consider the unitary GZ pattern in Figure~\ref{fig:unitary-GZ-pattern}. 
	The double diamond on the left contributes $\U(2)\ox_{\U(1)}\U(2)$,
	or 
	$S^3 \x \U(2) \iso (S^3)^2 \times S^1
	$ by \Cref{rmk:opportunism}
	and the decomposition $\U(2) = \big\{\mn\diag(z,z\-): z \in \U(1)\big\} \.\SU(2)$.
The long component on the right contributes 
	$\U(2)/\U(1) \cong S^3$ as in \Cref{eg:diamonds}. 
	Finally, the largest component contributes 
	$\U(4) \ox_{\U(2)} \U(3)/\U(2)$.
	The factors associated to the remaining connected components 
	are $\U(1)$, as explained  in \Cref{eg:tori}. 
	All told, the associated GZ fiber is diffeomorphic to 
	\[
	(S^1)^7 \,\x\, (S^3)^3 \,\x\,
	\U(4) \ox_{\U(2)} \U(3)/\U(2)
	\mathrlap,
	\]
	where the last factor can equivalently be written
	as a biquotient
	$ \U(2)\backslash\big(\U(4) \times \U(3) \big)/\U(2)$. 
\end{example}

\subsection{Circle factors}\label{sec:circles}

Some known observations concerning circle actions on
and circle factors of 
\GZ fibers 
fall more or less immediately, 
and in more transparent form, out of our analysis.

\begin{discussion}[The determinant and circle extraction]\label{rmk:circle-action}
The block-inclusions 
$\U(d) \lt \U(d+1)$
(standard by \Cref{thm:only-standard})
figuring in the iterated balanced products
\eqref{eq:description-v2} and
\eqref{eq:top-row-description-v2}
evidently preserve determinants.
Any unitary factor $Y$ as in \eqref{eq:description-v2}
corresponding to a component of the \GZ pattern
not ending in the top row thus admits a well-defined determinant map to $\U(1)$
given by $\defm{\det}\:g_1 \ox \cdots \ox g_r \lmt \det(g_1)\cdots\det(g_r)$.
There is a natural right action of $\U(1)$ on $Y$
induced by the lower-right-corner block inclusion in $\U(M_r)$,
and thus
\quation{
\begin{aligned}
\label{eq:unitary-circle-factor-splitting}
{\textstyle\det\-(1)} \x \U(1) 
& \lt
\U(M_1) \ox_{\U(m_1)}
\cdots \ox_{\U(m_{r-1})} 
\U(M_r),\\
(g_1 \ox \cdots \ox g_r,z) 
&\lmt 
g_1 \ox \cdots \ox g_r\cdot z 
\end{aligned}
}
is a well-defined diffeomorphism,
with inverse $x \lmt \big(x\.(\det x)\-,\det x\big)$.
Thus every component of the \GZ pattern not ending in the top row
contributes a circle factor to the fiber.

As for the other factor, $\det\-(1)$,
note that the natural map 
$\prod \SU(M_j) \inc \prod \U(M_j) \epi Y$
surjects onto it,
as can be seen by iteratively replacing factors, 
using the left and right $\U(1)$-actions on each $\U(m_j)$
induced by the lower-right-corner inclusion:
$g_1 \ox g_2 \ox \cdots
= g_1 \.(\det g_1)\- \otimes (\det g_1)g_2 \ox \cdots$
and so on.
One can show a similar surjection if $Y$
is of the form \eqref{eq:top-row-description-v2},
corresponding to a factor ending in the top row,
in which case the $\U(1)$-action on the right is quotiented out.
Thus 
 $ {\textstyle\det\-(1)} $ equals
\quation{\label{eq:SU-factor}
  \SU(M_1) \ox_{\SU(m_1)}
  \cdots \ox_{\SU(m_{r-1})} 
  \SU(M_r)
\ \quad\quad\mbox{or}\quad\quad\ 
Y = 
  \SU(M_1) \ox_{\SU(m_1)}
  \cdots \ox_{\SU(m_{r-1})} 
  \SU(M_r)/\mr{S}K_n
  \mathrlap.
}
\end{discussion}

These observations apply equally 
to the unitary factor $F_\U$ 
of an orthogonal GZ fiber
As noted in \Cref{ex:torus}, 
width-$2$ white diamonds in a GZ pattern
also contribute $\SO(2)$ components to the corresponding fiber.

\begin{notation}
For a unitary \GZ system,
let $\defm t = \defm{ t(p)}$ be the number of components 
not ending in the top row
of the \GZ pattern 
corresponding to $p$.
(In other words,
this is the number of \mtrap-shapes of width $1$
in the first $n$ rows.)
For an orthogonal \GZ system,
let $t(p)$ be the number of such 
components in the positive part of the GZ pattern,
plus the number of diamonds of width $2$
in the white components.
\end{notation}

Considering all factors of $\Psi\-(p)$ together,
the product of the determinant maps
is a map onto the torus $T^{t(p)}$
which we will call $\defm{\Det_p}$.
Thus we have the following:

\begin{proposition}\label{thm:circle-factor}
 There is for each $p \in \D$ a $T^{t(p)}$-equivariant diffeomorphism
 $T_p \iso \Det_p\-(1) \x T^{t(p)}$.
\end{proposition}

\nd In the unitary case,
the splitting off of a $T^{t}$ factor
is result of Cho--Kim--Oh~\cite[Thm.~6.11]{ChoKimOh2020}.

\begin{corollary}
 The fiber $\Psi\-(p)$ over a point $p$ in the interior
 of a unitary or orthogonal GZ polytope $\Delta_{\l^{(n+1)}}$
 is a torus.
\end{corollary}

\nd In the unitary case,
this is result of Cho--Kim--Oh~\cite[Thm.~6.14]{ChoKimOh2020}
and also of Bouloc--Miranda--Zung~\cite[\SS5(b)]{BoulocMirandaZung2018}.

\begin{discussion}[The topological Thimm trick and coordinates]
We compare  the torus action of \Cref{rmk:circle-action} with
the traditional torus action obtained by the Thimm trick
  (see Cho--Kim--Oh~\cite{ChoKimOh2020}
  and Pabiniak~\cite[\SS3.5]{Pabiniak2012}).
  In the notation from \Cref{rmk:stabilizer-bundle},
  each centralizer $\smash{Z_{H_j}(L_j)}$
  acts on $\Psi\-(p)$ by $z_j\.\genpt \ceq g_j z_j g_j\- \. \genpt$,
  where the action on the left is the coadjoint action of $G_{n+1}$
  and the map $\vp_n\o\cdots\o\vp_j$ has been suppressed.
Putting these all together, we have an action of $\prod Z_{H_k}(L_k)$ on $\Psi\-(p)$.%
\footnote{\ 
    It is clear from the definition that this action is continuous and
    does not depend on the chosen representative $g_j \in G_j/L_j$.
    To see the actions,
    for varying $j$, commute, note 
    that $z_\ell \. g_j z_j g_j^{-1} \genpt
    = 
    \check g_\ell z_\ell \check g_\ell^{-1} (g_j z_j g_j^{-1} \genpt)$
    for $\check g_\ell^{-1} $
    such that 
    $\check g_\ell^{-1} (g_jz_jg_j^{-1} \genup{\ell} ) = \lup{\ell} $;
    if $\ell > j$, then
    $\check g_\ell^{-1} = g_\ell^{-1} (g_jz_jg_j^{-1})^{-1}$
    works.
    On the other hand,
    $z_j \. g_\ell z_\ell g_\ell^{-1} \genpt$ is just
    $g_j z_j g_j^{-1} g_\ell z_\ell g_\ell^{-1}\genpt$ 
    since $\Phi^{n+1-j}(g_\ell z_\ell g_\ell^{-1}\genpt) = \genup j$.
}

Using \Cref{rmk:a} to convert from $\genpt \in \Psi\-(p)$
to $(g_k L_k) \in \W^r_p$,
the entries of $z_j\.(g_k L_k)$
are given by
  \[
	  \case{ 
	  		g_k L_k, & k < j,\\
	  		g_j z_j g_j\- g_k L_k, & k \geq j,
		}
  \]
where again $\vp_{k-1}\cdots\vp_j$ is suppressed. 
The coordinate $k = j$ gives $\big(z_j\.(g_k L_k)\big)_j = g_j z_j L_j$,
but it is useful to absorb $k = j$  into the $k \geq j$ case
to avoid case distinctions.
Converting this action on $\W^r_p$ into an action on $T_p$
using \Cref{rmk:explicit-conversion},
one finds that only the~$\smash{j\th}$ tensor-coordinate is changed
under the action of $z_j \in Z_{H_j}(L_j)$, 
so that 
the point
  \eqref{eq:tensor-conversion} maps to
  \quation{\label{eq:circle-tensor}
	a_\a 	g_\a			\, \ox\,
	a_{\a+1}  g_\a\- g_{\a+1} \,\ox \,
	\cdots \,\ox \,
	a_{j-1} g_{j-2}\- g_{j-1} \,\ox \,
	a_{j}  g_{j-1}\- g_{j}\textcolor{red}{ z_j } \,\ox \,
	a_{j+1} g_j\- g_{j+1} \,\ox \,
	\cdots \,\ox\,
	a_n g_{n-1}\-g_n L_n \mathrlap.
	}
  Factoring $H_k$ and $L_k$ respectively as $\Direct \U(d_i)$
  and $\Direct K_i$ as in \eqref{eq:J-U} and \eqref{eq:K-U},
  one evidently has $Z_{H_j}(L_j) = \Direct Z_{\U(d_i)}(K_i)$,
 with
 $Z_{\U(d_i)}(K_i) = Z\big(\U(d_i)\big)$ the diagonal copy of $\U(1)$
  if $\U(d_i) = K_i$,
  and otherwise $Z_{\U(d_i)}(K_i) \iso \U(1) \x Z(K_i)$
  if $K_i = [1] \+ \U(d_i - 1)$ (or $\U(d_i -1) \+ [1]$),
  where the $\U(1)$ factor corresponds to and acts on the $r_i$
  coordinate of $\xip{j+1}$ in \eqref{eq:Y}.
  If we restrict our attention to factors $Z_{\U(1)}\big(\U(0)\big)$,
  then on splitting $T_p$ into direct factors of the form in \eqref{eq:description-v2},
  we see from \eqref{eq:circle-tensor} that 
  the relevant $\U(1)$ act only on the right of the final tensor-component 
  of each direct factor,
  so this is the same action of $T^{t(p)}$
  we introduced in \Cref{rmk:circle-action}:\footnote{\ 
  An interesting note is that the action of $\smash{Z_{H_j}(L_j)}$
	is typically noncanonical, as $L_j$ depends on $\xip{j+1}$,
	which there are various reasonable choices of,
	but the action of $T^t$ is uniquely defined.}
\end{discussion}


\bthm\label{thm:Thimm}
The torus action of \Cref{rmk:circle-action}
is that given by the Thimm trick.
\ethm


\brmk\label{rmk:bmz-circles}
Bouloc--Miranda--Zung~\cite[(5.9)]{BoulocMirandaZung2018}
give a description of
multi-diamond singularities of unitary GZ systems as a product of unitary
groups (\emph{cf}.~\Cref{eg:diamonds}),
in particular observing that
every direct factor corresponding to a diamond can be written
as $\U(n) \cong \SU(n) \times \U(1)$.
This is of course a special case of \Cref{thm:circle-factor}.
\ermk

\brmk\label{rmk:cko-torus}
Cho--Kim--Oh obtained \Cref{thm:circle-factor}
for a unitary GZ fiber,
observing the other factor $Y$
is a simply-connected manifold~\cite[Theorem 6.11]{ChoKimOh2020}
(we will recover this more simply in \Cref{thm:pi-1-unitary}).
They obtain this decomposition in two steps,
first observing the~$t$ GZ functions
corresponding to the top vertices of each connected component 
in the GZ pattern with $k_t \leq n$ 
generate an action of $T^t$ on the GZ fiber,
and second using toric degenerations 
to show this $T^t$-action is free.
Our recovery of this decomposition is more elementary,
since we do not require a toric degeneration,
and our decomposition is also finer,
since 
\Cref{thm:component-product} and 
\eqref{eq:unitary-circle-factor-splitting} 
show the simply-connected space $Y$
is the direct product of the spaces \eqref{eq:SU-factor}.
On the other hand,
their description has much more geometric content.
\ermk

\section{The local structure of a coadjoint orbit}\label{sec:local}

We have just noted
that 
Cho--Kim--Oh obtain the torus factor of \Cref{thm:circle-factor},
in the unitary case,
as a consequence of a known toric degeneration,
whose time-$1$ gradient-Hamiltonian flow induces a quotient map
from $\Orb$ to a certain toric variety $X_0$.
Our approach in \Cref{sec:degeneration} 
will allow us to go in the other direction, 
defining a map onto a toric manifold 
in both the unitary and the orthogonal cases without reference to geometry.
Although the idea of the this map is strikingly simple, 
in fact a triviality fiberwise,
to show that the induced function 
on $\Orb$ is continuous
we will need to finally examine more than one GZ fiber at time, 
leading to local expressions over a face of the polytope 
in \Cref{sec:face},
and over the union of two faces in \Cref{sec:ray}.
These expressions in hand,
we conclude the proof in \Cref{sec:consistency}.

\subsection{Toric degeneration}\label{sec:degeneration}
We now attempt to relate the maps $\Det_p\: \Psi\-(p) \lt T^{t(p)}$ 
of \Cref{sec:circles} for varying $p$.

\begin{discussion}
  We will refer to the interior of any face of a \GZ polytope $\D$ as an \defd{open face}\footnote{\ 
	even though among these only the interior $\D^\o$
	is actually an open subset of $\D$
  }.
As the \GZ pattern of any two points in the same open face $\defm E$
of the \GZ polytope $\D = \D_{\l^{(n+1)}}$ are the same,
the function $t\: \D \lt \Z$
is constant on each open face.
Since each open face is a convex subset of $\R^{\dim \D}$, hence contractible,
and the fibers over each point of $E$ are diffeomorphic,
the restriction of $\Psi$ to $E$ 
is a trivial bundle,
but we will provide an explicit diffeomorphism with the product bundle
in \Cref{sec:face}.


Suppose now that $q \in \D$ lies in 
some codimension-one open face $\defm F \subn \ol E$.
We will show in \Cref{sec:ray} 
that the inclusion $F \longinc \ol E$ 
induces a map $f\:\Psi\-(p) \lt \Psi\-(q)$
and a homomorphism $T^f\: T^{t(p)} \lt T^{t(q)}$
rendering commutative the square
\quation{\label{eq:det-square}
\begin{gathered}
\xymatrix@C=1.5em{
\Psi\-(p) \ar[d]_(.45){\vphantom{T^f}f} 	\ar[r]^(.55){\Det_p}	& T^{t(p)}\ar[d]^(.45){T^{\smash{f}} \vphantom{f}}\\
\Psi\-(q) 		\ar[r]_(.55){\Det_q}	& T^{t(q)}\mathrlap.
	}
\end{gathered}
}

These maps, taken altogether,
express every $T^{t(p)}$ for $p \in \D$
as a quotient of $\defm T = T^{\dim \D}$
in a well-defined and coherent way, 
so that $\defm{ X_0} \ceq \Union_{p \in \D} T^{t(p)}$
can naturally be identified with a quotient of $T \x \D$, 
and hence inherits a $T$-action with orbit space $X_0/T = \D$.
The $T$-space $X_0$ underlies a toric variety
appropriately known as a \defd{Gelfand--Zeitlin toric variety}.
The various determinant maps compile into a surjective function
$\defm\Det\: \Orb \lt  X_0$
such that the composition with the projection $X_0 \lt \D$ 
taking $T^{t(p)}$ to $\{p\}$
is the GZ system $\Psi\: \Orb \lt \D$,
and we will show the following:
\end{discussion} 

\begin{theorem}\label{thm:det}
The map
$\Det\: \mc O_{\l^{(n+1)}} \lt X_0$
is well defined and continuous.
\end{theorem}

This gives a sort of underlying
topological model for a toric degeneration,
simply the open mapping cylinder of the map
$S^1 \x \Orb \epi \Orb \xtoo{\smash{\Det}} X_0$,
which fibers over the open disc in such a way that the generic fiber is $\Orb$
and the special fiber is $X_0$.

\begin{remark}\label{rmk:degeneration}
Nishinou--Nohara--Ueda~\cite{NishinouNoharaUeda2010}
construct flat toric degenerations $(X_t)_{t \in \C}$ 
with generic fiber $\mc O_{\l^{(n+1)}}$ and special fiber $X_0$
and prove using a limit argument
that the Hamiltonian flow map $\phi$ from a subset of $X_{1}$ to $X_0$
continuously extends to be defined on all of~$X_1$.
The construction of $\phi$ relies
on solving a system of ODEs,
and so is not easy to directly compare
with what we have done here,
but it is natural to expect that $\phi$ and $\Det$
agree up to a translation of each fiber $X_0|_p$ of $X_0$
by a continuously varying element of $T^{t(p)}$.
\end{remark}

\subsection{A local model over a face}\label{sec:face}

To prove the claims leading to \Cref{thm:det},
we are forced to compare nearby factors.

\begin{discussion}\label{rmk:face}
  Recall that 
the final pullback $F_\alpha^{n+1} = \Psi\-(p)$
of
 \Cref{fig:big-diagram}
 could be given as the iterated fiber product \eqref{eq:Omega},
 determined by the system of maps
\quation{\label{eq:point-zigzag}
  \cdots
  \to G_k/H_k \from G_k/L_k \to G_{k+1}/H_{k+1} \from \cdots \to 
  G_n/H_n \from G_n/L_n \to G_{n+1}/H_{n+1}
  \mathrlap,
}
where the leftward maps were projections and 
the rightward maps were described in \Cref{thm:iotak}.
The groups $H_k$
as described in \Cref{rmk:unitary-row} and \Cref{rmk:orthogonal-row}
were the same for all $p \in E$ since we chose diagonal representatives $\lup k$.
Similarly,
by an appropriate choice of $\xip{k+1}$
in \Cref{rmk:xi-explicit} and \Cref{rmk:xi-SO-explicit},
we could take the $L_k$ to all similarly be the same.
Since the $\xip{k}$ are continuous functions in 
$p = \smash{(\lup{n+1},\lup n, \cdots, \lup 1)}$,
the $a_k = a_k(p)$ are as well.
It follows that if we take the union of the zig-zags \eqref{eq:point-zigzag},
where the objects are the same but we recall from \Cref{rmk:explicit-modification}
that the maps $G_k/L_k \lt G_{k+1}/H_{k+1}$ depend on $p \in E$,
then we obtain a system
\quation{\label{eq:stair-system}
  \cdots \to E\x G_k/H_k \from E\x G_k/L_k \to E\x G_{k+1}/H_{k+1} \from \cdots \to 
  E\x G_{n+1}/H_{n+1}
}
of continuous, fiber-preserving maps of product bundles over $E$,
whose ultimate pullback $\defm{\W^r_E}$, restricted to any $p \in E$,
is $\W^r_p$.
Then $\W^r_E$ is diffeomorphic to $\Psi\-(E)$,
since the evaluation map $\big(p,(g_k L_k)\big) \lmt \vp_n(g_n)\xi_{n+1}(p)$ is
smooth and bijective on fibers.

In \Cref{rmk:modification}, 
we moved from $\W^r_p$ to $\W^c_p$ 
by replacing the original maps $G_k/L_k \lt G_{k+1}/H_{k+1}$
by
$g L_k \lmt \smash{a_{k+1}\vp_k(g)a_{k+1}\- H_{k+1}}$.
Since each map $(p,g_{k+1}H_{k+1}) \lmt \big(p,a_{k+1}(p)g_{k+1}H_{k+1}\big)$
is a diffeomorphism, the diffeomorphism type of the pullback is unaffected
if we insert these maps after the existing maps $E \x G_k/L_k \lt E \x G_{k+1}/H_{k+1}$,
and this shows the resulting final pullback~$\defm{\W^c_E}$ 
is diffeomorphic in a fiber-preserving
fashion to $\W^r_E$ and hence to $\Psi\-(E)$.

Then \Cref{thm:explicit-description}, applied fiberwise,
yields a diffeomorphism from $\W^c_E$ to
\[
\defm {T_E} \ceq
(E \x H_\alpha) \,\ox_{E \x L_\alpha}\, (E \x H_{\alpha+1}) \,\ox_{E \x L_{\alpha+1}}\, \cdots \,\ox_{E \x L_{n-1}}\, (E \x H_{n})/(E\x L_{n})\mathrlap,
\]
where the maps $E \x L_k \lt E \x H_k$ are inclusions
and the maps $E \x L_k \lt E \x H_{k+1}$ are given by 
$(p,g) \lmt \big(p,a_{k+1}(p)\vp_k(g)a_{k+1}(p)\-\big)$.
As discussed in \Cref{rmk:explicit-conversion},
although the elements $a_{k+1}(p) \in G_{k+1}$ continuously vary, 
the effect of the conjugation maps $L_k \lt H_{k+1}$
on unitary block components is the same for all $p \in E$.
Hence there is a natural identification $T_E \iso E \x T_p$,
and thus we also have an identification $\Psi\-(E) \iso E \x T_p$.
\end{discussion}

%
\begin{discussion}[Coordinate (in)dependence]
  It follows from \Cref{rmk:face} that the expression given in
\eqref{eq:hat-tensor-conversion} for the diffeomorphism 
$\W^c_E \isoto T_E$ does not depend on the point $p$ either, 
despite the appearance of the $a_k(p)$ in the expression.
From \eqref{eq:tensor-conversion}, however,
we see the same cannot be said of the composite diffeomorphism
$\W^r_E \isoto T_E$
where the expression does vary pointwise 
(nor, hence, can it be claimed of the total composite
$\Psi\-(E) \simto \W^r_E \simto \W^c_E \simto T_E$)
even though over each $p \in E$ it restricts to a diffeomorphism $\W^r_p \isoto T_p$.
\end{discussion}

\subsection{A local model over two faces}\label{sec:ray}

Now we would
like to understand $\Psi\-(E \union F)$ 
for $\ol F$ a codimension-one face of $\ol E$.
Intuitively speaking, 
since the fibers over $E$ and over $F$ are all the same,
we expect that all we really need to understand 
is what happens over a ray from $F$ into $E$.
This turns out to be essentially correct.

\begin{discussion}\label{rmk:ray}
Select a point $p \in F$ and suppose that $\mu$ is the eigenvalue labelling
the component merged from two components in the pattern of $E$.
We will consider the intersection with $\ol E$ of the affine $2$-plane $\defm P$
in which the $\lup k_j$ corresponding to all components 
except those merging to the $\mu$-component of $F$ are fixed.
Let $\smash{\defm{\l^+}}$ and $\defm{\l^-}$ denote the greatest and least of these
and $\defm{\mu^+}$ and~$\smash{\defm{\mu^-}}$ the labels on the merging components.
Then we are focusing attention on the region
$\l^+ > \mu^+ \geq \mu^- > \l^-$ of $P$
and interested in what happens as $|\mu^\pm|$ approach one another or $0$.
The numbers $\defm\mu = \frac 1 2(\mu^+ + \mu^-)$ and $\defm t = \frac 1 2(\mu^+ - \mu^-)$
also serve as coordinates on this region,
and we for now fix $\mu$ as well and consider what happens as $t \to 0$,
which is to say $\mu^\pm \to \mu$.
So we consider a closed interval $\defm{\theinterval} \subn P$
as parameterized by $t \in [0,\e_\mu]$.
\end{discussion}

\begin{discussion}[The effect on passing to a codimension-one face
  on $H_k$]\label{rmk:merge-H}

The block-diagonal group $H_k$
which is $H_k(q) = \Stab_{\U(k)} \l^{(k)}(q)$ for $p \neq q \in \theinterval$
is a subgroup of $\defm{\ol H_k} = \Stab_{\U(k)} \l^{(k)}(p)$.
For those rows in which the GZ patterns of $E$ and $F$ agree, 
we have $H_k = \ol H_k$.
On the other hand, in a row where they disagree,
we have a case distinction.
\bitem
\item
First assume that the corresponding block factors are standard unitary
(which is always the case for a unitary GZ fiber).
If $\mu^+$ occurs $\defm{m^+}$ times and $\mu^-$ $\defm{m^-}$ times 
in $\l^{(k)}(q)$,
then $\mu$ occurs $m^+ + m^-$ times 
in $\l^{(k)}(p)$, 
and the inclusion $H_k \longinc \ol H_k$ 
is the identity on the blocks corresponding to eigenvalues
other than $\mu^+$, $\mu$, and $\mu^-$,
and on these blocks, it is 
$ \U(m^+) \+ \U(m^-)  \longinc \U(m^+ + m^-) $.
%
\eitem

For an orthogonal GZ fiber,
consulting \eqref{eq:sequence-O-even},
\eqref{eq:interlacing-O-odd}, and
\eqref{eq:interlacing-O-even},
there are several more possibilities for the nontrivial
blocks of the inclusion $H_k \longinc \ol H_k$.
\bitem
\item
If $k$ is even 
with $\U_\ell < H_k$ nonstandard
with label $\mu^-$
and $\U(d_{\ell-1})$ standard with label $\mu^+$,
the block map is $\U(d_{\ell-1}) \+ \U(d_\ell)^o \longinc
\U(d_{\ell -1} + d_\ell)^{\id \oplus\, o}$.

\item
If $k$ is even,
$\U_{\ell} = \U(1) < H_k$ corresponds to $\mu^- < 0$,
and $\mu^+ > -\mu^-$ corresponds to $\U(d_{\ell-1})$,
then as $|\mu^+| - |\mu^-| \to 0$,
the block map is 
$\U(d_{\ell -1}) \+ \U(1) \longinc \U(d_{\ell-1} + 1)^o$.

\item
If $-\mu^- = 0$
and $\U_\ell$ corresponds to $\mu^+$,
the block map is $\U_\ell \+ \SO(d_0) \longinc \SO(2d_\ell + d_0)$.
\eitem
\end{discussion}

\begin{discussion}[The effect on passing to a codimension-one face
  on $L_k$]\label{rmk:merge-L}
The possibilities for the shapes of the merging 
($\mu^+$- and $\mu^-$-labelled) components
in rows $k$ and $k+1$
are precisely the following:
\[
\mbox{
\rpara\wtrap,\qquad
\mtrap\wtrap,\qquad
\mtrap\rpara,\qquad
\rpara\rpara,\qquad
\parallelogram\parallelogram,\qquad
\parallelogram\mtrap,\qquad
\wtrap\mtrap,\qquad
\wtrap\parallelogram
}
.\footnote{\ 
  There are two possible answers each to three independent questions:
  (1) Is the slant separating the two shapes left- or right-facing?
  (2) (resp., (3)) Is the left (resp., right) shape a 
  parallelogram or a trapezoid?
}
\]
In six of the cases,
we have $L_k(q) = L_k \leq \defm{\ol L_k} = L_k(p)$
for our usual choices of orbit representatives $\xip{k+1}(q)$.
In the cases
	 \parallelogram\mtrap\ and \wtrap\mtrap,
the usual representatives do not make this true,
but in these cases one can select 
a different smooth family of representatives $'\xip{k+1}(q)$
per the caution in \Cref{rmk:xi-explicit}
such that again $L_k \leq \ol L_k$.
(We suppress the details of this analysis.)
These groups are then connected by the commutative diagram
\[
\xymatrix{
\cdots \ar[r] & H_k \ar[d]& \ar[l] L_k \ar[r] \ar[d]& H_{k+1} \ar[d]&\ar[l] L_{k+1} \ar[r] \ar[d]&     \cdots\\
\cdots \ar[r] & \ol H_k & \ar[l] \ol L_k \ar[r] & \ol H_{k+1} &\ar[l] \ol L_{k+1} \ar[r] & \cdots \mathrlap,
}
\]
inducing a smooth map
\quation{\label{eq:ray-cylinder}
  T_q
  	=
  H_\alpha \,\ox_{L_\alpha}\, 
  H_{\alpha+1} \,\ox_{L_{\alpha+1}}\, 
  \cdots \,\ox_{L_{n-1}}\, 
  H_{n}/L_{n}
\,	 \lt  		\,
  \ol H_\alpha \,\ox_{\ol L_\alpha}\, 
  \ol H_{\alpha+1} \,\ox_{\ol L_{\alpha+1}}\, 
  \cdots \,\ox_{\ol L_{n-1}}\, 
  \ol H_{n}/\ol L_{n} 
  	= 
  T_p
}
for $q \in \g \less \{p\}$.
Our goal is to identify $\Psi\-(\g)$ with the mapping cylinder\footnote{\
  Recall that the \emph{mapping cylinder} of a continuous map $f\: X \lt Y$
  of topological spaces is the quotient of the 
  disjoint union of $X \x [0,1]$ and $Y$
  by the equivalence relation setting $(x,1) \sim f(x)$.
}
of this map.

Let us agree to write $\defm{{\theinterval} \cdot {G_k/H_k}}$ 
for the mapping cylinder of the
quotient map $G_k/H_k \longepi G_k/\ol H_k$
and similarly write $\defm{{\theinterval} \cdot {G_k/L_k}}$.
Then the continuity of $a_k(q)$ 
gives us a continuous map ${\theinterval} \cdot {G_k/L_k} \lt {\theinterval} \cdot {G_{k+1}/L_{k+1}}$
taking $\big[q,gL_k(q)\big] \lmt \big[q,\vp_k(g)a_{k+1}(q)\-H_k(q)\big]$,
and putting these 
together, 
we obtain a ``${\theinterval}\cdot$ variant'' of \eqref{eq:stair-system}.
We claim that the final pullback $\defm{\smash{\W^r_{\theinterval}}}$
of the staircase ending in this zig-zag is 
indeed homeomorphic to $\Psi\-( \g)$.
Indeed, $\smash{\Omega^r_{\theinterval}}$ is equipped 
with a natural map to ${\theinterval}$
and for each $q \in {\theinterval}$ 
there is a diffeomorphism $\defm{\w_q}\: \smash{\Omega^r_q} \lt \Psi\-(q)$.
As an iterated fiber product,
$\smash{\Omega^r_{q}}$ comprises
those lists $\big((g_k L_k(q)\big) \in \prod_{k=1}^n G_k/L_k(q)$
satisfying the equations $\vp_k(g_k)a_{k+1}\- H_{k+1} = g_{k+1}H_{k+1}$
for $1 \leq k \leq n-1$;
under this identification,
$\w_q$ is given by
$\big(g_k L_k(q)\big) \lmt \smash{\Ad^*_{\vp_n(g_n)} \xip{n}(q)}
$.
As $q \lmt \xip{n}(q)$ is continuous, 
it follows $\defm\w\: \smash{\big[q,\big(g_k L_k(q)\big)\big]} \lmt \w_q\big(g_k L_k(q)\big)$
is a continuous bijection $\smash{\Omega^r_{\theinterval}} \lt \Psi\-(\g)$.
Since $\smash{\Omega^r_{\theinterval}}$ is compact 
and $\Psi\-(\g)$ Hausdorff,
$\w$ is a homeomorphism.

As expected, replacing the maps ${\theinterval}\.G_k/L_k \lt {\theinterval}\.G_{k+1}/H_{k+1}$
with
\[
  \big[q,gL_k(q)\big] \lmt \big[q,a_{k+1}(q)\vp_k(g)a_{k+1}(q)\-H_k(q)\big]
\]
yields a space $\smash{\defm{\W^c_{\theinterval}}}$
homeomorphic to $\smash{\W^r_{\theinterval}}$.
To state the appropriate variant of \Cref{thm:explicit-description},
let~$\defm{{\theinterval}\.H_k}$ denote the mapping cylinder 
$[0,\e_\mu]\mnn \x \mnn H_k \,\union \,\{0\}\mnn \x\mnn \ol H_k$
of $H_k \longinc \ol H_k$, 
and similarly for~$\defm{{\theinterval}\.L_k}$.
Then as expected, we find a homeomorphism of $\W^c_{\theinterval}$ with
\[
\defm{T_{\theinterval}} \ceq
{\theinterval}\.H_\alpha \,\ox_{{\theinterval}\.L_\alpha}\, {\theinterval}\.H_{\alpha+1} \,\ox_{{\theinterval}\.L_{\alpha+1}}\, \cdots \,\ox_{{\theinterval}\.L_{n-1}}\, {\theinterval}\.H_{n}/{\theinterval}\.L_{n}
\mathrlap.\footnote{\ 
	One can then transfer over the smooth structure if one wishes.
}
\]

%
%

As we observed in \Cref{rmk:explicit-conversion},
for the maps $L_k \lt H_{k+1}$
determining $T_p$ for $p \in F$ do not depend on $p$,
so the $T_{q}$ are all the same for $t > 0$
and $T_{\theinterval}$
is itself a mapping cylinder,
that of the map 
\eqref{eq:ray-cylinder}.
Now we note that this whole discussion fixed, 
but did not depend upon the particular value of $\mu$,
but groups $H_k$, $\ol H_k$, $L_k$, and $\ol L_k$
depended only whether or not $t$ was $0$,
and again,
the maps between them,
induced by the $a_k(q)$, 
do not depend on $\mu$
which are continuous in $q$ within the region we are considering.
We may select a closed rectangular neighborhood 
$R_P = [\mu_0,\mu_1] \x [0,\e_P]$
of $p$ within $P$ in $\mu t$-coordinates,
where $\defm{\e_P}$ is the minimum value 
${\e_{\mu'}}$ for $\mu' \in [\mu_0,\mu_1]$.
Then varying $\mu$ in the above argument
and replacing each $e_\mu$ with $\e_P$,
we have homeomorphisms of $\Psi\-(R_P)$ with 
what we may call $\W^r_{R_P}$, $\W^c_{R_P}$, and $T_{R_P}$.
Again by the fact observed in \Cref{rmk:explicit-conversion}
that the maps $L_k \lt H_{k+1}$ do not depend on $q \in E$
and $\ol L_k \lt \ol H_{k+1}$ do not depend on $p \in F$,
we see $T_{R_P}$ is homeomorphic to $[\mu_0,\mu_1] \x T_{\theinterval'}$.
(To be sure,
as in \Cref{rmk:explicit-conversion},
the expressions for the maps $\W^r_{\theinterval} \lt T_{\theinterval}$
do indeed depend on $\mu$,
even though $\smash{\W^c_{\theinterval}}$, $\smash{T_{\theinterval}}$,
and the homeomorphism $\smash{\W^c_{\theinterval} \lt T_{\theinterval}}$
do not.)

Now we may equally well vary $P$ and hence $R_P$ 
by moving the $\lup k _j$ labeling the non-merging
components of the GZ pattern over $E$ and $F$
in such a way as to not collide.
By the independence of $T_p$ and $T_q$ from these choices
we have just observed,
we finally find the following.
\end{discussion}

\begin{theorem}\label{thm:ray}
For $U$ a small enough neighborhood of $p \in \ol E$,
we have homeomorphisms 
\[
\Psi\-(U) \homeoto \W^r_U \homeoto \W^c_U \homeoto T_U\homeoto (U \inter F) \x T_{\theinterval}
\mathrlap,
\]
where $T_{\theinterval}$ is the mapping cylinder of the map \eqref{eq:ray-cylinder}.
\end{theorem}


\subsection{Consistency and continuity of the determinant}\label{sec:consistency}

Now we can show the map $\Det$ of \Cref{sec:degeneration}
is well-defined and consistent.

\begin{proof}[Proof of \Cref{thm:det}]
By \Cref{rmk:face},
over an open face $E$ of $ \D$,
we may define a determinant map 
by
\[
	\Psi\-(E) 
		\isoto
	E \x T_q
		\longepi
	T_q
		\os\Det\lt
	T^{t(q)}\mathrlap,
\]
where we have seen the particular point $q \in E$ chosen does not affect
the map,
and by \Cref{rmk:ray}, 
for $F$ the interior of a face $\ol F \subn \ol E$
and $U \sub \ol E$ a sufficiently small neighborhood of $p \in F$,
we have $\Psi\-(U) \homeo (U \inter F) \x T_{\theinterval}$,
where $T_{\theinterval}$ is the mapping cylinder of the map $\defm f$ of \eqref{eq:ray-cylinder}.
Fixing $\mu$ with $\{p\} = \theinterval \inter F$ and a $q \in \theinterval \inter E$,
we claim there is a commutative square
as in \eqref{eq:det-square},
inducing a map of mapping cylinders 
from $T_{\theinterval}$ to the mapping cylinder $X_0|_{\theinterval}$ of $T^f$.
Granted this, composing through
\[
	\Psi\-(U) \lt (U \inter F) \x T_{\theinterval} \longepi T_{\theinterval} \lt X_0|_{U}
\]
gives a continuous, well-defined map $\Det_U$
agreeing with the previously defined $\Det_p$ for $p \in F$ and
$\Det_q$ for $q \in E$.
As $p \in F$ was arbitrary,
this defines $\Det$ continuously on $\Psi\-(E \union F)$.

To see the commutativity of \eqref{eq:det-square},
we note two combinatorial possibilities.
First, if the GZ pattern associated to $F$ arises from that of $E$ by joining
two components not meeting the top row,
then $L_n = \ol L_n$,
but the maps $\U(m^+) \+ \U(m^-) \lt \U(m^+ + m^-)$ on the 
components of $H_k$ corresponding to merging components induce
$\deg(g) \+ \deg(h) \lmt \det(g \+ h) = \det(g)\det(h)$.
Thus in this case the map $T^f\: T^{t(p)} \lt T^{t(q)}$ is of the form
$\id \x \mu \x \id$ where $\mu$ multiplies the $i\th$ and $(i+1)\st$ $\U(1)$ coordinates.
Second, if in the pattern for $F$, the $i\th$ component of the pattern of $E$
becomes joined with the top row, the map 
$H_n/L_n \lt \ol H_n/\ol L_n$ restricted to
the relevant components
is of the form $\U(m^+) \+ \U(m^-)/\U(m') \lt \U(m^+ + m^-)/\U(m')$ (or symmetrically swapping~$+$ and~$-$).
The effect on determinants is that the coordinate coming from the $i\th$ factor,
labeled in the pattern by $\l^+$ and corresponding to $\U(m^+)$ in the domain of this map,
drops out, so in this case
$T^f\: T^{t(p)} \lt T^{t(q)}$
projects out the $i\th$ coordinate.
Either way, the square commutes.


We now inductively extend the definition over all of $\D$.
To begin, 
for a vertex $p \in \D$,
we have $\Det_p$ defined on $\Psi\-(p)$ already by \Cref{thm:circle-factor}.
Suppose we have already proven $\Det$ is 
defined and continuous on the union $\defm{\Orb^{\leq m}} \sub \Orb$
of $\Psi$-preimages of faces of $\D$ 
of dimension $\leq m$,
in such a way that the definition over each open face is that of \Cref{sec:face}.
Each dimension-$(m+1)$ face~$\ol E$ is bounded by dimension-$m$ faces;
let $\ol F$ be one of them.
We can then use the first paragraph to define $\Det$
continuously on  $\Psi\-(E \union F)$.
We can equally well do this on $\Psi\-(E \union F')$ 
for any other face $\ol F'$ of $\ol E$.
If all these extensions over $E$ agree,
then we will have extended $\Det$ in a well-defined and continuous 
manner to $\Orb^{\leq m} \union \Psi\-(E)$,
and we can do this independently
for each $E$ of dimension $m+1$ to complete the induction.

To see the extensions agree,
note that the process of \Cref{sec:ray} 
depends on a choice of representatives $\xip{k+1}$,\footnote{\ 
  Indeed, in the situation \wtrap\mtrap\wtrap,
  specializing to \longrpara\wtrap\ 
  by letting the eigenvalues of \wtrap\mtrap\ approach one another
  and to \wtrap\longpara\ 
  by letting the eigenvalues of \mtrap\wtrap\ approach one another
  require different choices of $\xip{k+1}$.}
  but that by \Cref{rmk:xi-prime} and \Cref{rmk:xi-SO-explicit},
  these choices are related by $'\xip{k+1} = \vp_k(h) \xip{k+1}$ 
  for an $h = \Direct u_i \in \Direct_{i=1}^\ell \SU(d_i) \+ \SO(d_0)$.
  We have seen in \Cref{rmk:change-of-xi}
  and \eqref{eq:tensor-conversion}
that this coordinate change
takes the tensor-coordinate $a_k g_{k-1}\- g_k$ 
to $a_k g_{k-1}\- g_k h\-$
and does not affect the other coordinates.
This means $\Det$ as computed with respect to $'\xip{k+1}$ and $\xip{k+1}$
is the same: only the circle component of $\Det$ 
corresponding to the pattern component containing the \mtrap-shape
is affected by the transition, and this component of $\Det$ 
is multiplied by $\det(u_i) = 1$.
\end{proof}

%

%
%
%

\section{Algebro-topological invariants}\label{sec:invariants}

In this section we apply the results of \Cref{sec:product}
to homotopy groups and cohomology rings.

\subsection{Homotopy groups}\label{sec:homotopy}
From the balanced product presentation of \Cref{thm:explicit-description}
and the factorization of \Cref{thm:component-product}
into terms of the form \eqref{eq:description-v2}
and \eqref{eq:top-row-description-v2},
we can easily recover the first three homotopy groups 
of a Gelfand--Zeitlin fiber.


\bprop[Cho--Kim--Oh~{\cite[Prop.~6.13]{ChoKimOh2020}}]\label{thm:pi-1-unitary}
The fundamental group $\pi_1\big( \Psi\-(p)\big)$ of a unitary \GZ fiber
is free abelian of rank $t(p)$, the number of components of the \GZ pattern
not ending in the top row,
whereas $\pi_2\big(\Psi\-(p)\big)$ is trivial.
\eprop
\bpf
By \Cref{thm:circle-factor},
$\Psi\-(p)$ is the direct product of a torus $T^{t(p)}$
and a factor $\Det_p\-(1)$,
which by \Cref{thm:component-product} 
is a direct product of factors $Y_i$ of the form \eqref{eq:SU-factor}.
The result follows from the long exact homotopy sequences
of
$\prod \SU(m_j) \to \prod \SU(M_j) \to Y_i$.
%
\epf

The same sequence recovers $\pi_3(Y_i)$
as the cokernel of
$\pi_3\big(\mn \prod \SU(m_j) \lt \prod \SU(M_j)\big)$.
	Since 
	$\SU(m) \longinc \SU(n)$ 
	induces an isomorphism on $\pi_3 \iso \Z$
	for $n \geq m \geq 2$,
	one finds the following.

\begin{proposition}\label{thm:pi-3-U}
The third homotopy group $\pi_3\big(\Psi\-(p)\big)$ of a unitary \GZ fiber
is free abelian.
Its rank is the number of components
of the graph obtained from the \GZ pattern
by deleting rows of width $1$ from each component,
then deleting any components still meeting the top row.
\end{proposition}

Similar reasoning also applies in the case of an orthogonal \GZ fiber 
$\Psi\-(p)$.
Recall from \Cref{def:FSO} that 
$\Psi\-(p)$ is the direct product of a unitary part $F_\U$
that behaves like a unitary Gelfand--Zeitlin fiber 
and another factor $F_\SO$.

\begin{proposition}\label{thm:pi-3-O}
	The factor $F_\SO$ 
	splits as a direct product of a torus $T \iso \SO(2)^t$ and
	a space $Y$ with $\pi_1(Y) \iso (\Z/2)^s$
	and $\pi_2(Y) \iso \Z^f$
	and $\pi_3(Y) \iso \Z^g$,
		where $t$, $s$, and $f$ are determined
		from the white components of the associated \GZ pattern:
	$t$ is the number of local maxima of width $2$,
	$f$ is the number of local minima of width $2$,
	possibly including the top row,
	and $s$ is determined by 
	counting the remaining white components
	remaining after
	the union of all diamonds of width $2$
	is deleted,
	then any resulting component 
	containing $2$ or more vertices in the top row
	is also deleted.
	We determine $g$
	as a sum of numbers $h$
	given in terms of the description of factors $Y_i$
	as in \eqref{eq:description-v2} and \eqref{eq:top-row-description-v2}:
	if $N_4$ is the number of $j$ such that $M_j = 4$,
	while $n_{\geq 3}$ is the number with $m_j \geq 3$,
	then $h = r + N_4 - n_{\geq 3}$.
\end{proposition}
\begin{proof}
	We noted in \Cref{eg:diamonds}
	that each diamond of width $2$ contributes an $\SO(2)$ factor.
	It is then an exercise to examine each $Y_i$ separately
	using the long exact homotopy sequence of
	of the bundle
	$\prod \SO(m_j) \to \prod \SO(M_j) \to Y_i$.
	The key facts are that $\pi_2$ of a Lie group is trivial,
	$\pi_3$ of a simple Lie group is $\Z$,
and the subgroup inclusions
$\SO(m) \longinc \SO(n)$
induce an isomorphism on $\pi_1 \iso \Z/2$ for $m \geq 3$
and the surjection $\Z \longepi \Z/2$ for $m = 2$.
%
\end{proof}

\subsection{Cohomology of unitary Gelfand--Zeitlin fibers}\label{sec:cohomology-U}

In this section we compute the cohomology of the fibers of 
Gelfand--Zeitlin systems on unitary coadjoint orbits.
Our result is summarized in the following theorem.
	
\begin{theorem}\label{thm:main-U}
The integer cohomology ring of a unitary GZ fiber is isomorphic to the exterior algebra
\[
\phantom{|z_{2d-1,i}| = 2d-1\mathrlap.}
\ext[z_{2d-1,i}],\qquad
|z_{2d-1,i}| = 2d-1\mathrlap.
\]
where the generators $z_{2d-1,i}$ enumerate all \mtrap-shapes 
in the associated \GZ pattern with width $d$,
i.e., $d$ vertices on the long edge.
\end{theorem}

We demonstrate this result with an example before proceeding to the proof.

\begin{example} \label{eg:intro example cohom}
Enumerating \mtrap-shapes in the GZ pattern in Figure~\ref{fig:unitary-GZ-pattern},
we see that the integral cohomology ring of the associated GZ fiber is 
\[
\ext[z_{1,1},
	z_{1,2},
	z_{1,3},
	z_{1,4},
	z_{1,5},
	z_{1,6},
	z_{1,7},
	z_{3,1},
	z_{3,2},
	z_{3,3},
	z_{5,1},
	z_{5,2},
	z_{7,1}],
\qquad
|z_{m,j}| = m\mathrlap.
\]
Note that isolated vertices in the top row of the GZ pattern do not 
count as \mtrap-shapes of width $1$.
\end{example}

We now use the Serre exact sequence
compute the cohomology of the space~$F^{n+1}_1 \iso \Psi\-(p)$
in \Cref{fig:big-diagram}.
We will write $\defm{G_k} = \U(k)$
and $\defm{F^k} \ceq F^k_1$,
and $\defm{\iota}$ for the embedding $G_k/L_k \lt G_{k+1}/H_{k+1}$
of \Cref{thm:iotak}.
All cohomology will be taken with $\Z$ coefficients. 

\begin{lemma}\label{thm:unitary-collapse}
The Serre spectral sequence of each bundle
$F^{k+1} \to F^k$ collapses at $E_2$.
\end{lemma}
\begin{proof}
Let
$\defm{\chi_k}\: \U(k)/H_k \lt BH_k$ be a classifying map 
for the principal $H_k$-bundle $\U(k) \to \U(k)/H_k$,
and write 
 $\defm{j_k}\: F^k \lt F^k_k = \U(k)/H_k$ 
 for the composition of horizontal in \Cref{fig:big-diagram}. 
We prove by induction on $k$ that $\chi_k \o j_k$ 
induces the zero map in reduced cohomology. 
This is so for $k=1$ because $F^1$ is a point.
For the induction step, 
factor $j_{k+1}$ as the composite of 
$\defm{\wt\jmath_k}\: F^{k+1} \lt F^{k+1}_k = \U(k)/L_k$
and $\iota$
and let $\defm{\wt \chi_k}\: \U(k)/L_k \lt BL_k$ be a 
classifying map for the principal $L_k$-bundle $\U(k) \lt \U(k)/L_k$.
If we model $BL_k$ as $EH_k/L_k$,
then we have a map of $H_k/L_k$-bundles, 
shown in the bottom two rows of \eqref{eq:HkLk-universal-U}.
\quation{\label{eq:HkLk-universal-U}
\begin{aligned}
	\xymatrix@C=.5em@R=2em{						&
	G_{k+1}/H_{k+ 1} 	\ar[r]^(.55){\chi_{k+1}}			&
BH_{k+1}							\\
	F^{k+1}
				 \ar[r]^(.475){\wt{\jmath}_k}
				 \ar[ru]^{j_{k+1}}
				 \ar[d]									
									& 
	G_k/L_{k}		\ar[d]
				\ar[u]^\iota
				\ar[r]^(.525){{\wt\chi_k}}	 & 
				\ar[u]
	BL_{k}			\ar[d]				   \\
	F^{k}			\ar[r]_(.425){j_k}	
								&		
	G_k/H_{k}		\ar[r]_(.55){\chi_k}		& 
	BH_{k}
	}
\end{aligned}
}
Consider the induced diagram in reduced cohomology.
Since each
$\H B\U(m+1) \lt \H B\U(m)$ is a surjection,
so is $\H(BH_k) \lt \H(BL_k)$.
By the inductive assumption,
the bottom map $j_k^* \o \chi_k^*$ is trivial,
so it follows
$\wt\jmath_k^* \o \wt\chi_k^*$ is trivial and hence 
so is $j_{k+1}^* \o \chi_{k+1}^*$,
concluding the induction.

For the collapse,
note that
every exterior generator $\defm{z}$ of the ring $\H(H_k/L_k)$ 
transgresses in the spectral sequence of $BL_k \to BH_k$
to some $\tau z$, 
as one sees comparing the spectral sequences of the factors
$B\U(m) \os{=}\to B\U(m)$ and $B\U(m) \to B\U(m+1)$.
But the map of $H_k/L_k$-bundles on the bottom in \eqref{eq:HkLk-universal-U}
induces a map of {\SSS}s,
so $z$
transgresses to $(\chi_k \o j_k)^*\tau z = 0$
in the spectral sequence of $F^{k+1} \to F^k$.
\end{proof}


\begin{proof}[Proof of \Cref{thm:main-U}]
We show by induction on $k$ that
$
	\H(F^{k}) \iso \Tensor_{\ell=1}^{k-1} \H(H_\ell/L_\ell)
	$
as graded rings,
where $H_k/L_k$ were identified in \Cref{rmk:xi-explicit}.
For $k=1$, this is because the product is empty and $F^1$ is a point.
For the induction step, 
the Serre spectral sequence of $F^{k+1} \to F^k$ collapses at $E_2$
by \Cref{thm:unitary-collapse},
so that associated graded algebra of $\H(F^{k+1})$
is $ \H(F^{k}) \otimes \H(H_{k}/L_{k})$.
This is a free abelian group by the inductive hypothesis,
so there is no additive extension problem
and in particular no $2$-torsion.
The lifts of a set of exterior generators for $\H(H_{k}/L_{k})$ 
along the fiber restriction
$\H(F^{k+1}) \lt \H(H_{k}/L_{k})$ 
thus square to zero and 
hence generate an exterior subalgebra $
\defm\Lambda \iso \H(H_{k}/L_{k})$ within $\H(F^{k+1})$.
A basis of $\Lambda$ also forms a basis for $\H(F^{k+1})$ 
as an $\H(F^{k})$-module, by the collapse 
and absence of an additive extension problem,
so $\H(F^{k+1})$ is the tensor product of $\Lambda$ 
and $\H(F^k)$, concluding the induction.
\end{proof}

Knowing the cohomological structure 
also makes available a very natural homological description.
Recall that the K\"unneth theorem 
and the multiplication of a topological group $G$,
induce a graded ring structure on $H_*(G;k)$
for any choice of coefficient ring $k$,
%
%
called the \emph{Pontrjagin ring},
and a homomorphism $H \lt G$ makes $H_*(G;k)$ a module over $H_*(H;k)$.
This particularly holds for an inclusion $\U(m) \lt \U(M)$
and the quotient map $\U(m) \lt {*}$.
Using \Cref{thm:main-U}
and counting ranks, one finds the following.
%

\begin{proposition}\label{thm:homology-U}
Let a space $X$ be given
as on the left in 
\eqref{eq:description-v2} or
\eqref{eq:top-row-description-v2}
be given, writing $\U(q)$ for the trivial group $\U(0) = 1$
in the former case and for $K_n$ in the latter.
Then the projection of the fiber bundle
\[
\prod_{p=1}^{r-1} \U(m_p) \x \U(q)
\,\lt\,
\prod_{p=1}^r \U(M_p) 
\,\lt\,
	\U(M_1) \ox_{\U(m_1)}
	\cdots \ox_{\U(m_{r-1})} 
	\U(M_r)/\U(q) = X
\]
induces a surjection in homology and an identification
\[
H_*(X) \iso 
H_*\U(M_1) \ox_{H_*\U(m_1)}
\cdots \ox_{H_*\U(m_{r-1})} 
H_*\U(M_r) \ox_{H_*\U(q)} \Z
\mathrlap.
\]
\end{proposition}
\nd This is similar to the reasoning that allowed us to compute $\pi_3$
in \Cref{thm:pi-3-U}.

\subsection{Cohomology of orthogonal Gelfand--Zeitlin fibers}\label{sec:cohomology-O}

The computation of the cohomology 
of a \GZ fiber is more interesting in the orthogonal case.

\begin{notation}\label{def:zero-component}
	We extend the diagram of \Cref{fig:big-diagram}
	for an orthogonal \GZ fiber
	to $\a = 1$ by taking $G_k = H_k = L_k = \SO(1)$
	and setting $\defm{F^k} = F^k_1$.
	We write $\defm F = F^{n+1}$ for the \GZ fiber itself.
	We decomposing $H_k$ and $L_k$
	respectively as $H_k^\U \x H_k^\SO$
	and $L_k^\U \x L_k^\SO$
	with $\defm{H_k^\U}$ and $\defm{L_k^\U}$
	products of unitary groups and $H_k^\SO$ and $L_k^\SO$
	orthogonal.
\end{notation}

	Both $F_\U$ and $F_\SO$ themeselves
	admit tower descriptions as in \Cref{fig:big-diagram},
	given by taking
	$F_k^k$ and $F_k^{k+1}$
	to respectively be $\SO(k)/H_k^\U$ and $\SO(k)/L_k^\U$
	in the former case and 
	$\SO(k)/H_k^\SO$ and $\SO(k)/L_k^\SO$
	in the latter.
Because the quotients $H_k^\U/L_k^\U$ 
are products of odd-dimensional
spheres, \Cref{thm:main-U}
and Propositions
\ref{thm:circle-factor}, \ref{thm:homology-U}, and \ref{thm:pi-3-U} and their proofs
apply to give the following.

\begin{proposition}\label{thm:homotopy-U}
	Given an orthogonal \GZ fiber $F$,
	the factor $F_\U$
	is itself the product of a torus~$T$
	and a space $Y$ with $\pi_1(Y) = 0 = \pi_2(Y)$.
	The dimension of~$T$ is the number of~\mtrap-shapes
	of width~$1$ in the positive part of the \GZ pattern.
\end{proposition}

\begin{theorem}\label{thm:main-O-nonzero}
	Given an orthogonal \GZ fiber $F$,
	the cohomology ring
	$\H(F_\U;\Z)$ is
	an exterior algebra 
$
	\ext[z_{2d-1,i}]
$,
	where the generators $z_{2d-1,i}$ 
	of degree $2d-1$ are enumerated by
	black \mtrap-shapes of width~$d$
	in the \emph{positive part} of the \GZ pattern.
\end{theorem}

Since $\H(F_\U)$ is free abelian, 
the K\"unneth theorem
shows $\H(F)$ is isomorphic as a 
graded ring to $H(F_\U) \ox \H(F_\SO)$,
so we focus exclusively on $F_\SO$
for the rest of the section.

  We have seen in \Cref{rmk:Stiefel}
  that every Stiefel manifold is realized as a GZ fiber
and in \Cref{rmk:Stiefel}
  that $F_\SO$ is the total space of tower of bundles of real Stiefel manifolds
  $V_p = \SO(M_p)/\SO(m_p)$,
\begin{equation}\label{eq:iter-Stiefel}
\begin{aligned}
\xymatrix{
	V_{r} \ar[d]^{i_r} & 
	V_{r-1}\ar[d]^{i_{r-1}}&
	V_{r - 2}\ar[d]^{i_{r-2}}&
	\cdots &
	V_3 \ar[d]^{i_3} &
	V_2 \ar[d]^{i_2}	\\
	F_\SO \ar[r]&
	F^{\ell_{r-1}}\ar[r] &
	F^{\ell_{r-2}} \ar[r]&
	\cdots \ar[r]&
	F^{\ell_3} \ar[r] & 
	F^{\ell_2} \ar[r] &
	V_1 \mathrlap,
}
\end{aligned}
\end{equation}
which will play the same role for  $\H(F_\SO)$
that \Cref{fig:big-diagram} did for the unitary \Cref{thm:main-U}.

\begin{lemma}\label{thm:orthogonal-collapse}
  In the tower \eqref{eq:iter-Stiefel},
  the integral cohomology \SSS of each bundle 
  $V_p \to F^{\ell_p} \to F^{\ell_{p-1}}$
  collapses.
\end{lemma}

Thus most of what we have to say about 
the cohomology of $F_\SO$ is a consequence 
of facts about the cohomology of a Stiefel manifold, 
which we now discuss.

\begin{examples}\label{ex:SO}
  Following \Cref{ex:diamond,ex:Stiefel},
  we see any $\SO(M)$ and more generally any $V_{M-m}(\R^M) = \SO(M)/\SO(m)$
arises as a \GZ fiber,
so a general description of the ring $\H(F_\SO)$ is at least as hard
as that of $\H V_k(\R^m)$.
	The simplest non-sphere Stiefel manifolds, $V_2(\R^m)$
	already	give an example of the failure of the analogue
	of \Cref{thm:unitary-collapse}
	in the orthogonal case.
	Indeed, the {\SSS} of the associated sphere bundle
	$
	\frac{\SO(m-1)}{\SO(m-2)} \to \frac{\SO(m)}{\SO(m-2)} \to\frac{\SO(m)}{\SO(m-1)} 
	$
	supports a differential
	taking $z = 1 \otimes [S^{m-2}]$
	to $\chi(S^{m-1})$ times $s = [S^{m-1}] \ox 1$,
	which is nonzero when $m$ is odd.
It follows that
	the cohomology ring 
	is  $\H\big(V_2(\R^m)\big) \iso \ext[z,s]$
	when $m$ is even,
	for any coefficients,
	whereas when $m = 2k+1$ is odd,
	the cohomology rings over $\Zf$ and $\Z$ are respectively
	$ \ext[zs] $
	and\,
	$\sfrac{\Z[s] \ox \ext [z]}{(2s,s^2,sz)} $.

	For another relatively simple example,
	the cohomology rings of $\SO(4)$
	over $\Zf$ and $\Z$ are respectively
	$\ext [z,q]	$
	and
	$\sfrac{\Z[s] \ox \ext [z,q]}{(2y,y^2,sq)} $,
	where $|y| = 2$ and $|z| = |q| = 3$.
\end{examples}

	\v{C}adek, Mimura, and Van\v{z}ura~\cite{cadekmimuravanzura}
	have found a presentation of $\H V_k(\R^M)$
	which we will summarize some only some aspects of \Cref{thm:HV},
	as it is rather complicated.\footnote{\ 
	  Even integral cohomology ring of $\SO(k)$ is involved,
	  and has been known in full only since
	  the late 1980s~\cite{SOSpinPittie}. 
	}

\begin{notation}\label{thm:mod-2}
In the Bockstein long exact sequence of cohomology groups of a space $X$
arising from the short exact sequence
$0 \to \Z \to \Z \to \Z/2 \to 0$ of coefficient rings,
we write $\defm{\beta}\: H^{*-1}(X;\F_2) \lt \H(X)$
for the connecting map
and $\defm\rho\: \H(X) \lt \H(X;\F_2)$ for the reduction map,
which is a ring homomorphism.
The first Steenrod square is 
$\defm{\Sq^1} = \rho\b\: H^{*-1}(X;\F_2) \lt H^{*}(X;\F_2)$.
We also write $\defm{\oH(X)} = \H(X)/\,\mnn\mathrm{torsion}$.
\end{notation}


\begin{definition}\label{def:simple}
	Let $k$ be $\Z$ or a finite quotient thereof,
	$A$ a commutative graded algebra,
	and $M$ a submodule of $A$ free over $k$.
	We say $M$ \defd{admits a simple system of generators} 
	if there are elements $v_j \in M$
	such that
	the distinct monomials $v_{j_1} \cdots v_{j_\ell}$ 
	of degree $0$ or $1$ in each $v_j$
	form a $k$-basis of $M$.\footnote{\ 
		For example, the polynomial ring $\Z[x]$
		admits the simple system of generators $\smash{x^{2^j}}$.
	If $A$ itself admits a simple system of generators,
	then its ring structure is determined entirely
	by the existence of the simple system of generators and
	the values of the squares $\smash{v_j^2}$.
	}
\end{definition}

\begin{proposition}[\v{C}adek--Mimura--Van\v{z}ura%
								~{\cite{cadekmimuravanzura}}]%
								\label{thm:HV}
	The cohomology ring $\H\big(\SO(M)/\SO(m);\F_2\big)$
	admits a system of one simple generator $\defm{v_i}$ 
	of each degree $i \in [m,M-1]$
	with $v_i^2 = v_{2i}$ for $2i < M$ and $v_i^2 = 0$ for $2i \geq M$
	and $\Sq^1 v_{2i} = v_{2i+1}$ for $2i + 1 < M$.
	All torsion in $\H\big(\SO(M)/\SO(m)\big)\mn$ is $2$-torsion
	and the torsion ideal has a canonical free abelian complement
	admitting a simple system of generators.
	This simple system 
	consists of 
	the following nontorsion generators indexed by degree:
	\[
	\defm{q_{4j-1}}\quad \mbox{whenever } [2j-1,2j+1] \sub [m,M],
	\qquad\quad
	\defm{u_m} \quad \mbox{if $m$ is even},
	\qquad\quad
	\defm{\chi_{2\ell-1}} \quad\mbox{if } M = 2\ell\mathrlap.
	\]
	The squares of these generators are $2$-torsion,
	and particularly the squares $u_m^2$ and $\chi_{2\ell-1}^2$ are $0$.
	
	One has $\rho \chi_{2\ell-1} = v_{2\ell-1}$.
	Given $N < M$, the 
	maps induced by the inclusion $\SO(N)/\SO(m) \lt \SO(M)/\SO(m)$
	are surjective on the $v_i$, the~$q_{4j-1}$, and~$u_m$,
	but~$\chi_{2\ell-1}$ is not in the image.
	When $M = 2\ell$ is even and $M > n > m$, 
	the map induced in cohomology 
	by the quotient map $\SO(M)/\SO(m) \lt \SO(M)/\SO(n)$
	does preserve $\chi_{2\ell-1}$.
	The transgression of $\chi_{2\ell-1}$
	in the \SSS of the universal principal $\SO(2\ell)$-bundle
	is a non-torsion class $\defm{e_{2\ell}} \in H^{2\ell}B\SO(2\ell)$,
	the \defd{Euler class},
	which is annihilated by restriction along the maps 
	$B\SO(N) \lt B\SO(2\ell)$ for $N < 2\ell$.
\end{proposition}

We will show that, similarly,
both $\oH(F_\SO)$ and $\H(F_\SO;\F_2)$
admit simple systems of generators
and $\H(F_\SO)$ has only $2$-torsion.
To state the result precisely 
it will help to expand our graphical lexicon.

\begin{notation}\label{def:traps2}
	We refer to as a \defd{\MMtrap}-shape 
	any three consecutive rows $k$, $k+1$, $k+2$ in 
	an orthogonal \GZ pattern
	such that the full subgraph on the white components
	in these rows is made up of a
	\mtrap-shape in rows $k$ and $k+1$ and an
	\mtrap-shape in rows $k+1$ and $k+2$ .
	For either of the shapes \mtrap{} and \MMtrap, 
	we will use subscript labels 
	to refer to the width $\ell$, counted in vertices, 
	of the bottom row of the shape,
	or the parity $\defm{\bar \ell} = \ell \mod 2$ of this width.
	
	We refer to as a (white) \defd{\hex}-shape the full subgraph of a white component
	of an orthogonal \GZ pattern bounded by (and including) two rows
	whose widths are consecutive local maxima.
	If the top row is of locally maximum width,
	it \emph{does not count} as a hexagon.
	To a \hex-shape of width $M$ vertices at its longest row
	and $m$ vertices at its top row, 
	we assign the 
	\defd{associated Stiefel manifold} $\SO(M)/\SO(m) = V_{M-m}(\R^M)$,
	visible as a fiber in \eqref{eq:iter-Stiefel}.
\end{notation}

\begin{examples}
	There are two 
	\MMtrap-shapes in \Cref{fig:SO4},
	in the upper half of the white diamond,
	and precisely one \MMtrap-shape
	in each member of the family in 
	\Cref{fig:V2R3,fig:V2R4,fig:V2R5}
	yielding GZ fibers $V_2(\R^m)$.
	The pattern in \Cref{fig:SO4}
	has a~\defm{\mtrap$_{\bar 0}$}-shape,
	more specifically a \mtrap$_{4}$-shape,
	in rows $4$ and $5$,
	and a \MMtrap$_3$-shape above in rows $5$, $6$, $7$.
	 In each of the patterns in \Cref{fig:Stiefel}
	 and in \Cref{fig:orthogonal-GZ-pattern}
	there is precisely one \hex-shape,
	and in \Cref{fig:white-large} there are five.
\end{examples}

\begin{theorem}\label{thm:main-O}
	All torsion in the cohomology of $F_\SO$
	(\Cref{def:zero-component})
	is $2$-torsion.
	If $\defm{V_p}$ are the Stiefel manifolds associated
	to the white \hex-shapes of the \GZ-pattern in \Cref{def:traps2}
	and \eqref{eq:iter-Stiefel},
	one has graded group isomorphisms
	\[
	  {\H(F_\SO)} \iso \H\Big(\prod_{p=1}^r V_p\Big)
	  \qquad\mbox{and}\qquad
	  \H(F_\SO;\F_2) \iso \Tensor_{p=1}^r \H(V_p;\F_2)
	\] 
and a graded ring isomorphism
	\[
	  \phantom{
			  |q_{4\smash{m'_j}-1}| = 4\smash{m'_j}-1\mathrlap,
  }
  \oH(F_\SO) \iso 
  \Tensor \oH(V_p) \iso
	  \ext[z_{s_i},q_{4\smash{s'_j}-1}], 
	  \qquad
			  |z_{s_i}| = s_i, \quad 
			  |q_{4\smash{s'_j}-1}| = 4\smash{s'_j}-1\mathrlap,
	\]
	where the generators $q_{4\smash{s'_j}-1}$
	are indexed by the white \MMtrap$_{2\smash{s'_j}+1}$-shapes
	in the \GZ pattern
	and the generators $z_{s_i}$ are indexed by the
	set of those white \mtrap$_{s_i+1}$-shapes 
	not contained in any \MMtrap$_{\bar 1}$-shape.\footnote{\
		For ease of presentation, 
		we allow even-dimensional exterior generators,
		meaning central elements squaring to $0$.
	}
\end{theorem}

Some, but not all, of the multiplicative structure of the $\H(V_p)$
is visible in $\H(F_\SO)$. 
If we are willing to discard multiplicative information altogether,
we have an additive characterization that is even simpler.

\begin{corollary}
	As a graded group, $\H(F_\SO)$ is isomorphic to the cohomology
	of a direct product of Stiefel manifolds $V_2(\R^{2s_j'+1})$
	indexed by \MMtrap$_{2s_j'+1}$-shapes
	and of spheres $S^{s_i}$ indexed by \mtrap$_{s_i+1}$-shapes not contained
	in any \MMtrap$_{\bar 1}$-shape.
\end{corollary}

Before diving into the proof, we give one last example.

\begin{example}
	Consider an orthogonal GZ pattern with white component
	outlined in Figure \ref{fig:white-large}, 
	where rows of locally maximal or minimal width are labeled by their widths.
	One has an expression $F^{22}_\SO \iso X \ox_{\SO(1)} Y = X \x Y$,
	where $X$ and $Y$ are the spaces corresponding to the subpatterns
	joined at the pinch point of width $1$.
	
	Counting \MMtrap$_{\bar 1}$-shapes,
	we find a \MMtrap$_3$-shape (a triangle)
	contributing a $V_2(\R^3) = \SO(3)$ factor
	in cohomology and a disjoint \mtrap$_3$-shape
	contributing an $S^2$ factor.
	Counting all \mtrap-shapes, we 
	find two \mtrap$_3$'s and an \mtrap$_2$,
	so $\H(X;\F_2)$ admits a simple system of generators
	of orders $2$, $2$, $1$.
	For the ring structure, the balanced product description
	\[
	X \iso \SO(1) \ox_{\SO(1)} \SO(3) \ox_{\SO(2)} \SO(3) / \SO(1)
		\iso \SO(3) \ox_{\SO(2)} \SO(3)\mathrlap.
	\]
	gives the sequence of Stiefel manifolds 
	$ \SO(3) /{\SO(2)} \iso S^2$ and $\SO(3) / \SO(1) \iso \SO(3)$
	corresponding to the \hex-shapes.
	In this case the diffeomorphism $g \ox h \lmt \big(g\SO(2),gh\big)$
	gives us $S^2 \x \SO(3)$.

	\begin{figure}[h]
		\centering
		\begin{tikzpicture}[scale =.25,line cap=round,line join=round,>=triangle 45,x=1cm,y=1cm,scale=1.5, every node/.style={scale=1}]
		\begin{scriptsize}
		
		\draw [line width=0.5mm] (0/2,-2*.9)--(2/2,0*.9);
		\draw [line width=0.5mm] (1/2,1*.9)--(2/2,0*.9);
		\draw [line width=0.5mm] (1/2,1*.9)--(2/2,2*.9);
		\draw [line width=0.5mm] (0/2,4*.9)--(2/2,2*.9);
		\draw [line width=0.5mm] (0/2,4*.9)--(7/2,11*.9);
		\draw [line width=0.5mm] (6/2,12*.9)--(7/2,11*.9);
		\draw [line width=0.5mm] (6/2,12*.9)--(7/2,13*.9);
		\draw [line width=0.5mm] (5/2,15*.9)--(7/2,13*.9);
		\draw [line width=0.5mm] (5/2,15*.9)--(6/2,16*.9);
		\draw [line width=0.5mm] (2/2,20*.9)--(6/2,16*.9);
		\draw [line width=0.5mm] (2/2,20*.9)--(-2/2,20*.9);
		
		\draw [line width=0.5mm] (0/2,-2*.9)--(-2/2,0*.9);
		\draw [line width=0.5mm] (-1/2,1*.9)--(-2/2,0*.9);
		\draw [line width=0.5mm] (-1/2,1*.9)--(-2/2,2*.9);
		\draw [line width=0.5mm] (0/2,4*.9)--(-2/2,2*.9);
		\draw [line width=0.5mm] (0/2,4*.9)--(-7/2,11*.9);
		\draw [line width=0.5mm] (-6/2,12*.9)--(-7/2,11*.9);
		\draw [line width=0.5mm] (-6/2,12*.9)--(-7/2,13*.9);
		\draw [line width=0.5mm] (-5/2,15*.9)--(-7/2,13*.9);
		\draw [line width=0.5mm] (-5/2,15*.9)--(-6/2,16*.9);
		\draw [line width=0.5mm] (-2/2,20*.9)--(-6/2,16*.9);

		\draw [line width=0.5mm] (-2/2,0*.9)--(2/2,0*.9);
		\draw [line width=0.5mm] (-1/2,1*.9)--(1/2,1*.9);
		\draw [line width=0.5mm] (-2/2,2*.9)--(2/2,2*.9);
		\draw [line width=0.5mm] (-7/2,11*.9)--(7/2,11*.9);
		\draw [line width=0.5mm] (-6/2,12*.9)--(6/2,12*.9);
		\draw [line width=0.5mm] (-7/2,13*.9)--(7/2,13*.9);
		\draw [line width=0.5mm] (-6/2,14*.9)--(6/2,14*.9);
		\draw [line width=0.5mm] (5/2,15*.9)--(-5/2,15*.9);
		\draw [line width=0.5mm] (-6/2,16*.9)--(6/2,16*.9);
		\draw [line width=0.5mm] (-4/2,18*.9)--(4/2,18*.9);
		\draw [line width=0.5mm] (-2/2,20*.9)--(2/2,20*.9);
		
		\path (0/2,-2*.9)		node[label=right:{$1$}] {};
		\path (2/2,0*.9)		node[label=right:{$3$}] {};
		\path (1/2,1*.9)		node[label=right:{$2$}] {};
		\path (2/2,2*.9)		node[label=right:{$3$}] {};
		\path (0/2,4*.9)		node[label=right:{$1$}] {};
		\path (7/2,11*.9)		node[label=right:{$8$}] {};
		\path (6/2,12*.9)		node[label=right:{$7$}] {};		
		\path (7/2,13*.9)		node[label=right:{$8$}] {};		
		\path (6/2,14*.9)		node[label=right:{$7$}] {};		
		\path (5/2,15*.9)		node[label=right:{$6$}] {};	
		\path (6/2,16*.9)		node[label=right:{$7$}] {};
		\path (4/2,18*.9)		node[label=right:{$5$}] {};
		\path (2/2,20*.9)		node[label=right:{$3$}] {};
		\end{scriptsize}
		\end{tikzpicture}
		\caption{Outline of a white component of a GZ pattern on a non-regular coadjoint orbit of $\SO(22)$;
			the numerical labels are row width
			in number of vertices.}
		\label{fig:white-large}
	\end{figure}
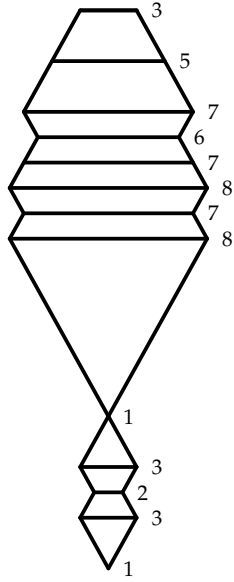
	
	Similarly, reading up from the bottom of the other component $Y$
	for \hex-shapes,
	we find the Stiefel manifolds $S^7 \iso \SO(8)/\SO(7)$,\ 
	$V_2 = \SO(8)/\SO(6)$, and $V_3 = \SO(7)/\SO(3)$.
	In the balanced product description
	\[
	Y	\iso \SO(1) \ox_{\SO(1)} \SO(8) \ox_{\SO(7)} \SO(8) \ox_{\SO(6)} \SO(7)/\SO(3)\mathrlap,
	\]
	the first Stiefel manifold $S^7$ splits off 
	to leave $\SO(8) \ox_{\SO(6)} \SO(7)/\SO(3)$,
	a $V_3$-bundle over $V_2$.
	Reading the remaining two \hex-shapes,
	corresponding to $V_2$ and $V_3$,
	for \MMtrap$_{\bar 1}$-shapes,
	we find a \MMtrap$_7$ 
	and a \MMtrap$_5$,
	respectively giving factors $\H V_2(\R^7)$
	and $\H V_2(\R^5)$.
	Then looking for \mtrap-shapes not contained in these,
	one encounters,
	reading from the bottom, 
	two \mtrap$_8$'s and a \mtrap$_7$,
	respectively contributing two $S^7$ factors and an $S^6$.
	Thus $Y$ has the cohomology of
	$S^7 \x S^7 \x S^6 \x V_2(\R^7) \x V_2(\R^5)$
	additively over $\Z$ and $\F_2$, 
	and multiplicatively over $\Z$ 
	after quotienting out $2$-torsion.
\end{example}

\begin{proof}[Proof of \Cref{thm:orthogonal-collapse}]
A result of Borel~\cite[Prop.~8.2]{borelthesis} 
states that given a bundle
$F \to E \to B$ of
compact CW complexes
such that all torsion in $\H(F)$ and $\H(B)$ is $2$-torsion
and the rational and mod-$2$ {\SSS}s collapse,
the integral spectral sequence also collapses
without an extension problem,
so that as graded groups, one has
$\H(E) \iso E_\infty = E_2 \iso \H(B \times F)$.
We will show that each \SSS collapses over $\Q$ and $\F_2$,
and then it will follow by induction 
starting with \Cref{thm:HV} and the K\"unneth theorem
that all torsion in each $\H(F^{\ell_p})$ is $2$-torsion
and each \SSS collapses integrally.

Fix $p \in [1,r]$.
Letting $N$ be the largest of the $M_t$ for $t \in [1,p]$, 
we 
observe 
that $\prod_t \SO(M_t)$ is a subgroup of $\SO(M)^{p}$,
where $M_{\ell_r} = 2\ell$,
and that the iterated group multiplication
$(g_1,\ldots,g_p) \lmt g_1 \cdots g_p$ on $\SO(N)$
restricts to a map $\prod_t \SO(M_t) \lt \SO(N)$.
By the definition of the balanced product, 
this multiplication descends to a well-defined map
\eqn{
	\defm{\mu_p}\: F^{\ell_p} =
\SO(M_1) \ox_{\SO(m_1)}
\SO(M_2) \ox_{\SO(m_2)} \cdots \ox_{\SO(m_{p-1})} 
\SO(M_p)/\SO(m_p)
	&\lt \SO(N)/\SO(m_p)\mathrlap,\\
	g_1 \ox \cdots \ox g_p\SO(m_p) &\lmt g_1 \cdots g_p\SO(m_p)\mathrlap.
}
It is clear that 
following 
the fiber inclusion
\eqn{
\defm{i_p}\: \SO(M_p)/\SO(m_p) = V_p 	&\lt  F^{\ell_p},\\
g_p\SO(m_p) 				&\lmt 1 \ox \cdots 1 \ox g_p\SO(m_p).
}
with $\mu_p$ 
gives the standard injection
$\SO(M_p) / {\SO(m_p)} \lt \SO(N) / {\SO(m_p)}$
discussed in \Cref{thm:HV}.
Thus each generator $v_i$
lies in the image of $i_p^*$, 
hence $i_p^*$ is surjective onto $\H(V_p;\F_2)$
and the mod-$2$ \SSS collapses.
This already gives the mod-$2$ isomorphism.

To prove the rational collapse,
note each $q_{4j-1}$ is also in the image of $i_p^*$,
as is $\smash{u_{m_p}}$ when $m_p$ is even.
It remains only to 
find a lift in $H^{2s-1}(\smash{F^{\ell_p}};\Q)$ of
$\chi_{2s-1} \in H^{2s-1}\big(\SO(M_p)/\SO(m_p);\Q\big)$ 
if $M_p = 2s$ is even.
Abbreviate the balanced product of the first $p-1$ factors 
in \eqref{eq:top-row-description-v2} as $\defm X$,
so that $F^{\ell_p}$ is $X \ox_{\SO(m_{p-1})} \SO(2s)/\SO(m_p)$.
Since the actions of $\SO(m_{p-1})$ and $\SO(m_p)$ we
are quotienting by are free,
$F^{\ell_p}$
is homotopy equivalent to the homotopy orbit space
\[
\defm{F'} 
	\ceq
X 
	\ox_{\SO(m_{p-1})} 
\big(E\SO(2s) \x \SO(2s)\ox_{\SO(m_p)} E\SO(2s)\big) 
\mathrlap,
\]
where the action of $\SO(m_{p-1})$ on 
the right-hand factor
is by $g\.(e,h \ox e') \ceq (eg\-,gh \ox e')$.
We can then replace $V_p$ by 
$\defm{V'} 
	\ceq 
\SO(2s) \ox_{\SO(m_p)} E\SO(2s)$,
replace 
$i_p$ by the obvious map 
$\defm{i'}\:
h \ox e' \lmt x_0 \otimes (e_0,h \ox e')$
for fixed basepoints $x_0 \in X$ and $e_0 \in E\SO(2s)$,
and replace the quotient map $\SO(2s) \lt V_p$ by 
$\defm{\pi'}\: h \lmt h \ox e_0$.

There are natural quotient maps 
from $F'$ to 
\[
\defm{F''} \ceq X \ox_{\SO(m_{p-1})} 
\big(E\SO(2s) \x \SO(2s) \ox_{\SO(2s-1)} E\SO(2s)\big)
\]
and from $V'$
to $\defm {V''} \ceq \SO(2s) \ox_{\SO(2s-1)} E\SO(2s) \hmt S^{2s-1}$.
Dropping the $X$ coordinate in $F''$ yields a natural
quotient map $\defm\varpi$
to 
$\defm E 
	\ceq
E\SO(2s)\ox_{\SO(m_{p-1})}\SO(2s) \ox_{\SO(2s-1)} E\SO(2s)$.
These maps then fit together as follows:
\[
\xymatrix{
	\SO(2s) 	\ar[r]^(.535){{\pi'}}			&
	V' 		\ar[r]^(.45){\defm\pi}	
			\ar[d]_{i'} 				&
	V''	 	\ar[d]^{\defm{i}}	
	\\ 							&	
	F'	 	\ar[r]_{\defm{\varpi'}}			&
	F'' 		\ar[r]_{\varpi}				&
	E\mathrlap.
}
\]
To lift $\chi_{2s-1} \in H^{2s-1}(V_p;\Q) \iso H^{2s-1}(V';\Q)$
to $F^{\ell_p} \hmt F'$, 
it is enough by commutativity of the square
to find a lift to $F''$ along $i \o \pi$.
Now by \Cref{thm:HV},
the fundamental class of $V'' \hmt S^{2s-1}$
is the unique lift of $\chi_{2s-1}$
under $\pi^*$ and also the unique lift of 
$\defm{\wt\chi} = (\pi')^*\chi_{2s-1} \in H^{2s-1}\big(\SO(2s);\Q\big)$
under $(\pi \o \pi')^*$,
so it will be enough to see $\wt \chi$
lies in the image of the fiber restriction 
along $\varpi \o i \o \pi \o \pi'$.

For this, we note that this map
is a fiber inclusion of the
nonprincipal $\SO(2s)$-bundle 
$E \lt B\SO(m_{p-1}) \x B\SO(2s-1)$
with projection 
$e \ox g \ox e' \lmt \big(e\SO(m_{p-1}),\SO(2s-1)e'\big)$.
We will show that $\wt\chi$ transgresses to $0$
in the \SSS of this bundle,
so that it lies in the image of the fiber restriction.
Following Eschenburg~\cite{eschenburg1992biquotient}, 
consider the following square:
\[
\xymatrix{
	E =\dsp E\SO(2s) \ox_{\SO(m_{p-1})} \SO(2s) \ox_{\SO(2s-1)} E\SO(2s)
	\ar[r]^(.475){\defm{\wt\kappa}}\ar@<1em>[d]&
	\dsp E\SO(2s)\ox_{\SO(2s)} \SO(2s) \ox_{\SO(2s)} E\SO(2s) 
	\hmt B\SO(2s)
	\ar[d]\\
	{\phantom{ E = {} }}
	B\SO(m_{p-1}) \x B\SO(2s-1) \ar[r]_{\defm\kappa} &
	B\SO(2s) \x B\SO(2s)\mathrlap.
}
\]
The right vertical map
$e \ox g \ox e' \lmt \big(\SO(2s)e,e'\SO(2s)\big)$
is the projection of another $\SO(2s)$-bundle
and by inspection the
horizontal map
$e \ox g \ox e' \lmt e \ox g \ox e'$
makes the square a map of $\SO(2s)$-bundles
such that 
the map $\kappa$ of base spaces is up to homotopy that
induced functorially by the subgroup inclusions.
The projection of the right bundle is,
up to homotopy, the diagonal map on $B\SO(2s)$,
which induces the cup product in cohomology,
and so in the \SSS of the right bundle,
$\wt\chi \in E_2^{0,2s-1} 
		\iso H^{2s-1}\big(\SO(2s);\Q\big)$
		transgresses to 
the image in $E_{2s}^{2s,0}$ 
of the element 
$1 \ox e_{2s} - e_{2s } \ox 1$
of $E_2^{2s,0} 
	\iso 
	\Direct_{a+b = 2s} H^a\big( B\SO(2s);\Q\big) \ox_\Q H^b \big(B\SO(2s);\Q\big)$.
The map of fiber bundles induces a map 
$\smash{\big(E_r(\wt\kappa)\big)_r}$
of spectral sequences,
so in the left spectral sequence, 
$\wt\chi$
transgresses to the $E_{2s}^{2s,0}(\wt\kappa)$-image of this element.
But 
$e_{2s}$ lies in the kernel of 
$\H \big(B\SO(2s);\Q\big) \lt \H \big(B\SO(t);\Q\big)$ for $t < 2s$
by \Cref{thm:HV},
and $2s-1,m_{p-1} < M_p = 2s$,
so $\smash{E_2^{2s,0}(\wt\kappa) }
	= 
\smash{	H^{2s}\big(\kappa;H^0(\id_{\SO(2s)};\Q)\big)}$
annihilates this image 
and hence $\wt\chi$ transgresses to $0$ as claimed.
It is thus the image of some $\defm x$ in $H^{2s-1}(E;\Q)$
and our desired lift in $H^{2s-1}(F^{\ell_p};\Q)$ is $(\varpi\o\varpi')^*x$.
\end{proof}

\begin{proof}[Proof of \Cref{thm:main-O}]
  The integral statement follows from Borel's
  theorem cited in the proof of \Cref{thm:orthogonal-collapse}.
  For the statement about $\smash{\oH(F_\SO)}$,
 we show
 inductively that each $\smash{\oH(V_p)}$ is exterior.
 Particularly, $\smash{\oH(V_1)}$ is by \Cref{thm:HV},
 so inductively assume the same of $\smash{\oH(F^{\ell_{p-1}})}$.
The \SSS of 
$V_p \to F^{\ell_p} \to F^{\ell_{p-1}}$
modulo torsion collapses
at $E_2 = \oH(F^{\ell_{p-1}}) \ox \oH(V_p)$,
so we need only to lift a system of exterior generators from $\oH(V_p)$.
For these we 
take the various images $\smash{\mu_p^*q_{4j-1}}$, 
and additionally $(\varpi\o\varpi')^*x$ if $M_p$ is even\footnote{\ 
	Since there is no longer any $2$-torsion,
	for the odd-degree lifts we could really have taken anything.
} 
and
$\mu_p^*u_{m_p}$
if $u_{m_p}$ is even, 
which works because $u_{m_p}^2 = 0$ in $H^{2m_p}\big(\SO(N)/\SO(m_p)\big)$.

It remains to enumerate these generators by trapezoids.
This follows on examining each $V_p = \SO(M_p)/\SO(m_p)$ separately.
Write $\defm{M'}$ for $M_p$ itself
if it is odd and otherwise for $M_p-1$,
and similarly $\defm{m'}$ for $m_p$ itself
if it is odd and otherwise for $m_p+1$.
Then there is a bundle tower
\[
\xymatrix@C=1.25em@R=2em{
	\frac{\SO(m')}{\SO(m_p)} 			\ar[d]& 
	\frac{\SO(m'+2)}{\SO(m')}			\ar[d]&
	\frac{\SO(m'+4)}{\SO(m'+2)}			\ar[d]&
	\cdots 							&
	\frac{\SO(M'-2)}{\SO(M'-4)}			\ar[d]&
	\frac{\SO(M')}{\SO(M'-2)}			 \ar[d]\\
	V_p						\ar[r]& 
	\frac{\SO(M_p)}{\SO(m')}			\ar[r] & 
	\frac{\SO(M_p)}{\SO(m'+2)}			\ar[r]& 
	\cdots 						\ar[r]& 
	\frac{\SO(M_p)}{\SO(M'-4)}			\ar[r] & 
	\frac{\SO(M_p)}{\SO(M'-2)}			\ar[r] & 
	\frac{\SO(M_p)}{\SO(M')} 			 \mathrlap.
}
\]
Borel showed~\cite[Prop.~10.4]{borelthesis}
that the cohomology of $V_p$
and of the direct product of the fibers and the base
in this diagram are isomorphic.
Quotienting out torsion,
the cohomology of each fiber and the base space
become exterior, on generators,
displayed respectively from left to right,
\[
\phantom{,,}
\mathllap(u_{m_p},) \qquad\quad
q_{2m'+1}, \qquad\quad
q_{2m'+5}, \qquad\quad
\cdots \qquad\quad
q_{2M'-7}, \qquad\quad
q_{2M'-3}, \qquad\quad
(\chi_{M_p-1}\mathrlap{,)}
\]
where $u_{m_p}$ is contributed if and only if the last fiber is nontrivial,
which is to say, $m_p$ is even,
and $\chi_{M_p-1}$ is contributed if and only if the base is nontrivial,
which is to say, $M_p$ is even.
Now the fibers $\SO(M'-2t)/\SO(M'-2t-2)$
correspond to \MMtrap-shapes of (odd) width $M'-2t$
and $\SO(m')/\SO(m_p)$ and
$\SO(M_p)/\SO(M')$ correspond to the possible leftover
\mtrap-shapes of respective widths $m'$ and $M_p$ 
on the top and bottom of the white \hex-shape 
associated to $V_p$.
\end{proof}

We do not hope to give a simple presentation 
for $\H(F_\SO)$ or even $\H(F_\SO;\F_2)$,
but the mod-$2$ homology
admits a nice description
following the same reasoning as
\Cref{thm:homology-U}.

\begin{proposition}\label{thm:homology-O}
	The quotient map
	\[
	\pi: \prod_{p=1}^r \SO(M_p) 
	\lt
	\SO(M_1) \ox_{\SO(m_1)}
	\cdots \ox_{\SO(m_{r-1})} 
	\SO(M_r)/\SO(m_r) = F_\SO
	\]
	induces a surjection in mod-$2$ homology and an identification
	\[
	H_*(F_\SO;\F_2) \iso 
	H_*\big(\SO(M_1);\F_2\big) \ox_{H_*(\SO(m_1);\F_2)}
	\cdots \ox_{H_*(\SO(m_{r-1});\F_2)} 
	H_*\big(\SO(M_r);\F_2\big) \ox_{H_*(\SO(m_r);\F_2)} \F_2
	\mathrlap.
	\]
Dually, $\pi$ induces an injection in mod-$2$ cohomology.
\end{proposition}

\brmk
One even has $\pi^*$ an injection integrally if none of the $m_p$
except possibly $m_r$ are even, because
the lifts of $q_{4j-1}$ and $\chi_{2s-1}$ for $M_p = 2s$ inject under $\pi^*$;
but $\pi^*u_{m_p}$ is $2$-torsion for~$m_p$ even.
\ermk

\begin{example}\label{rmk:xex2}
	Recall that we promised in \Cref{rmk:xex}
	a counterexample to the most optimistic
	possible claim about the structure
	the cohomology ring of $X = \SO(9) \ox_{\SO(5)} \SO(6)/\SO(2)$.
	Pulling back along the quotient map $\mu_1$ to $V_1 = \SO(9) / \SO(5)$
	we get four generators $v'_8,v'_7,v'_6,v'_5$.
	There is no map to $V_2 = \SO(6)/\SO(2)$,
	but pulling back along the multiplication map $\mu_2$ to $\SO(9)/\SO(3)$
	we get three more generators $v''_5,v''_4,v''_3$.
	These seven classes form a simple system
	of generators for $\H(X;\F_2)$
	whose respective images in $\H\big(\SO(9) \x \SO(6);\F_2\big)$
	under the embedding $\pi^*$ of \Cref{thm:homology-O}
	are
	\[
		v_8 \ox 1,
			\qquad\!
		v_7 \ox 1,
			\qquad\!
		v_6 \ox 1,
			\qquad\!
		v_5 \ox 1,
			\qquad\!
		v_5 \ox 1 + 1 \ox v_5,
			\qquad\!
		v_4 \ox 1 + 1 \ox v_4,
			\qquad\!
		v_3 \ox 1 + 1 \ox v_3\mathrlap.
	\]
In particular, we see $(v''_4)^2 = v'_8$.
Since $v''_4$ lies in the image of $\mu_2^*$ but not in that of $\mu_1^*$
and $v'_8$ lies in the image of $\mu_1^*$ but not of $\mu_2^*$,
there is no way to arrange a ring isomorphism
$\H(F_\SO;\F_2) \iso \H(V_1 \x V_2;\F_2)$.

Integrally, we have a simple system of generators
$q'_{15},\ q'_{11},\ \chi''_5,\ q''_7$ for the free abelian summand,
and in particular $\smash{\ol{H}{}^8(F_\SO)} = 0$.
Now $v''_4 = \Sq^1 v''_3 = \rho \b v''_3$
and $v'_8 = \Sq^1 v'_7 = \rho \b v'_7$
are reductions of integral classes,
and $\rho$ is a ring homomorphism injective on torsion,
so it follows that $(\b v''_3)^2 = \b v'_7$ in $H^8(F_\SO)$.
Thus we cannot have a ring isomorphism $\H(F_\SO) \iso \H(V_1 \x V_2)$.
\end{example}

\begin{remark}\label{rmk:biquotient-cohomology}
	The description of a GZ fiber as a product of biquotients 
	(hence a biquotient itself) from \Cref{thm:biquotient} 
	makes available a different closed-form expression for the cohomology of such spaces.
	Singhof~\cite{singhof1993} observed that
	a biquotient $K \lq G / H$ of compact, connected Lie groups
	can be expressed up to homotopy as the homotopy pullback 
	of the diagram $BK \to BG \from BH$,
	to which Munkholm's general Eilenberg--Moore collapse result~\cite{munkholm1974emss}
	applies to yield 
	\begin{equation}\label{eq:Tor}
	\H(K\lq G/H;k) \iso \Tor^*_{\H(BG;k)}\big(\H(BK;k),\H(BH;k)\big)
	\end{equation}
	as graded modules over
	any principal ideal domain $k$
	of characteristic $\neq 2$ 
	such that the cohomology rings of $BG$, $BK$, $BH$ are polynomial.
	Vitali Kapovitch~\cite{kapovitch2002biquotients},
	building on work of Eschenburg~\cite{eschenburg1992biquotient},
	showed \eqref{eq:Tor} holds as a ring isomorphism for $k = \Q$,
	and the first author recently found,
	using $A_\infty$-algebraic techniques and a product due to Franz
	on a two-sided bar construction,
	that the isomorphism 
	\eqref{eq:Tor} is multiplicative so long as
	$2$ is a unit of $k$~\cite{carlsonfranz}.
	
	Since $\H\U(n)$ is free abelian
	and $\H\SO(n)$ has torsion only of order two,
	\eqref{eq:Tor} determines the cohomology 
	ring of any unitary GZ fiber 
	(or orthogonal GZ fiber for $k =\Zf$).
	For a general biquotient, 
	knowing \eqref{eq:Tor} greatly simplifies life 
	since the Tor in question is algorithmically computable using
	a finitely generated Koszul complex,
	but in the present case, 
	it is easier still
	to simply read the ring structure off the GZ pattern.
\end{remark}

\bs
{\footnotesize

}

	\nd\footnotesize{\textsc{Department of Mathematics,\\ 
					Imperial College London}\\
					\url{j.carlson@imperial.ac.uk}
				}
				
	\smallskip
	
	\nd\footnotesize{\textsc{Department of Mathematics \& Statistics,\\
					McMaster University}\\
					\url{lanej5@math.mcmaster.ca}
				}
\end{document}